\documentclass[a4paper]{amsart}
\usepackage{mathrsfs}
\usepackage{amssymb}
\usepackage{amsrefs}
\newcommand{\Complex}{\mathbb C} 
\newcommand{\Real}{\mathbb R}
\newcommand{\ddbar}{\overline\partial}
\newcommand{\pr}{\partial}
\newcommand{\ol}{\overline}
\newcommand{\Td}{\widetilde}
\newcommand{\norm}[1]{\lefy\Vert#1\right\Vert}
\newcommand{\abs}[1]{\left\vert#1\right\vert}
\newcommand{\set}[1]{\left\{#1\right\}}
\newcommand{\To}{\rightarrow}
\newcommand{\La}{\triangle}
\theoremstyle{plain}
\newtheorem{thm}{Theorem}[section]

\newtheorem{lem}[thm]{Lemma}
\newtheorem{prop}[thm]{Proposition}

\theoremstyle{definition}

\theoremstyle{remark}
\newtheorem{rem}[thm]{Remark}
\numberwithin{equation}{section}
\renewcommand{\labelenumi}{(\alph{enumi})}
\begin{document}
\title[]{The second coefficient of the asymptotic expansion of the weighted Bergman kernel on $\Complex^n$}
\author[]{Chin-Yu Hsiao}
\address{Universit{\"a}t zu K{\"o}ln,  Mathematisches Institut,
    Weyertal 86-90,   50931 K{\"o}ln, Germany}
\email{chsiao@math.uni-koeln.de}  

\begin{abstract}
Let $\phi\in C^\infty(\Complex^n)$ be a given real valued function. We assume that $\pr\ddbar\phi$ is non-degenerate of constant signature $(n_-,n_+)$ on $\Complex^n$. When $q=n_-$, it is well-known that the Bergman kernel for $(0,q)$ forms with respect to the $k$-th weight $e^{-2k\phi}$, $k>0$, admits a full asymptotic expansion in $k$.
In this paper, we compute the trace of the second coefficient of the asymptotic expansion on the diagonal. 
\end{abstract} 

\maketitle \tableofcontents

\section{Introduction and statement of the main result} 

Let $L$ be a holomorphic line bundle over a Hermitian manifold $(M,\Theta)$, where $\Theta$ is a smooth positive $(1,1)$-form on $M$, and let $L^k$ be the $k$-th tensor power of $L$. 
Let $\Box^{(q)}_k$ be the Gaffney extension of the Kodaira Laplacian acting on $(0,q)$ forms with values in $L^k$. The Bergman kernel is the distribution kernel of the orthogonal projection onto ${\rm Ker\,}\Box^{(q)}_k$ in the $L^2$ space. We assume that the curvature of $L$ is non-degenerate of constant signature $(n_-,n_+)$ on $M$ and let $q=n_-$. 
When $M$ is compact,
Catlin~\cite{Cat97} and Zelditch~\cite{Zel98} established the asymptotic expansion of the diagonal of the Bergman kernel for $q=n_-=0$ and Berman-Sj\"{o}strand~\cite{BS05}, Ma-Marinescu~\cite{MM06} established the asymptotic expansion of the Bergman kernel for $q=n_-\geq0$.  When $M$ is complete and $L$ is uniformly positive on $M$ with $\sqrt{-1}R^{K^*_M}$ and $\partial\Theta$ bounded below, where $R^{K^*_M}$ is the curvature of the bundle of $(n,0)$ forms, Ma-Marinescu~\cite{MM08a} obtained the asymptotic expansion of the Bergman kernel for $q=n_-=0$. More generally, if $M$ is any complex manifold and $\Box^{(q)}_k$ has  $O(k^{-n_0})$ small spectral gap on an open set $D\Subset M$ (see Definition 1.5 in~\cite{HM11}, for the precise meaning of $O(k^{-n_0})$ small spectral gap), then it is known by a recent result (see Theorem 1.6 in~\cite{HM11}) that the Bergman kernel admits a full asymptotic expansion in $k$ on $D$. The coefficients of these expansions turned out to be deeply related to various problem in complex geometry (see  e.g.\ \cite{Do:01}, \cite{Fine08}, \cite{Fine10}).

The first four coefficients of the expansion of the Bergman kernel for $q=n_-=0$ on the diagonal were computed by
Lu \cite{Lu00}. The method of Lu is to construct appropriate peak sections as in \cite{Tian}, using 
H{\"o}rmander's $L^2$-method. Ma-Marinescu~\cite{MM10} 
calculated the first three coefficients of the expansion of the kernel of Berezin-Toeplitz quantization on the diagonal by using kernel calculations on $\Complex^n$. The author~\cite{Hsiao12} 
gave a new method to calculated the first three coefficients of the expansion of the kernel of Berezin-Toeplitz quantization on the diagonal by using microlocal analysis. All these results are concern $q=n_-=0$. 

In this paper, we give for the first time a formula of the second coefficient of the expansion of the Bergman kernel for $q=n_->0$ on $\Complex^n$. We calculate the trace of the second coefficient of the expansion of the Bergman kernel for $q=n_->0$ when $M=\Complex^n$ and $L$ is the trivial line bundle $\Complex$ endowed with the metric $\abs{1}^2=e^{-2\phi}$, where 
$\phi\in C^\infty(\Complex^n)$ is a given real valued function with $\pr\ddbar\phi$ is non-degenerate of constant signature $(n_-,n_+)$ on $\Complex^n$. There are two ingredients of our approach: the phase function version of the asymptotic expansion of the Bergman kernel and the method of stationary phase. Even through the calculation is quite complicate, the arguments in this paper are simple. 

After the paper was completed Wen Lu~\cite{Lu12} informed the author that he also obtained the formula for the coefficient $b_1$ by using the method of Ma-Marinescu~\cite{MM07}.  Moreover, Wen Lu~\cite{Lu12} obtained the formula for $b_1$ on general compact complex manifolds.

{\small\emph{
\textbf{Acknowledgements.}
A large part of this paper has been carried out when the author was a postdoc fellow at Chalmers University of Technology during December 2008-December 2010 supported by the Swedish research council. I am grateful to Bo Berndtsson and Robert Berman for several useful conversations and to
the Department of Mathematics, Chalmers University of Technology and the University of G\"{o}teborg for
offering excellent working conditions. I would like to thank Chun-Chung Hsieh at the institute of Mathematics of Academia Sinica, Taiwan, for several interesting discussions on the material of this paper. Furthermore, the author is grateful to George Marinescu for 
comments and useful suggestions on an early draft of the manuscript
}} 

\subsection{Notations} 

Let $\Omega$ be a $C^\infty$ paracompact manifold equipped with a smooth density of integration. We let $T(\Omega)$ and $T^*(\Omega)$ denote the tangent bundle of $\Omega$ and the cotangent bundle of $\Omega$ respectively.
The complexified tangent bundle of $\Omega$ and the complexified cotangent bundle of $\Omega$ will be denoted by $\Complex T(\Omega)$
and $\Complex T^*(\Omega)$ respectively. We write $<\ ,\ >$ to denote the pointwise duality between $T(\Omega)$ and $T^*(\Omega)$.
We extend $<\ ,\ >$ bilinearly to $\Complex T(\Omega)\times\Complex T^*(\Omega)$.
Let $E$ be a $C^\infty$ vector bundle over $\Omega$. The fiber of $E$ at $x\in\Omega$ will be denoted by $E_x$.
Let $Y\subset\subset\Omega$ be an open set. From now on, the spaces of
smooth sections of $E$ over $Y$ and distribution sections of $E$ over $Y$ will be denoted by $C^\infty(Y;\, E)$ and $\mathscr D'(Y;\, E)$ respectively.
Let $\mathscr E'(Y;\, E)$ be the subspace of $\mathscr D'(Y;\, E)$ whose elements have compact support in $Y$. Put $C^\infty_0(Y;\, E)=C^\infty(Y;\, E)\bigcap\mathscr E'(Y;\, E)$.
We let $L^2(Y;\, E)$ denote the $L^2$ space of sections of $E$ over $Y$. 

We shall denote the real coordinates by $x_j$, $j=1,\ldots,2n$, and the complex coordinates by $z=(z_1,\ldots,z_n)$, $z_j=x_{2j-1}+ix_{2j}$, $j=1,\ldots,n$.  
Let $\Lambda^{1,0}T(\Complex^n)$ and $\Lambda^{0,1}T(\Complex^n)$ denote the holomorphic tangent bundle 
and the anti-holomorphic tangent bundle of $\Complex^n$ respectively. We take the Hermitian metric $(\ |\ )$ on $\Complex T(\Complex^n)$ such that 
$(\frac{\pr}{\pr z_j}\ |\ \frac{\pr}{\pr z_k})=(\frac{\pr}{\pr\ol z_j}\ |\ \frac{\pr}{\pr\ol z_k})=\delta_{j,k}$, $j, k=1,\ldots,n$, $\Lambda^{1,0}T(\Complex^n)\bot\Lambda^{0,1}T(\Complex^n)$, where $\frac{\pr}{\pr z_j}=\frac{1}{2}(\frac{\pr}{\pr x_{2j-1}}-i\frac{\pr}{\pr x_{2j}})$, $\frac{\pr}{\pr\ol z_j}=\frac{1}{2}(\frac{\pr}{\pr x_{2j-1}}+i\frac{\pr}{\pr x_{2j}})$, $j=1,\ldots,n$, and $\delta_{j,k}=1$ if $j=k$, $\delta_{j,k}=0$ if $j\neq k$.
For $p, q\geq0$, $p, q\in\mathbb Z$, let $\Lambda^{p,q}T^*(\Complex^n)$ be the bundle of $(p,q)$ forms of $\Complex^n$.
We say that a multiindex $J=(j_1,\ldots,j_q)\in\{1,\ldots,n-1\}^q$ has length $q$ and write $\abs{J}=q$. We say that $J$ is strictly increasing if $1\leqslant  j_1<j_2<\cdots<j_q\leqslant  n-1$. For multiindices $J=(j_1,\ldots,j_q)$, $K=(k_1,\ldots,k_p)$,  we define $dz_K\wedge d\ol z_J:=dz_{k_1}\wedge\cdots\wedge dz_{k_p}\wedge d\ol z_{j_1}\wedge\cdots\wedge d\ol z_{j_q}$.  We take the Hermitian metric $(\ |\ )$ on $\Lambda^{p,q}T^*(\Complex^n)$ so that
$\{dz_K\wedge d\ol z_J: \mbox{$\abs{K}=p$, $\abs{J}=q$, $K$, $J$ are strictly increasing}\}$
is an orthonormal frame for $\Lambda^{p,q}T^*(\Complex^n)$. Let $T\in\mathscr L(\Lambda^{p,q}T^*(\Complex^n), \Lambda^{p,q}T^*(\Complex^n))$. Then the trace of $T$ is given by 
\begin{equation} \label{s1-e-1} 
{\rm Tr\,}T=\sideset{}{'}\sum_{\abs{K}=p,\abs{J}=q}(T(dz_K\wedge d\ol z_J)\ |\ dz_K\wedge d\ol z_J),
\end{equation}
where $\sum'$ means that the summation is performed only over strictly increasing multiindices.
Thus ${\rm Tr\,}T=0$ if 

\begin{equation} \label{s1-e0} 
(T(dz_K\wedge d\ol z_J)\ |\ dz_K\wedge d\ol z_J)=0
\end{equation} 
for all strictly increasing multiinindices $K$, $J$, $\abs{K}=p$, $\abs{J}=q$. 

If $w\in\Lambda^{0,1}T^*_z(\Complex^n)$, let
$w^{\wedge, *}: \Lambda^{0,q+1}T^*_z(\Complex^n)\To \Lambda^{0,q}T^*_z(\Complex^n)$ be the adjoint of left exterior multiplication
$w^\wedge: \Lambda^{0,q}T^*_z(\Complex^n)\To \Lambda^{0,q+1}T^*_z(\Complex^n)$.
That is,
\begin{equation} \label{s1-e1}
(w^\wedge u\ |\ v)=(u\ |\ w^{\wedge, *}v),
\end{equation}
for all $u\in\Lambda^{0,q}T^*_z(\Complex^n)$, $v\in\Lambda^{0,q+1}T^*_z(\Complex^n)$.
Notice that $w^{\wedge, *}$ depends anti-linearly on $w$.

Let $E$, $F$ be $C^\infty$ vector bundles over a smooth manifold $M$. We say that a $k$-dependent function $f(x, y, k)\in C^\infty(M\times M;\, \mathscr L(E_y, F_x))$ is negligible if for every compact set $K\subset M\times M$ and for all $N>0$ and multiindice $\alpha$, $\beta$, there is a constant $c_{N,\alpha,\beta,K}>0$ independent of $k$ such that for $k$ sufficiently large, $\abs{\pr^{\alpha}_x\pr^{\beta}_yf(x, y, k)}\leq c_{\alpha,\beta,K}k^{-N}$, $(x, y)\in K$.  Let $b(x, y, k)\in C^\infty(M\times M;\, \mathscr L(E_y, F_x))$ be a $k$-dependent smooth function. We write 
\[b\sim\sum^\infty_0b_j(x, y)k^{-N_0-j}\]
in $C^\infty(M\times M;\, \mathscr L(E_y, F_x))$, $b_j(x, y)\in C^\infty(M\times M;\, \mathscr L(E_y, F_x))$, $j=0,1,\ldots$, if for all $M_0\in\mathbb N$, every compact set $K\subset M\times M$ and for all multiindice $\alpha$, $\beta$, there is a constant $c_{M_0,\alpha,\beta,K}>0$ independent of $k$ such that for $k$ sufficiently large,
\begin{equation*}
\abs{\pr^{\alpha}_x\pr^{\beta}_y(b-\sum^{M_0}_0b_j(x,y)k^{-N_0-j})}
\leq  c_{M_0,\alpha,\beta,K}k^{-N_0-M_0-1},
\end{equation*} 
$(x, y)\in K$.

\subsection{The asymptotic expansion of the Bergman kernel} 

Let $\phi(z)\in C^\infty(\Complex^n; \Real)$. In this work we assume that 
$\left(\frac{\pr^2\phi}{\pr\ol z_j\pr z_k}\right)^n_{j,k=1}$ is non-degenerate of constant signatute $(n_-, n_+)$. That is, the number of negative eigenvalues of $\left(\frac{\pr^2\phi}{\pr\ol z_j\pr z_k}\right)^n_{j,k=1}$ is $n_-$ and $n_-+n_+=n$. 

We take $dm=2^ndx_1dx_2\cdots dx_{2n}$ as the volume form on $\Complex^n$.
Let $(\ |\ )$ be the inner product on $C^\infty_0(\Complex^n;\, \Lambda^{0,q}T^*(\Complex^n))$
defined by
\begin{equation} \label{s1-e2}
(f\ |\ g)=\int_{\Complex^n}(f(z)\ |\ g(z))(dm),\quad f,g\in C^\infty_0(\Complex^n;\, \Lambda^{0,q}T^*(\Complex^n)).
\end{equation} 
For $k>0$, let $(\ |\ )_k$ be the inner product on $C^\infty_0(\Complex^n;\, \Lambda^{0,q}T^*(\Complex^n))$
defined by
\begin{equation} \label{s1-e3}
(f\ |\ g)_k=\int_{\Complex^n}(f(z)\ |\ g(z))e^{-2k\phi(z)}(dm),\quad f,g\in C^\infty_0(\Complex^n;\, \Lambda^{0,q}T^*(\Complex^n)).
\end{equation}  
Let $L^2_q(\Complex^n)$ and $L^2_{q,k}(\Complex^n)$ be the completions of $C^\infty_0(\Complex^n;\, \Lambda^{0,q}T^*(\Complex^n))$ with respect to $(\ |\ )$ and $(\ |\ )_k$ respectively. We extend the $L^2$ inner products $(\ |\ )$ and $(\ |\ )_k$ to $L^2_q(\Complex^n)$ and $L^2_{q,k}(\Complex^n)$ respectively. 

Let $\ddbar:C^\infty(\Complex^n;\, \Lambda^{0,q}T^*(\Complex^n))\To C^\infty(\Complex^n;\ \Lambda^{0,q+1}T^*(\Complex^n))$
be the part of the exterior differential operator which maps forms of type $(0,q)$ to forms of type $(0,q+1)$. 
We extend 
$\ddbar$ to $L^2_{q,k}(\Complex^n)$ by 
\begin{equation}\label{e:ddbar1}
\ddbar_k:{\rm Dom\,}\ddbar_k\subset L^2_{q,k}(\Complex^n)\To L^2_{q+1,k}(\Complex^n)\,,
\end{equation}
where ${\rm Dom\,}\ddbar_k:=\{u\in L^2_{q,k}(\Complex^n);\, \ddbar u\in L^2_{q+1,k}(\Complex^n)\}$, where 
$\ddbar u$ is defined in the sense of distributions. 
We write 
\begin{equation}\label{e:ddbar_star1}
\ol{\pr}^{\,*}_k:{\rm Dom\,}\ol{\pr}^{\,*}_k\subset L^2_{q+1,k}(\Complex^n)\To L^2_{q,k}(\Complex^n)
\end{equation}
to denote the Hilbert space adjoint of $\ddbar_k$ in the $L^2$ space with respect to $(\ |\ )_k$.
Let $\Box^{(q)}_k$ denote the Gaffney extension of the Kodaira Laplacian given by  
\begin{equation} \label{s1-e4} 
\Box^{(q)}_k=\ddbar_k\ol{\pr}^*_k+\ol{\pr}^*_k\ddbar_k:{\rm Dom\,}\Box^{(q)}_k\subset L^2_{q,k}(\Complex^n)\To L^2_{q,k}(\Complex^n),
\end{equation} 
where
\[{\rm Dom\,}\Box^{(q)}_k=\set{s\in L^2_{q,k}(\Complex^n);\, s\in{\rm Dom\,}\ddbar_k\cap{\rm Dom\,}\ol{\pr}^{\,*}_k,\ \ddbar_ku\in{\rm Dom\,}\ol{\pr}^{\,*}_k,\ \ol{\pr}^{\,*}_ku\in{\rm Dom\,}\ddbar_k}\,.\] 
By a result of Gaffney \cite[Prop.\,3.1.2]{MM07}, $\Box^{(q)}_k$ is a positive self-adjoint operator.
Let 
\begin{equation} \label{s1-e5} 
\Pi^{(q)}_k:L^2_{q,k}(\Complex^n)\To{\rm Ker\,}\Box^{(q)}_k 
\end{equation} 
be the Bergman projection, i.e. the orthogonal projection onto ${\rm Ker\,}\Box^{(q)}_k$ with respect to $(\ |\ )_k$ and
let $\Pi^{(q)}_k(z, w)$ be the distribution kernel of $\Pi^{(q)}_k$ with respect to the volume form $dm$. Since $\Box^{(q)}_k$ is elliptic, it is not difficult to see that 
\[\Pi^{(q)}_{k}(z,w)\in C^\infty\big(\Complex^n\times\Complex^n;\, \mathscr L(\Lambda^{0,q}T^*(\Complex^n),\Lambda^{0,q}T^*(\Complex^n)\big).\]
We write
\begin{equation*} 
\Pi^{(q)}_ku(z)=\int_{\Complex^n}\Pi^{(q)}_k(z, w)u(w)dm(w),
\end{equation*}
$u\in C^\infty_0(\Complex^n;\, \Lambda^{0,q}T^*(\Complex^n))$. 

It is well-known that $\Box^{(q)}_k$ has spectral gap $\geq Ck$, for $k$ large, where $C>0$ is a constant independent of $k$. 
From this observation and Theorem 4.12, Theorem 4.14 in~\cite{HM11}, we deduce the following 

\begin{thm} \label{s1-t1} 
Let $0\leq q\leq n$. If $q\neq n_-$, then $e^{-k\phi(z)}\Pi^{(q)}_k(z, w)e^{k\phi(w)}$ is negligible. If $q=n_-$, then
\begin{equation} \label{s1-e6} 
\Pi^{(q)}_ku(z)=\int_{\Complex^n}e^{ik\psi(z,w)}e^{k(\phi(z)-\phi(w))}b(z, w, k)u(w)dm(w)+Ru, 
\end{equation} 
for $z\in\Complex^n$, $u\in L^2_{q,k}(\Complex^n)$, where 
\[b\sim\sum^\infty_0b_j(z,w)k^{n-j}\]
in $C^\infty(\Complex^n\times\Complex^n; \mathscr L(\Lambda^{0,q}T^*_w(\Complex^n),\Lambda^{0,q}T^*_z(\Complex^n)))$, 
\[b_j\in C^\infty(\Complex^n\times\Complex^n; \mathscr L(\Lambda^{0,q}T^*_w(\Complex^n),\Lambda^{0,q}T^*_z(\Complex^n))),\] 
$j=0,1,\ldots$, 
\begin{equation*}
Ru=\int_{\Complex^n}e^{k(\phi(z)-\phi(w))}r(z, w, k)u(w)dm(w), 
\end{equation*}
$r(z, w, k)$ is negligible and $\psi\in C^\infty(\Complex^n\times\Complex^n)$, $\psi(z, z)=0$, $\psi(z, w)=-\ol\psi(w, z)$, ${\rm Im\,}\psi(z, w)\geq c\abs{z-w}^2$, $c>0$. For $z=w$, we have 
\begin{equation} \label{s1-e7}
\frac{\pr\psi}{\pr\ol z}=i\frac{\pr\phi}{\pr\ol z},\quad\frac{\pr\psi}{\pr z}=-i\frac{\pr\phi}{\pr z},\quad\frac{\pr\psi}{\pr\ol w}=-i\frac{\pr\phi}{\pr\ol z},\quad\frac{\pr\psi}{\pr w}=i\frac{\pr\phi}{\pr z}.
\end{equation}
Moreover, 
\begin{equation} \label{s1-e8} 
\sum^n_{j=1}\Bigr(i\frac{\pr\psi(z,w)}{\pr\ol z_j}+\frac{\pr\phi(z)}{\pr\ol z_j}\Bigr)
\Bigr(-i\frac{\pr\psi(z,w)}{\pr z_j}+\frac{\pr\phi(z)}{\pr z_j}\Bigr)
\end{equation} 
vanishes to infinity order on $z=w$. Furthermore, the Taylor expansion of the phase 
$\psi(z, w)$ is uniquely determined at each point of $z=w$.

In particular, 
\begin{equation} \label{s1-e9} 
\Pi^{(q)}_k(z, z)\sim k^nb_0(z, z)+k^{n-1}b_1(z, z)+k^{n-2}b_2(z, z)+\cdots.
\end{equation}
\end{thm} 

The leading term $b_0(z, z)$ is essentially well-known (see formula (1.24) in Ma-Marinescu \cite{MM06})

\begin{thm} \label{s1-t2} 
Let $q=n_-$. For $p\in\Complex^n$, we assume that 
$\lambda_j(p)$, $j=1,\ldots,n$, are the eigenvalues of $\left(\frac{\pr^2\phi}{\pr\ol z_j\pr z_k}(p)\right)^n_{j,k=1}$ with respect to $(\ |\ )$ and that
$\lambda_j(p)<0$ if $1\leq j\leq n_-$. Let $U_1(p),\ldots,U_n(p)$ be an orthonormal basis of $\Lambda^{1,0}T_p(\Complex^n)$ such that $<\pr\ddbar\phi(p),U_s(p)\wedge\ol U_t(p)>=\delta_{s,t}\lambda_s(p)$, $s,t=1,\ldots,n$. Let $\ol U^*_j(p)$, $j=1,\ldots,n$, denote the orthonormal basis of $\Lambda^{0,1}T^*_p(\Complex^n)$, which is dual to $\ol U_j(p)$, $j=1,\ldots,n$. Then, 
\begin{equation} \label{s1-e10}
b_0(p,p)=\abs{\lambda_1(p)}\abs{\lambda_2(p)}\cdots\abs{\lambda_n(p)}\pi^{-n}\prod_{j=1}^{j=n_-}\ol U^*_j(p)^\wedge \ol U^*_j(p)^{\wedge, *}. 
\end{equation}
\end{thm} 

\subsection{The main result}  

In order to state our result precisely, we have to introduce some notations and definitions. 
Let 
\[F:\Lambda^{1,0}T(\Complex^n)\To \Lambda^{1,0}T(\Complex^n)\otimes\Lambda^{p,q}T^*(\Complex^n)\] 
be a linear operaor, where $p, q\in\mathbb Z$, $p, q\geq0$. We write $F=\left(F_{j,k}\right)^n_{j,k=1}$, $F_{j,k}\in \Lambda^{p,q}T^*(\Complex^n)$, $1\leq j,k\leq n$, 
\begin{equation} \label{s1-e11}
F\frac{\pr}{\pr z_k}=\sum^n_{j=1}\frac{\pr}{\pr z_j}\otimes F_{j,k},
\end{equation}
$k=1,\ldots,n$. We have 
\begin{equation} \label{s1-e12} 
(FU\ |\ V)=\sum^n_{j,k=1}u_k\ol v_jF_{j,k}\in\Lambda^{p,q}T^*(\Complex^n),
\end{equation} 
where $U=\sum^n_{k=1}u_k\frac{\pr}{\pr z_k}, V=\sum^n_{j=1}v_j\frac{\pr}{\pr z_j}\in\Lambda^{1,0}T(\Complex^n)$. 

Let $T:\Lambda^{1,0}T(\Complex^n))\To\Lambda^{1,0}T(\Complex^n)\otimes\Lambda^{r,t}T^*(\Complex^n)$ be another linear operator, where $r, t\in\mathbb Z$, $r, s\geq0$. We write $T=\left(T_{j,k}\right)^n_{j,k=1}$, $T_{j,k}\in\Lambda^{r,t}T^*(\Complex^n)$, $j, k=1,\ldots,n$, as in 
\eqref{s1-e12}. Then 
\[TF:\Lambda^{1,0}T(\Complex^n)\To\Lambda^{1,0}T(\Complex^n)\otimes\Lambda^{p+r,q+t}T^*(\Complex^n)\] 
is the linear operator defined by $TF\frac{\pr}{\pr z_k}=\sum^n_{s,j=1}\frac{\pr}{\pr z_s}\otimes(T_{s,j}\wedge F_{j,k})$, $k=1,\ldots,n$. 

We assume that $F$ is smooth. That is, 
\[F:C^\infty(\Complex^n;\, \Lambda^{1,0}T(\Complex^n))\To C^\infty(\Complex^n;\, \Lambda^{1,0}T(\Complex^n)\otimes\Lambda^{p,q}T^*(\Complex^n)).\]  
Let
\begin{equation} \label{s1-e13} 
\ddbar F: C^\infty(\Complex^n;\, \Lambda^{1,0}T(\Complex^n))\To C^\infty(\Complex^n;\, \Lambda^{1,0}T(\Complex^n)\otimes\Lambda^{p,q+1}T^*(\Complex^n))
\end{equation} 
be the smooth linear operator defined by $\ddbar F\frac{\pr}{\pr z_k}=\sum^n_{j=1}\frac{\pr}{\pr z_j}\otimes(\ddbar F_{j,k})$, $k=1,\ldots,n$. Similarly, 
\begin{equation} \label{s1-e14} 
\pr F: C^\infty(\Complex^n;\, \Lambda^{1,0}T(\Complex^n))\To C^\infty(\Complex^n;\, \Lambda^{1,0}T(\Complex^n)\otimes\Lambda^{p+1,q}T^*(\Complex^n))
\end{equation} 
is the smooth linear operator defined by $\pr F\frac{\pr}{\pr z_k}=\sum^n_{j=1}\frac{\pr}{\pr z_j}\otimes(\pr F_{j,k})$, $k=1,\ldots,n$.

Let 
\begin{equation} \label{s1-e15} 
M_\phi:C^\infty(\Complex^n; \Lambda^{1,0}T(\Complex^n))\To C^\infty(\Complex^n;\, \Lambda^{1,0}T(\Complex^n))
\end{equation}
be the smooth linear map defined by 
\begin{equation} \label{s1-e16} 
(M_\phi U\ |\ V)=<\pr\ddbar\phi, U\wedge\ol V>,
\end{equation} 
$U, V\in C^\infty(\Complex^n;\, \Lambda^{1,0}T(\Complex^n))$. Note that $M_\phi=\left(\frac{\pr^2\phi}{\pr \ol z_j\pr z_k}\right)^n_{j,k=1}$ in the sense of \eqref{s1-e11}. We write $M_\phi^{-1}:C^\infty(\Complex^n; \Lambda^{1,0}T(\Complex^n))\To C^\infty(\Complex^n;\, \Lambda^{1,0}T(\Complex^n))$ to denote the inverse of $M_\phi$. 

We recall that we work with the assumption that $M_\phi$ is non-degenerate of constant signature $(n_-, n_+)$. For $z\in\Complex^n$, we can diagonalize $M_\phi(z)$, i.e. we can find an orthonormal basis $\set{U_j}^n_{j=1}$
 of $\Lambda^{1,0}T(\Complex^n)$ such that $M_\phi(z)U_j(z)=\lambda_j(z)U_j(z)$, $j=1,\ldots,n$, $\lambda_j(z)\in\Real$, $j=1,\ldots,n$. From now on, we assume that 

\begin{align} \label{s1-e17} 
&M_\phi(z)U_j(z)=\lambda_j(z)U_j(z),\ \ U_j(z)\in\Lambda^{1,0}T_z(\Complex^n),\ \,(U_j\ |\ U_j)=1,\ \ j=1,\ldots,n,\nonumber \\
&\lambda_j(z)<0,\ j=1,\ldots,n_-,\ \lambda_j(z)>0,\ j=n_-+1,\ldots,n.
\end{align}
Let $W_+$ be the subbundle of $\Lambda^{1,0}T(\Complex^n)$ spanned by $\set{U_{n_-+1},\ldots,U_n}$ and 
let $W_-$ be the subbundle of $\Lambda^{1,0}T(\Complex^n)$ spanned by $\set{U_{1},\ldots,U_{n_-}}$. 
We take the Hermitian metric $(\ |\ )_{\abs{\phi}}$ on $\Lambda^{1,0}T(\Complex^n)$ 
such that $W_+\bot W_-$, $(U\ |\ V)_{\abs{\phi}}=(M_\phi U\ |\ V)$ if $U, V\in W_+$, $(U\ |\ V)_{\abs{\phi}}=-(M_\phi U\ |\ V)$ if $U, V\in W_-$.
The Hermitian metric $(\ |\ )_{\abs{\phi}}$ on $\Complex T(\Complex^n)$ induces a Hermitian metric on
$\Lambda^{p, q}T^*(\Complex^n)$ also denoted by $(\ |\ )_{\abs{\phi}}$. 

The two form $\pr\ddbar\phi$ induces a connection $D_\phi$ on the bundle $\Lambda^{1,0}T(\Complex^n)$:
\begin{equation} \label{s1-e18}
\begin{split}
&D_\phi=d+\theta:C^\infty(\Complex^n;\, \Lambda^{1,0}(\Complex^n))\To C^\infty(\Complex^n;\, \Complex T^*(\Complex^n)\otimes\Lambda^{1,0} T(\Complex^n)),\\
&D_\phi(\xi)=\sum^n_{j=1}(d\xi_j)\otimes\frac{\pr}{\pr z_j}+\sum_{1\leq j,k\leq n}\theta_{j,k}\xi_k\otimes\frac{\pr}{\pr z_j},
\end{split}
\end{equation}
where $\xi=\sum^n_{j=1}\xi_j\frac{\pr}{\pr z_j}$ and $\theta=h^{-1}\pr h=\left(\theta_{j,k}\right)^n_{j,k=1}$, $h=\left(\frac{\pr^2\phi}{\pr\ol z_j\pr z_k}\right)^n_{j,k=1}$. We call $\theta$ the connection matrix for $\pr\ddbar\phi$. 
The curvature of the connection $D_\phi$ is given by 
\begin{equation} \label{s1-e19} 
\begin{split}
&\Theta_\phi=\ddbar\theta=\left(\ddbar \theta_{j,k}\right)^n_{j,k=1}=\left(\Theta_{j,k}\right)^n_{j,k=1},\\
&\Theta_\phi: C^\infty(\Complex^n;\, \Lambda^{1,0}T(\Complex^n))\To C^\infty(\Complex^n;\, \Lambda^{1,1}T^*(\Complex^n)\otimes\Lambda^{1,0}T(\Complex^n)),\\
&\xi=\sum^n_{j=1}\xi_j\frac{\pr}{\pr z_j}\To\sum^n_{j,k=1}\Theta_{j,k}\xi_k\otimes\frac{\pr}{\pr z_j}.
\end{split}
\end{equation}

For $j=1,\ldots,n$, define 
\begin{align} \label{s1-e20} 
&\delta_j(k)=0\ \ \mbox{if $\{j,k\}\subset\{1,\ldots,q\}$ or $\{j,k\}\subset\{q+1,\ldots,n\}$}\nonumber \\
&\delta_j(k)=1\ \ \mbox{otherwise}.
\end{align} 
Define 
\begin{equation} \label{s1-e21} 
\begin{split}
&Q:\Lambda^{1,0}T(\Complex^n)\To\Lambda^{1,0}T^*(\Complex^n)\otimes\Lambda^{1,0}T(\Complex^n), \\
&<(QU_j\ |\ U_k), U_s>=\Bigr(\frac{\abs{\lambda_j}\delta_k(j)+\abs{\lambda_s}\delta_k(s)}{\abs{\lambda_k}+\abs{\lambda_j}\delta_k(j)+\abs{\lambda_s}\delta_k(s)}
-\delta_k(j)\delta_k(s)\times\\
&\frac{\abs{\lambda_j}^2\abs{\lambda_s}^2}{(\abs{\lambda_j}+\abs{\lambda_k}+\abs{\lambda_s})^2}\bigr(\frac{1}{\abs{\lambda_j}+\abs{\lambda_k}}+\frac{1}{\abs{\lambda_k}+\abs{\lambda_s}}\bigr)^2\Bigr)<(\pr M_\phi U_j\ |\ U_k), U_s>.
\end{split}
\end{equation} 
It is not difficult to see that the definition \eqref{s1-e21} is independent of the choices of eigenvectors $U_1,\ldots,U_n$. 

Define
\begin{equation} \label{s1-e22}
R=\Theta_\phi-(\ddbar M_\phi^{-1})Q:\Lambda^{1,0}T(\Complex^n)\To \Lambda^{1,1}T^*(\Complex^n)\otimes\Lambda^{1,0}T(\Complex^n).
\end{equation} 

Put $e_j=\frac{1}{\sqrt{\abs{\lambda_j}}}U_j$, $j=1,\ldots,n$, where $U_j$, $j=1,\ldots,n$, are as in \eqref{s1-e17}.
The main result of this work is the following 

\begin{thm} \label{t-main} 
Under the assumptions and notations above, let $q=n_-$. For $p\in\Complex^n$ and for $b_1$ in \eqref{s1-e9}, we have 
\begin{equation} \label{s1-e26-I} 
\begin{split}
&{\rm Tr\,}b_1(p, p)\\
&=2^n\abs{\lambda_1(p)}\abs{\lambda_2(p)}\cdots\abs{\lambda_n(p)}\pi^{-n}\times\\
&\Bigr(\sum_{1\leq k\leq q, q+1\leq j\leq n, 1\leq s\leq n }a_{j,k,s}(p)<(\pr M_\phi U_j\ |\ U_k), U_s >^2(p)\\
&+\frac{1}{4}\sum^n_{j,k=1}(1+\delta_j(k)\frac{\abs{\lambda_j}-\abs{\lambda_k}}{\abs{\lambda_j}+\abs{\lambda_k}})< (Re_j\ |\ e_k)_{\abs{\phi}}, \ol e_j\wedge e_k>(p)\\
&-\sum^n_{j,k=1}\delta_j(k)\frac{\abs{\lambda_j}}{\abs{\lambda_j}+\abs{\lambda_k}}{\rm Re\,}((Qe_j\ |\ e_j)\ |\ (\pr M_\phi e_k\ |\ e_k))_{\abs{\phi}}(p)\\
&+\frac{1}{2}\abs{\sum^n_{j=1}(Qe_j\ |\ e_j)}^2_{\abs{\phi}}(p)\Bigr),
\end{split}
\end{equation}
where for $j, k, s=1,\ldots,n$, 
\begin{equation} \label{s1-e24} 
a_{j,k,s}(p)=\frac{\delta_k(j)\delta_k(s)\abs{\lambda_s(p)}}{2(\abs{\lambda_j(p)}+\abs{\lambda_k(p)})^2(\abs{\lambda_j(p)}+\abs{\lambda_k(p)}+\abs{\lambda_s(p)})^2},
\end{equation}
\end{thm}



\begin{rem} \label{s1-r1} 
It is straight forward to see that the right side of \eqref{s1-e26-I} is real (see \eqref{s5-e17}) and is independent of the choices of eigenvectors $U_1,\ldots,U_n$.
\end{rem}

\section{The Taylor expansion of $\psi(z,w)$ at $z=w$} 

From now on, we assume that $q=n_-$.
The goal of this work is to compute ${\rm Tr\,}b_1(p,p)$, for $p\in\Complex^n$. We may assume that $p=0$ and by taking unitary transformation, we can assume that near $0$, 
\[\phi(z)=a_0+\sum^n_{j=1}(a_jz_j+\ol a_j\ol z_j)+\sum^n_{j=1}\lambda_j\abs{z_j}^2+O(\abs{z}^3),\]
where $\lambda_1,\ldots,\lambda_n$ are eigenvalues of $M_\phi$, $a_0\in\Real$, $a_j\in\Complex$, $j=1,\ldots,n$. Put 
$\Td\phi=\phi-a_0-\sum^n_{j=1}(a_jz_j+\ol a_j\ol z_j)$ and let $\Td\Pi^{(q)}_k$ be the Bergman projection with respect to $e^{-2k\Td\phi}$ as in \eqref{s1-e5}. It is easy to see that 
\[\Pi^{(q)}_k=e^{ka_0+k\sum^n_{j=1}a_jz_j}\circ\Td\Pi^{(q)}_k\circ e^{-ka_0-k\sum^n_{j=1}a_jz_j}.\]
Thus, the coefficients of the asymptotic expansion of the kernel of $\Pi^{(q)}_k$ on the diagonal are the same as 
the coefficients of the asymptotic expansion of the kernel of $\Td\Pi^{(q)}_k$. 

From the discussion above, we may assume that 
\begin{equation} \label{s2-e1} 
\phi(z)=\lambda_1\abs{z_1}^2+\lambda_2\abs{z_2}^2+\cdots+\lambda_n\abs{z_n}^2
+O(\abs{z}^3)
\end{equation}
near $z=0$. Suppose that $\lambda_j<0$, $j=1,\ldots,q$, and $\lambda_j>0$, $j=q+1,\ldots,n$. 
In this section, we are going to compute the Taylor expansion of 
$\psi(z,w)$ at $z=w=0$. 
We introduce some notations. For $\alpha=(\alpha_1,\ldots,\alpha_n)$, $\alpha_j\in\mathbb N\bigcup0$, $j=1,\ldots,n$, we put $\alpha'=(\alpha_1,\ldots,\alpha_q)$, $\alpha''=(\alpha_{q+1},\ldots,\alpha_n)$. We write $<\lambda',\alpha'>:=\sum^q_{j=1}\lambda_j\alpha_j$, 
$<\lambda'',\alpha''>:=\sum^n_{j=q+1}\lambda_j\alpha_j$, 
 $<\abs{\lambda'},\alpha'>:=\sum^q_{j=1}\abs{\lambda_j}\alpha_j$ and
$<\abs{\lambda''},\alpha''>:=\sum^n_{j=q+1}\abs{\lambda_j}\alpha_j$.  
The main goal of this section is to prove the following 

\begin{thm} \label{s2-t1} 
Under the assumptions and notations before, we have
\begin{equation} \label{s2-e2} 
\begin{split}
\psi(z,0)&=i\sum^n_{j=1}\abs{\lambda_j}\abs{z_j}^2\\
&+i\sum_{\abs{\alpha}+\abs{\beta}=3,(\alpha'',\beta')\neq0}\frac{<\lambda'',\alpha''>+<\lambda',\beta'>}{<\abs{\lambda''},\alpha''>+<\abs{\lambda'},\beta'>}\frac{\pr^3\phi}{\pr\ol z^\alpha\pr z^\beta}(0)\frac{\ol z^\alpha z^\beta}{\alpha!\beta!}\\
&+\frac{i}{2}\sum_{q+1\leq j,k\leq n,1\leq s\leq q}\frac{1}{\abs{\lambda_{j}}+\abs{\lambda_{k}}+\abs{\lambda_{s}}}
\Bigr(-\abs{\lambda_{j}}-\abs{\lambda_{k}}-\abs{\lambda_{s}} \\
&\quad +\frac{2\abs{\lambda_{j}}\abs{\lambda_{k}}}{\abs{\lambda_{j}}+\abs{\lambda_{s}}}+\frac{2\abs{\lambda_{j}}
\abs{\lambda_{k}}}{\abs{\lambda_{k}}+\abs{\lambda_{s}}}\Bigr)\frac{\pr^3\phi}{\pr z_{j}\pr z_{k}\pr\ol z_{s}}(0)z_jz_k\ol z_s \\
&+\frac{i}{2}\sum_{q+1\leq j\leq n,1\leq t,s\leq q}\frac{1}{\abs{\lambda_{j}}+\abs{\lambda_{t}}+\abs{\lambda_{s}}}
\Bigr(\abs{\lambda_{j}}+\abs{\lambda_{t}}+\abs{\lambda_{s}} \\
&\quad-\frac{2\abs{\lambda_{t}}\abs{\lambda_{s}}}{\abs{\lambda_{j}}+\abs{\lambda_{t}}}-\frac{2\abs{\lambda_{t}}
\abs{\lambda_{s}}}{\abs{\lambda_{j}}+\abs{\lambda_{s}}}\Bigr)\frac{\pr^3\phi}{\pr z_{j}\pr\ol z_{t}\pr\ol z_{s}}(0)z_j\ol z_t\ol z_s \\
&-\frac{i}{3}\sum_{q+1\leq j,k,s\leq n}\frac{\pr^3\phi}{\pr z_j\pr z_k\pr z_s}(0)z_jz_kz_s\\
&+\frac{i}{3}\sum_{1\leq j,k,s\leq q}\frac{\pr^3\phi}{\pr\ol z_j\pr\ol z_k\pr\ol z_s}(0)\ol z_j\ol z_k\ol z_s+O(\abs{z}^4),
\end{split}
\end{equation}
in some neighborhood of $0$. Moreover, we have  
\begin{equation} \label{s2-e3}
\begin{split}
&\frac{\pr^4\psi}{\pr\ol z_{j}\pr z_{j}\pr\ol z_{k}\pr z_{k}}(0,0) \\
&=\sum^q_{t=1}\Bigr(\frac{2i\abs{\lambda_{j}}}{\bigr(\abs{\lambda_{j}}+\abs{\lambda_{k}}\bigr)\bigr(\abs{\lambda_t}+
\abs{\lambda_{j}}\bigr)}\abs{\frac{\pr^3\phi}{\pr\ol z_t\pr z_{j}\pr z_{k}}(0)}^2\\
&\quad +\frac{2i\abs{\lambda_{t}}\abs{\lambda_{k}}\abs{\lambda_{j}}}{\bigr(\abs{\lambda_{j}}+\abs{\lambda_{k}}\bigr)
\bigr(\abs{\lambda_{t}}+\abs{\lambda_{k}}+\abs{\lambda_{j}}\bigr)^2} \\
&\quad\times\bigr(\frac{1}{\abs{\lambda_{j}}+\abs{\lambda_t}}+\frac{1}{\abs{\lambda_{j}}+\abs{\lambda_{k}}}\bigr)\abs{\frac{\pr^3\phi}{\pr\ol z_t\pr z_{j}\pr\ol z_{k}}(0)}^2 \\ 
&\quad+\frac{2i\abs{\lambda_j}}{(\abs{\lambda_j}+\abs{\lambda_t})(\abs{\lambda_j}+\abs{\lambda_k})}\frac{\pr^3\phi}{\pr z_t\pr\ol z_j\pr z_j}(0)\frac{\pr^3\phi}{\pr\ol z_t\pr\ol z_k\pr z_k}(0)\Bigr)\\
&\quad +\sum^n_{t=q+1}\Bigr(\frac{2i\abs{\lambda_{k}}}{\bigr(\abs{\lambda_{j}}+\abs{\lambda_{k}}\bigr)\bigr(\abs{\lambda_t}+
\abs{\lambda_{k}}\bigr)}\abs{\frac{\pr^3\phi}{\pr\ol z_t\pr z_{j}\pr z_{k}}(0)}^2\\
&\quad +\frac{2i\abs{\lambda_{t}}\abs{\lambda_{k}}\abs{\lambda_{j}}}{\bigr(\abs{\lambda_{j}}+\abs{\lambda_{k}}\bigr)
\bigr(\abs{\lambda_{t}}+\abs{\lambda_{k}}+\abs{\lambda_{j}}\bigr)^2} \\
&\quad\times\bigr(\frac{1}{\abs{\lambda_{j}}+\abs{\lambda_{k}}}+\frac{1}{\abs{\lambda_t}+\abs{\lambda_{k}}}\bigr)\abs{\frac{\pr^3\phi}{\pr\ol z_t\pr \ol z_{j}\pr z_{k}}(0)}^2 \\
&\quad+\frac{2i\abs{\lambda_k}}{(\abs{\lambda_k}+\abs{\lambda_t})(\abs{\lambda_j}+\abs{\lambda_k})}\frac{\pr^3\phi}{\pr z_t\pr\ol z_j\pr z_j}(0)\frac{\pr^3\phi}{\pr\ol z_t\pr\ol z_k\pr z_k}(0)\Bigr)\\
&\quad +i\frac{\lambda_{j}+\lambda_{k}}{\abs{\lambda_{j}}+\abs{\lambda_{k}}}\frac{\pr^4\phi}{\pr\ol z_{j}\pr z_{j}\pr \ol z_{k}\pr z_{k}}(0),
\end{split} 
\end{equation} 
where $q+1\leq j\leq n$, $1\leq k\leq q$,
\begin{equation} \label{s2-e4} 
\begin{split}
\frac{\pr^4\psi}{\pr\ol z_{j}\pr z_{j}\pr\ol z_{k}\pr z_{k}}(0,0)
&=\sum^q_{t=1}\Bigr(\frac{2i}{\bigr(\abs{\lambda_t}+\abs{\lambda_{j}}+
\abs{\lambda_{k}}\bigr)^2}\bigr(\abs{\lambda_t}+\abs{\lambda_{j}}+
\abs{\lambda_{k}}\\
&-\frac{\abs{\lambda_{j}}\abs{\lambda_{k}}}{\abs{\lambda_{j}}+\abs{\lambda_t}}-\frac{\abs{\lambda_{j}}\abs{\lambda_{k}}}{\abs{\lambda_{k}}+\abs{\lambda_t}}\bigr)\abs{\frac{\pr^3\phi}{\pr z_t\pr \ol z_{j}\pr\ol z_{k}}(0)}^2\Bigr)\\
&+i\frac{\pr^4\phi}{\pr\ol z_{j}\pr z_{j}\pr\ol z_{k}\pr z_{k}}(0),
\end{split}
\end{equation} 
where $q+1\leq j,k\leq n$, and
\begin{equation} \label{s2-e5} 
\begin{split}
\frac{\pr^4\psi}{\pr\ol z_{j}\pr z_{j}\pr\ol z_{k}\pr z_{k}}(0,0)
&=\sum^n_{t=q+1}\Bigr(\frac{2i}{\bigr(\abs{\lambda_t}+\abs{\lambda_{j}}+
\abs{\lambda_{k}}\bigr)^2}\bigr(\abs{\lambda_t}+\abs{\lambda_{j}}+
\abs{\lambda_{k}}\\
&-\frac{\abs{\lambda_{j}}\abs{\lambda_{k}}}{\abs{\lambda_{k}}+\abs{\lambda_t}}-\frac{\abs{\lambda_{j}}\abs{\lambda_{k}}}{\abs{\lambda_{j}}+\abs{\lambda_t}}\bigr)\abs{\frac{\pr^3\phi}{\pr\ol z_t\pr z_{j}\pr z_{k}}(0)}^2\Bigr)\\
&-i\frac{\pr^4\phi}{\pr\ol z_{j}\pr z_{j}\pr\ol z_{k}\pr z_{k}}(0),
\end{split} 
\end{equation} 
where $1\leq j,k\leq q$. 
\end{thm}

To prove Theorem~\ref{s2-t1}, we first need the following (see \cite{Hsiao08}, for a proof)

\begin{prop} \label{s2-p1} 
Under the assumptions above, we have 
\begin{equation} \label{s2-e6} 
\begin{split}
\psi(z,w)&=i\sum^n_{j=1}\abs{\lambda_j}\abs{z_j-w_j}^2+i\sum^n_{j=1}\lambda_j(\ol z_jw_j-\ol w_jz_j) \\
&\quad+O(\abs{(z, w)}^3)
\end{split}
\end{equation} 
near $z=w=0$. 
\end{prop} 

From \eqref{s2-e1} and \eqref{s2-e6}, we may write 
\begin{equation} \label{s2-e7} 
\begin{split}
\psi(z,w)&=i\sum^n_{j=1}\abs{\lambda_j}\abs{z_j-w_j}^2+i\sum^n_{j=1}\lambda_j(\ol z_jw_j-\ol w_jz_j) \\
&\quad+\psi_3(z, w)+\psi_4(z, w)+\cdots,
\end{split}
\end{equation} 
where $\psi_j(z, w)$ is a homogeneous polynomial of degree $j$ in $(z,w)$, $j=3,4,\ldots$, and 
\begin{equation} \label{s2-e8} 
\phi(z)=\sum^n_{j=1}\lambda_j\abs{z_j}^2+\phi_3(z)+\phi_4(z)+\cdots,
\end{equation} 
where $\phi_j(z)$ is a homogeneous polynomial of degree $j$ in $z$, $j=3,4,\ldots$. Now, using \eqref{s2-e7} and \eqref{s2-e8} in \eqref{s1-e8}, we get 
\begin{equation} \label{s2-e9} 
\begin{split}
&\sum^q_{j=1}\Bigr((2\lambda_j(z_j-w_j)+i\frac{\pr\psi_3}{\pr\ol z_j}+\frac{\pr\phi_3}{\pr\ol z_j}+i\frac{\pr\psi_4}{\pr\ol z_j}+\frac{\pr\phi_4}{\pr\ol z_j})\\
&\times(-i\frac{\pr\psi_3}{\pr z_j}+\frac{\pr\phi_3}{\pr z_j}-i\frac{\pr\psi_4}{\pr z_j}+\frac{\pr\phi_4}{\pr z_j})\Bigr) \\
&+\sum^n_{j=q+1}\Bigr((i\frac{\pr\psi_3}{\pr\ol z_j}+\frac{\pr\phi_3}{\pr\ol z_j}+i\frac{\pr\psi_4}{\pr\ol z_j}+\frac{\pr\phi_4}{\pr\ol z_j})\\
&\quad\times(2\lambda_j(\ol z_j-\ol w_j)-i\frac{\pr\psi_3}{\pr z_j}+\frac{\pr\phi_3}{\pr z_j}-i\frac{\pr\psi_4}{\pr z_j}+\frac{\pr\phi_4}{\pr z_j})\Bigr)
=O(\abs{(z,w)}^5). 
\end{split}
\end{equation}

We regroup the terms in \eqref{s2-e9} according to the order. Then, the order three and order four terms are the following

\begin{equation} \label{s2-e10} 
\begin{split}
 &T\psi_3(z, w)-\sum^q_{j=1}2i\abs{\lambda_j}w_j\frac{\pr\psi_3(z,w)}{\pr z_j}-
 \sum^n_{j=q+1}2i\abs{\lambda_j}\ol w_j\frac{\pr\psi_3(z,w)}{\pr\ol z_j} \\
&\quad=-\sum^q_{j=1}2\lambda_j(z_j-w_j)\frac{\pr\phi_3(z)}{\pr z_j}
-\sum^n_{j=q+1}2\lambda_j(\ol z_j-\ol w_j)\frac{\pr\phi_3(z)}{\pr\ol z_j},
\end{split}
\end{equation}
\begin{equation} \label{s2-e11} 
\begin{split}
 &T\psi_4(z, w)-\sum^q_{j=1}2i\abs{\lambda_j}w_j\frac{\pr\psi_4(z,w)}{\pr z_j}-
 \sum^n_{j=q+1}2i\abs{\lambda_j}\ol w_j\frac{\pr\psi_4(z,w)}{\pr\ol z_j} \\
&\quad=-\sum^q_{j=1}2\lambda_j(z_j-w_j)\frac{\pr\phi_4(z)}{\pr z_j}
-\sum^n_{j=q+1}2\lambda_j(\ol z_j-\ol w_j)\frac{\pr\phi_4(z)}{\pr\ol z_j} \\
&\quad-\sum^n_{j=1}\Bigr(i\frac{\pr\psi_3(z,w)}{\pr\ol z_j}+\frac{\pr\phi_3(z)}{\pr\ol z_j}\Bigr)
\Bigr(-i\frac{\pr\psi_3(z,w)}{\pr z_j}+\frac{\pr\phi_3(z)}{\pr z_j}\Bigr),
\end{split}
\end{equation} 
where
\begin{equation} \label{s2-e12} 
T=\sum^q_{j=1}2i\abs{\lambda_j}z_j\frac{\pr}{\pr z_j}+\sum^n_{j=q+1}2i\abs{\lambda_j}\ol z_j\frac{\pr}{\pr\ol z_j}.
\end{equation} 
Let 
\begin{equation} \label{s2-e13} 
\psi_3=\psi^0_3+\psi^1_3+\psi^2_3+\psi^3_3,
\end{equation}
where $\psi^j_3$ is a homogeneous polynomial of degree $j$ in 
$w$, $j=0,1,2,3$. Now, we write \eqref{s2-e10} according to the degree of homogenity in $w$, we get 
\begin{align}
&T\psi^0_3=-\sum^q_{j=1}2\lambda_jz_j\frac{\pr\phi_3}{\pr z_j}-\sum^n_{j=q+1}2\lambda_j\ol z_j\frac{\pr\phi_3}{\pr\ol z_j}, \label{s2-e14}\\
&T\psi^1_3=\sum^q_{j=1}2\lambda_jw_j(-i\frac{\pr\psi^0_3}{\pr z_j}+\frac{\pr\phi_3}{\pr z_j})+\sum^n_{j=q+1}2\lambda_j\ol w_j(i\frac{\pr\psi^0_3}{\pr \ol z_j}+\frac{\pr\phi_3}{\pr\ol z_j}), \label{s2-e15}\\
&T\psi^2_3=\sum^q_{j=1}2i\abs{\lambda_j}w_j\frac{\pr\psi^1_3}{\pr z_j}+\sum^n_{j=q+1}2i\abs{\lambda_j}\ol w_j\frac{\pr\psi^1_3}{\pr\ol z_j}.\label{s2-e16}
\end{align} 

We need 

\begin{lem} \label{s2-l1}  
We use the same notations as before. Let 
\begin{equation*}
g=\sum_{\alpha,\beta,\gamma,\delta}g_{\alpha,\beta,\gamma,\delta}\ol z^\alpha z^\beta\ol w^\gamma w^\delta, h=\sum_{\alpha,\beta,\gamma,\delta}h_{\alpha,\beta,\gamma,\delta}\ol z^\alpha z^\beta\ol w^\gamma w^\delta,
\end{equation*}
where $\alpha=(\alpha_1,\ldots,\alpha_n)$, $\alpha_j\in\mathbb N\bigcup0$, $j=1,\ldots,n$, and similar for $\beta$, $\gamma$, $\delta$. If $Tg=h$, then 
\begin{align} 
&h_{\alpha,\beta,\gamma,\delta}=0\ \ \mbox{if\ \  $(\alpha'', \beta')=0$},\label{s2-e17} \\
&g_{\alpha,\beta,\gamma,\delta}=
\frac{1}{2i}\frac{1}
{<\abs{\lambda''},\alpha''>+<\abs{\lambda'},\beta'>}
h_{\alpha,\beta,\gamma,\delta}\ \
\mbox{if\ \ $(\alpha'',\beta')\neq0$}. \label{s2-e18}
\end{align}
\end{lem} 

\begin{proof} 
From the definition of $T$ (see \eqref{s2-e12}), we can compute 
\begin{equation*}
T(g_{\alpha,\beta,\gamma,\delta}\ol z^\alpha z^\beta\ol w^\gamma w^\delta)=0
\end{equation*}
if $(\alpha'',\beta')=0$, and 
\begin{equation*}
T(g_{\alpha,\beta,\gamma,\delta}\ol z^\alpha z^\beta\ol w^\gamma w^\delta)=2i(<\abs{\lambda''},\alpha''>+<\abs{\lambda'},\beta'>)
\ol z^\alpha z^\beta\ol w^\gamma w^\delta g_{\alpha,\beta,\gamma,\delta}
\end{equation*}
if $(\alpha'',\beta')\neq0$. From this, the lemma follows.
\end{proof} 

We have $\phi_3=\sum\frac{\pr^3\phi}{\pr\ol z^\alpha\pr z^\beta}(0)\frac{\ol z^\alpha z^\beta}{\alpha!\beta!}$, where $\alpha!=\alpha_1!\alpha_2!\cdots\alpha_n!$ and similar for $\beta!$. From this, we can compute 
\begin{equation} \label{s2-e19} 
\begin{split}
T\psi^0_3&=-\sum^q_{j=1}2\lambda_jz_j\frac{\pr\phi_3}{\pr z_j}-\sum^n_{j=q+1}2\lambda_j\ol z_j\frac{\pr\phi_3}{\pr\ol z_j} \\
&\quad=-2\sum_{\abs{\alpha}+\abs{\beta}=3}(<\lambda'',\alpha''>+<\lambda',\beta'>)\frac{\ol z^\alpha z^\beta}{\alpha!\beta!}\frac{\pr^3\phi}{\pr\ol z^\alpha\pr z^\beta}(0). 
\end{split}
\end{equation} 
From \eqref{s2-e17}, \eqref{s2-e18} and \eqref{s2-e19}, we deduce that
\begin{equation} \label{s2-e20} 
\psi^0_3=i\sum_{\abs{\alpha}+\abs{\beta}=3,(\alpha'',\beta')\neq0}\frac{\pr^3\phi}{\pr\ol z^\alpha\pr z^\beta}(0)\frac{\ol z^\alpha z^\beta}{\alpha!\beta!}\frac{<\lambda'',\alpha''>+<\lambda',\beta'>}{<\abs{\lambda''},\alpha''>+<\abs{\lambda'},\beta'>}+g(z),
\end{equation}
where $\frac{\pr^{\abs{\alpha}+\abs{\beta}}g}{\pr\ol z^\alpha\pr z^\beta}=0$ if $(\alpha'',\beta')\neq0$. 

Now, we compute $\psi^1_3$. From \eqref{s2-e20}, it is straight 
forward to see that 
\begin{equation} \label{s2-e21}
\frac{\pr\psi^0_3}{\pr z_j}=i\sum_{\abs{\alpha}+\abs{\beta}=2}\frac{\ol z^\alpha z^\beta}{\alpha!\beta!}\frac{\pr^3\phi}{\pr\ol z^\alpha\pr z^\beta\pr z_j}(0)\frac{<\lambda'',\alpha''>+<\lambda',\beta'>+\lambda_j}{<\abs{\lambda''},\alpha''>+<\abs{\lambda'},\beta'>+\abs{\lambda_j}},
\end{equation} 
where $1\leq j\leq q$ and 
\begin{equation} \label{s2-e22} 
\frac{\pr\psi^0_3}{\pr\ol z_j}=i\sum_{\abs{\alpha}+\abs{\beta}=2}\frac{\ol z^\alpha z^\beta}{\alpha!\beta!}\frac{\pr^3\phi}{\pr\ol z^\alpha\pr z^\beta\pr\ol z_j}(0)\frac{<\lambda'',\alpha''>+<\lambda',\beta'>+\lambda_j}{<\abs{\lambda''},\alpha''>+<\abs{\lambda'},\beta'>+\abs{\lambda_j}},
\end{equation}  
where $q+1\leq j\leq n$. From \eqref{s2-e21} and \eqref{s2-e22}, we can check that 
\begin{equation*} 
\begin{split}
T\psi^1_3&=\sum^q_{j=1}2\lambda_jw_j(-i\frac{\pr\psi^0_3}{\pr z_j}+\frac{\pr\phi_3}{\pr z_j})
+\sum^n_{j=q+1}2\lambda_j\ol w_j(i\frac{\pr\psi^0_3}{\pr\ol z_j}+\frac{\pr\phi_3}{\pr \ol z_j})\\
&=\sum^q_{j=1}\sum_{\abs{\alpha}+\abs{\beta}=2}\frac{\ol z^\alpha z^\beta w_j}{\alpha!\beta!}\frac{\pr^3\phi}{\pr\ol z^\alpha\pr z^\beta\pr z_j}(0)\frac{4\lambda_j<\abs{\lambda''}, \alpha''>}{<\abs{\lambda''},\alpha''>+<\abs{\lambda'},\beta'>+\abs{\lambda_j}} \\
&\quad+\sum^n_{j=q+1}\sum_{\abs{\alpha}+\abs{\beta}=2}\frac{\ol z^\alpha z^\beta\ol w_j}{\alpha!\beta!}\frac{\pr^3\phi}{\pr\ol z^\alpha\pr z^\beta\pr \ol z_j}(0)\frac{4\lambda_j<\abs{\lambda'}, \beta'>}{<\abs{\lambda''},\alpha''>+<\abs{\lambda'},\beta'>+\abs{\lambda_j}}.
\end{split}
\end{equation*}
From this, we deduce that 
\begin{equation} \label{s2-e23}  
\begin{split}
\psi^1_3&=i\sum^q_{j=1}\sum_{\abs{\alpha}+\abs{\beta}=2,\alpha''\neq0}\Bigr(
\frac{\ol z^\alpha z^\beta w_j}{\alpha!\beta!}\frac{\pr^3\phi}{\pr\ol z^\alpha\pr z^\beta\pr z_j}(0)\\
&\times\frac{2\abs{\lambda_j}<\abs{\lambda''}, \alpha''>}{(<\abs{\lambda''},\alpha''>+<\abs{\lambda'},\beta'>)(<\abs{\lambda''},\alpha''>+<\abs{\lambda'},\beta'>+\abs{\lambda_j})}\Bigr) \\
&-i\sum^n_{j=q+1}\sum_{\abs{\alpha}+\abs{\beta}=2,\beta'\neq0}\Bigr(
\frac{\ol z^\alpha z^\beta\ol w_j}{\alpha!\beta!}\frac{\pr^3\phi}{\pr\ol z^\alpha\pr z^\beta\pr\ol z_j}(0) \\
&\times\frac{2\abs{\lambda_j}<\abs{\lambda'}, \beta'>}{(<\abs{\lambda''},\alpha''>+<\abs{\lambda'},\beta'>)(<\abs{\lambda''},\alpha''>+<\abs{\lambda'},\beta'>+\abs{\lambda_j})}\Bigr)+h(z,w),
\end{split}
\end{equation} 
where $\frac{\pr^{\abs{\alpha}+\abs{\beta}}h}{\pr\ol z^\alpha\pr z^\beta}=0$ 
if $(\alpha'',\beta')\neq0$.

Now, we compute $\psi^2_3$. From \eqref{s2-e23}, we can compute 
\begin{equation} \label{s2-e24} 
\begin{split}
&\sum^q_{j=1}2i\abs{\lambda_j}w_j\frac{\pr\psi^1_3}{\pr z_j} \\
&=\frac{1}{2}\sum_{q+1\leq k\leq n,1\leq j,s\leq q}\Bigr(\ol z_kw_jw_s\frac{\pr^3\phi}{\pr\ol z_k\pr z_j\pr z_s}(0)\frac{-4\abs{\lambda_j}\abs{\lambda_k}\abs{\lambda_s}}{\abs{\lambda_k}+
\abs{\lambda_j}+\abs{\lambda_s}}\\
&\quad\times\bigr(\frac{1}{\abs{\lambda_k}+\abs{\lambda_j}}+
\frac{1}{\abs{\lambda_k}+\abs{\lambda_s}}\bigr)\Bigr)\\
&+\sum_{1\leq j\leq q,q+1\leq k\leq n,\abs{\alpha}+\abs{\beta}=1}\Bigr(\ol z^\alpha z^\beta w_j\ol w_k\frac{\pr^3\phi}{\pr\ol z^\alpha\pr z^\beta\pr z_j\pr\ol z_k}(0)\\
&\times\frac{4\abs{\lambda_j}\abs{\lambda_k}(<\abs{\lambda'},\beta'>+\abs{\lambda_j})}{(<\abs{\lambda''},\alpha''>+<\abs{\lambda'},\beta'>+\abs{\lambda_j})(<\abs{\lambda''},\alpha''>+<\abs{\lambda'},\beta'>+\abs{\lambda_j}+\abs{\lambda_k})}\Bigr)
\end{split}
\end{equation} 
and 
\begin{equation} \label{s2-e25} 
\begin{split}
&\sum^q_{j=q+1}2i\abs{\lambda_j}\ol w_j\frac{\pr\psi^1_3}{\pr\ol z_j} \\
&=\frac{1}{2}\sum_{1\leq k\leq q,q+1\leq j,s\leq n}\Bigr(z_k\ol w_j\ol w_s\frac{\pr^3\phi}{\pr z_k\pr\ol z_j\pr \ol z_s}(0)\frac{4\abs{\lambda_j}\abs{\lambda_k}\abs{\lambda_s}}{\abs{\lambda_k}+\abs{\lambda_j}+\abs{\lambda_s}}\\
&\times\big(\frac{1}{\abs{\lambda_j}+\abs{\lambda_k}}+\frac{1}{\abs{\lambda_k}+\abs{\lambda_s}}\bigr)\Bigr)\\
&-\sum_{1\leq j\leq q,q+1\leq k\leq n,\abs{\alpha}+\abs{\beta}=1}\Bigr(\ol z^\alpha z^\beta w_j\ol w_k\frac{\pr^3\phi}{\pr\ol z^\alpha\pr z^\beta\pr z_j\pr\ol z_k}(0)\\
&\times\frac{4\abs{\lambda_j}\abs{\lambda_k}(<\abs{\lambda''},\alpha''>+\abs{\lambda_k})}{(<\abs{\lambda''},\alpha''>+<\abs{\lambda'},\beta'>+\abs{\lambda_k})(<\abs{\lambda''},\alpha''>+<\abs{\lambda'},\beta'>+\abs{\lambda_k}+\abs{\lambda_j})}\Bigr).
\end{split}
\end{equation}  
We rewrite the last equation in \eqref{s2-e24}:
\begin{equation} \label{s2-e26} 
\begin{split}
&\sum_{1\leq j\leq q,q+1\leq k\leq n,\abs{\alpha}+\abs{\beta}=1}\Bigr(\ol z^\alpha z^\beta w_j\ol w_k\frac{\pr^3\phi}{\pr\ol z^\alpha\pr z^\beta\pr z_j\pr\ol z_k}(0)\\
&\times\frac{4\abs{\lambda_j}\abs{\lambda_k}(<\abs{\lambda'},\beta'>+\abs{\lambda_j})}{(<\abs{\lambda''},\alpha''>+<\abs{\lambda'},\beta'>+\abs{\lambda_j})(<\abs{\lambda''},\alpha''>+<\abs{\lambda'},\beta'>+\abs{\lambda_j}+\abs{\lambda_k})}\Bigr) \\
&=\sum_{q+1\leq s,k\leq n,1\leq j\leq q}\ol z_s\ol w_kw_j\frac{4\abs{\lambda_k}\abs{\lambda_j}^2}{(\abs{\lambda_s}+\abs{\lambda_j})(\abs{\lambda_s}+\abs{\lambda_j}+\abs{\lambda_k})}\frac{\pr^3\phi}{\pr\ol z_s\pr z_j\pr\ol z_k}(0)\\
&\quad+\sum_{q+1\leq k\leq n,1\leq j,s\leq q} z_s\ol w_kw_j\frac{4\abs{\lambda_k}\abs{\lambda_j}}{\abs{\lambda_s}+\abs{\lambda_j}+\abs{\lambda_k}}\frac{\pr^3\phi}{\pr z_s\pr z_j\pr\ol z_k}(0) \\
&+\sum_{1\leq j\leq q,q+1\leq k\leq n,\abs{\alpha}+\abs{\beta}=1,(\alpha'',\beta')=0}\ol z^\alpha z^\beta\ol w_k w_j\frac{4\abs{\lambda_k}\abs{\lambda_j}}{\abs{\lambda_k}+\abs{\lambda_j}}\frac{\pr^3\phi}{\pr\ol z^\alpha\pr z^\beta\pr \ol z_k\pr z_j}(0).
\end{split}
\end{equation} 
Similarly, we rewrite the last equation in \eqref{s2-e25},
\begin{equation} \label{s2-e27} 
\begin{split} 
&-\sum_{1\leq j\leq q,q+1\leq k\leq n,\abs{\alpha}+\abs{\beta}=1}\Bigr(\ol z^\alpha z^\beta w_j\ol w_k\frac{\pr^3\phi}{\pr\ol z^\alpha\pr z^\beta\pr z_j\pr\ol z_k}(0)\\
&\times\frac{4\abs{\lambda_j}\abs{\lambda_k}(<\abs{\lambda''},\alpha''>+\abs{\lambda_k})}{(<\abs{\lambda''},\alpha''>+<\abs{\lambda'},\beta'>+\abs{\lambda_k})(<\abs{\lambda''},\alpha''>+<\abs{\lambda'},\beta'>+\abs{\lambda_k}+\abs{\lambda_j})}\Bigr) \\
&=\sum_{q+1\leq s,k\leq n,1\leq j\leq q}\ol z_s\ol w_kw_j\frac{-4\abs{\lambda_k}\abs{\lambda_j}}{\abs{\lambda_s}+\abs{\lambda_j}+\abs{\lambda_k}}\frac{\pr^3\phi}{\pr\ol z_s\pr z_j\pr\ol z_k}(0)\\
&\quad+\sum_{q+1\leq k\leq n,1\leq j,s\leq q} z_s\ol w_kw_j\frac{-4\abs{\lambda_k}^2\abs{\lambda_j}}{(\abs{\lambda_s}+\abs{\lambda_k})(\abs{\lambda_s}+\abs{\lambda_j}+\abs{\lambda_k})}\frac{\pr^3\phi}{\pr z_s\pr z_j\pr\ol z_k}(0) \\
&+\sum_{1\leq j\leq q,q+1\leq k\leq n,\abs{\alpha}+\abs{\beta}=1,(\alpha'',\beta')=0}\ol z^\alpha z^\beta\ol w_k w_j\frac{-4\abs{\lambda_k}\abs{\lambda_j}}{\abs{\lambda_k}+\abs{\lambda_j}}\frac{\pr^3\phi}{\pr\ol z^\alpha\pr z^\beta
\pr \ol z_k\pr z_j}(0).
\end{split}
\end{equation} 
Combining \eqref{s2-e26}, \eqref{s2-e27} with \eqref{s2-e24} and \eqref{s2-e25}, 
we obtain
\begin{equation} \label{s2-e28} 
\begin{split}
T\psi^2_3&=\sum^q_{j=1}2i\abs{\lambda_j}w_j\frac{\pr\psi^1_3}{\pr z_j}+\sum^n_{j=q+1}2i\abs{\lambda_j}\ol w_j\frac{\pr\psi^1_3}{\pr\ol z_j} \\
&=\frac{1}{2}\sum_{q+1\leq k\leq n,1\leq j,s\leq q}\Bigr(\ol z_kw_jw_s\frac{\pr^3\phi}{\pr\ol z_k\pr z_j\pr z_s}(0)\frac{-4\abs{\lambda_j}\abs{\lambda_k}\abs{\lambda_s}}{\abs{\lambda_k}+
\abs{\lambda_j}+\abs{\lambda_s}}\\
&\times\big(\frac{1}{\abs{\lambda_j}+\abs{\lambda_k}}+\frac{1}{\abs{\lambda_k}+\abs{\lambda_s}}\bigr)\Bigr)\\
&+\frac{1}{2}\sum_{1\leq k\leq q,q+1\leq j,s\leq n}\Bigr(z_k\ol w_j\ol w_s\frac{\pr^3\phi}{\pr z_k\pr\ol z_j\pr \ol z_s}(0)\frac{4\abs{\lambda_j}\abs{\lambda_k}\abs{\lambda_s}}{\abs{\lambda_k}+\abs{\lambda_j}+\abs{\lambda_s}}\\
&\times\big(\frac{1}{\abs{\lambda_j}+\abs{\lambda_k}}+\frac{1}{\abs{\lambda_k}+\abs{\lambda_s}}\bigr)\Bigr)\\
&+\sum_{q+1\leq s,k\leq n,1\leq j\leq q}\ol z_s\ol w_kw_j\frac{-4\abs{\lambda_k}\abs{\lambda_j}\abs{\lambda_s}}{(\abs{\lambda_s}+\abs{\lambda_j})(\abs{\lambda_s}+\abs{\lambda_j}+\abs{\lambda_k})}\frac{\pr^3\phi}{\pr\ol z_s\pr z_j\pr\ol z_k}(0)\\
&+\sum_{q+1\leq k\leq n,1\leq j,s\leq q}z_s\ol w_kw_j\frac{4\abs{\lambda_k}\abs{\lambda_j}\abs{\lambda_s}}{(\abs{\lambda_s}+\abs{\lambda_k})(\abs{\lambda_s}+\abs{\lambda_j}+\abs{\lambda_k})}\frac{\pr^3\phi}{\pr z_s\pr z_j\pr\ol z_k}(0).
\end{split}
\end{equation}
From this and Lemma~\ref{s2-l1}, we get 
\begin{equation} \label{s2-e29}
\begin{split}
\psi^2_3&=i\sum_{q+1\leq k\leq n,1\leq j,s\leq q}\Bigr(\ol z_kw_jw_s\frac{\pr^3\phi}{\pr\ol z_k\pr z_j\pr z_s}(0)\frac{\abs{\lambda_j}\abs{\lambda_s}}{\abs{\lambda_k}+
\abs{\lambda_j}+\abs{\lambda_s}}\\
&\times\big(\frac{1}{\abs{\lambda_j}+\abs{\lambda_k}}+\frac{1}{\abs{\lambda_k}+\abs{\lambda_s}}\bigr)\Bigr)\\
&-i\sum_{1\leq k\leq q,q+1\leq j,s\leq n}\Bigr(z_k\ol w_j\ol w_s\frac{\pr^3\phi}{\pr z_k\pr\ol z_j\pr \ol z_s}(0)\frac{\abs{\lambda_j}\abs{\lambda_s}}{\abs{\lambda_k}+\abs{\lambda_j}+\abs{\lambda_s}}\\
&\times\big(\frac{1}{\abs{\lambda_j}+\abs{\lambda_k}}+\frac{1}{\abs{\lambda_k}+\abs{\lambda_s}}\bigr)\Bigr)\\
&+i\sum_{q+1\leq s,k\leq n,1\leq j\leq q}\ol z_s\ol w_kw_j\frac{2\abs{\lambda_j}\abs{\lambda_k}}{(\abs{\lambda_s}+\abs{\lambda_j})(\abs{\lambda_s}+\abs{\lambda_j}+\abs{\lambda_k})}\frac{\pr^3\phi}{\pr\ol z_s\pr z_j\pr\ol z_k}(0)\\
&-i\sum_{q+1\leq k\leq n,1\leq j,s\leq q}z_s\ol w_kw_j\frac{2\abs{\lambda_j}\abs{\lambda_k}}{(\abs{\lambda_s}+\abs{\lambda_k})(\abs{\lambda_s}+\abs{\lambda_j}+\abs{\lambda_k})}\frac{\pr^3\phi}{\pr z_s\pr z_j\pr\ol z_k}(0)\\
&+r(z,w),
\end{split}
\end{equation} 
where $\frac{\pr^{\abs{\alpha}+\abs{\beta}}r}{\pr\ol z^\alpha\pr z^\beta}=0$ if 
$(\alpha'', \beta')\neq0$. Summing up, we get the following 
\begin{prop} \label{ps2-1.1} 
Under the assumptions and notations before, we have 
\begin{equation} \label{s2-e30} 
\begin{split}
&\psi(z,w)=i\sum^n_{j=1}\abs{\lambda_j}\abs{z_j-w_j}^2+i\sum^n_{j=1}\lambda_j(\ol z_jw_j-\ol w_jz_j)\\ 
&+i\sum_{\abs{\alpha}+\abs{\beta}=3,(\alpha'',\beta')\neq0}\frac{<\lambda'',\alpha''>+<\lambda',\beta'>}{<\abs{\lambda''},\alpha''>+<\abs{\lambda'},\beta'>}\frac{\pr^3\phi}{\pr\ol z^\alpha\pr z^\beta}(0)\frac{\ol z^\alpha z^\beta}{\alpha!\beta!}\\
&+i\sum^q_{j=1}\sum_{\abs{\alpha}+\abs{\beta}=2,\alpha''\neq0}\Bigr(
\frac{\ol z^\alpha z^\beta w_j}{\alpha!\beta!}\frac{\pr^3\phi}{\pr\ol z^\alpha\pr z^\beta\pr z_j}(0)\\
&\times\frac{2\abs{\lambda_j}<\abs{\lambda''}, \alpha''>}{(<\abs{\lambda''},\alpha''>+<\abs{\lambda'},\beta'>)(<\abs{\lambda''},\alpha''>+<\abs{\lambda'},\beta'>+\abs{\lambda_j})}\Bigr) \\
&-i\sum^n_{j=q+1}\sum_{\abs{\alpha}+\abs{\beta}=2,\beta'\neq0}\Bigr(
\frac{\ol z^\alpha z^\beta\ol w_j}{\alpha!\beta!}\frac{\pr^3\phi}{\pr\ol z^\alpha\pr z^\beta\pr\ol z_j}(0) \\
&\times\frac{2\abs{\lambda_j}<\abs{\lambda'}, \beta'>}{(<\abs{\lambda''},\alpha''>+<\abs{\lambda'},\beta'>)(<\abs{\lambda''},\alpha''>+<\abs{\lambda'},\beta'>+\abs{\lambda_j})}\Bigr)\\
&+i\sum_{q+1\leq k\leq n,1\leq j,s\leq q}\Bigr(\ol z_kw_jw_s\frac{\pr^3\phi}{\pr\ol z_k\pr z_j\pr z_s}(0)\frac{\abs{\lambda_j}\abs{\lambda_s}}{\abs{\lambda_k}+
\abs{\lambda_j}+\abs{\lambda_s}}\\
&\times\big(\frac{1}{\abs{\lambda_j}+\abs{\lambda_k}}+\frac{1}{\abs{\lambda_k}+\abs{\lambda_s}}\bigr)\Bigr)\\
&-i\sum_{1\leq k\leq q,q+1\leq j,s\leq n}\Bigr(z_k\ol w_j\ol w_s\frac{\pr^3\phi}{\pr z_k\pr\ol z_j\pr \ol z_s}(0)\frac{\abs{\lambda_j}\abs{\lambda_s}}{\abs{\lambda_k}+\abs{\lambda_j}+\abs{\lambda_s}}\\
&\times\big(\frac{1}{\abs{\lambda_j}+\abs{\lambda_k}}+\frac{1}{\abs{\lambda_k}+\abs{\lambda_s}}\bigr)\Bigr)\\
&+i\sum_{q+1\leq s,k\leq n,1\leq j\leq q}\ol z_s\ol w_kw_j\frac{2\abs{\lambda_j}\abs{\lambda_k}}{(\abs{\lambda_s}+\abs{\lambda_j})(\abs{\lambda_s}+\abs{\lambda_j}+\abs{\lambda_k})}\frac{\pr^3\phi}{\pr\ol z_s\pr z_j\pr\ol z_k}(0)\\
&-i\sum_{q+1\leq k\leq n,1\leq j,s\leq q}z_s\ol w_kw_j\frac{2\abs{\lambda_j}\abs{\lambda_k}}{(\abs{\lambda_s}+\abs{\lambda_k})(\abs{\lambda_s}+\abs{\lambda_j}+\abs{\lambda_k})}\frac{\pr^3\phi}{\pr z_s\pr z_j\pr\ol z_k}(0)\\
&+R(z,w)+O(\abs{(z,w)}^4),
\end{split}
\end{equation} 
where $R(z,w)=O(\abs{(z,w)}^3)$ and
$\frac{\pr^{\abs{\alpha}+\abs{\beta}}R}{\pr\ol z^\alpha\pr z^\beta}=0$ if $(\alpha'',\beta')\neq0$.
\end{prop}

Now, we are readay to compute $\frac{\pr^3\psi}{\pr\ol z^\alpha\pr z^\beta}(0,0)$, 
where $\abs{\alpha}+\abs{\beta}=3$ and $(\alpha'',\beta')=0$. We compute 
$\frac{\pr^3\psi}{\pr\ol z_{s_0}\pr z_{j_0}\pr z_{k_0}}(0,0)$, $1\leq s_0\leq q$, 
$q+1\leq j_0, k_0\leq n$. From $\psi(z, w)=-\ol\psi(w, z)$, we have 
\begin{equation} \label{s2-e31}
\frac{\pr^3\psi(z, w)}{\pr\ol z_{s_0}\pr z_{j_0}\pr z_{k_0}}=-\ol{\frac{\pr^3\psi(z, w)}{\pr w_{s_0}\pr\ol w_{j_0}\pr\ol w_{k_0}}}.
\end{equation} 
To compute $\frac{\pr^3\psi}{\pr\ol z_{s_0}\pr z_{j_0}\pr z_{k_0}}(0,0)$, it is equivalent to compute $\frac{\pr^3\psi}{\pr w_{s_0}\pr\ol w_{j_0}\pr\ol w_{k_0}}(0,0)$. From \eqref{s1-e7}, we know that 
\begin{equation} \label{s2-e32}
\frac{\pr\psi(z, z)}{\pr\ol w_{j_0}}=-i\frac{\pr\phi(z)}{\pr\ol z_{j_0}}.
\end{equation} 
Differentiate \eqref{s2-e32} with respect to $\ol z_{k_0}$, we get 
\begin{equation} \label{s2-e33} 
\frac{\pr^2\psi(z, z)}{\pr\ol w_{j_0}\pr\ol z_{k_0}}+\frac{\pr^2\psi(z, z)}{\pr\ol w_{j_0}\pr\ol w_{k_0}}=-i\frac{\pr^2\phi(z)}{\pr\ol z_{k_0}\pr\ol z_{j_0}}.
\end{equation}
Again, differentiate \eqref{s2-e33} with respect to $z_{s_0}$, we get 
\begin{equation} \label{s2-e34} 
\begin{split}
&\frac{\pr^3\psi}{\pr\ol w_{j_0}\pr\ol z_{k_0}z_{s_0}}(0,0)+\frac{\pr^3\psi}{\pr\ol w_{j_0}\pr\ol z_{k_0}\pr w_{s_0}}(0,0) \\
&+\frac{\pr^3\psi}{\pr\ol w_{j_0}\pr\ol w_{k_0}\pr z_{s_0}}(0,0)+\frac{\pr^3\psi}{\pr\ol w_{j_0}\pr\ol w_{k_0}\pr w_{s_0}}(0,0)
=-i\frac{\pr^3\phi}{\pr\ol z_{k_0}\pr\ol z_{j_0}\pr z_{s_0}}(0). 
\end{split} 
\end{equation}
From \eqref{s2-e34}, we have 
\begin{equation} \label{s2-e35} 
\begin{split}
\frac{\pr^3\psi}{\pr\ol w_{j_0}\pr\ol w_{k_0}\pr w_{s_0}}(0,0)
&=-i\frac{\pr^3\phi}{\pr\ol z_{k_0}\pr\ol z_{j_0}\pr z_{s_0}}(0)-\frac{\pr^3\psi}{\pr\ol w_{j_0}\pr\ol z_{k_0}z_{s_0}}(0,0) \\
&\quad -\frac{\pr^3\psi}{\pr\ol w_{j_0}\pr\ol z_{k_0}\pr w_{s_0}}(0,0)-
\frac{\pr^3\psi}{\pr\ol w_{j_0}\pr\ol w_{k_0}\pr z_{s_0}}(0,0). 
\end{split}
\end{equation} 
In view of \eqref{s2-e30}, we see that 
\begin{align*}
&\frac{\pr^3\psi}{\pr\ol w_{j_0}\pr\ol z_{k_0}z_{s_0}}(0,0)=\frac{-2i\abs{\lambda_{j_0}}\abs{\lambda_{s_0}}}{\bigr(\abs{\lambda_{j_0}}+\abs{\lambda_{k_0}}+\abs{\lambda_{s_0}}\bigr)
\bigr(\abs{\lambda_{k_0}}+\abs{\lambda_{s_0}}\bigr)}\frac{\pr^3\phi}{\pr\ol z_{k_0}\pr\ol z_{j_0}\pr z_{s_0}}(0)\\
&\frac{\pr^3\psi}{\pr\ol w_{j_0}\pr\ol z_{k_0}\pr w_{s_0}}(0,0)=
\frac{2i\abs{\lambda_{j_0}}\abs{\lambda_{s_0}}}{\bigr(\abs{\lambda_{j_0}}+\abs{\lambda_{k_0}}+\abs{\lambda_{s_0}}\bigr)
\bigr(\abs{\lambda_{k_0}}+\abs{\lambda_{s_0}}\bigr)}\frac{\pr^3\phi}{\pr\ol z_{k_0}\pr\ol z_{j_0}\pr z_{s_0}}(0)\\
&\frac{\pr^3\psi}{\pr\ol w_{j_0}\pr\ol w_{k_0}\pr z_{s_0}}(0,0)=
\frac{-2i\abs{\lambda_{j_0}}\abs{\lambda_{k_0}}}{\abs{\lambda_{k_0}}+\abs{\lambda_{s_0}}+\abs{\lambda_{j_0}}}\\
&\times\bigr(\frac{1}{\abs{\lambda_{j_0}}+\abs{\lambda_{s_0}}}+\frac{1}{\abs{\lambda_{k_0}}+\abs{\lambda_{s_0}}}\bigr)\frac{\pr^3\phi}{\pr\ol z_{k_0}\pr\ol z_{j_0}\pr z_{s_0}}(0).
\end{align*} 
Combining this with \eqref{s2-e35}, we have
\begin{equation} \label{s2-e36}
\begin{split}
&\frac{\pr^3\psi}{\pr\ol w_{j_0}\pr\ol w_{k_0}\pr w_{s_0}}(0,0)=\Bigr(-i-\frac{-2i\abs{\lambda_{j_0}}\abs{\lambda_{s_0}}}{\bigr(\abs{\lambda_{j_0}}+\abs{\lambda_{k_0}}+\abs{\lambda_{s_0}}\bigr)
\bigr(\abs{\lambda_{k_0}}+\abs{\lambda_{s_0}}\bigr)}\\
&\quad -\frac{2i\abs{\lambda_{j_0}}\abs{\lambda_{s_0}}}{\bigr(\abs{\lambda_{j_0}}+\abs{\lambda_{k_0}}+\abs{\lambda_{s_0}}\bigr)
\bigr(\abs{\lambda_{k_0}}+\abs{\lambda_{s_0}}\bigr)}\\
&\quad -\frac{-2i\abs{\lambda_{j_0}}\abs{\lambda_{k_0}}}{\abs{\lambda_{k_0}}+\abs{\lambda_{s_0}}+\abs{\lambda_{j_0}}}\bigr(\frac{1}{\abs{\lambda_{j_0}}+\abs{\lambda_{s_0}}}+\frac{1}{\abs{\lambda_{k_0}}+\abs{\lambda_{s_0}}}\bigr)\Bigr) \\
&\quad\times\frac{\pr^3\phi}{\pr\ol z_{k_0}\pr\ol z_{j_0}\pr z_{s_0}}(0)  \\
&=\frac{i}{\abs{\lambda_{j_0}}+\abs{\lambda_{k_0}}+\abs{\lambda_{s_0}}}
\Bigr(-\abs{\lambda_{j_0}}-\abs{\lambda_{k_0}}-\abs{\lambda_{s_0}} \\
&\quad +\frac{2\abs{\lambda_{j_0}}\abs{\lambda_{k_0}}}{\abs{\lambda_{j_0}}+\abs{\lambda_{s_0}}}+\frac{2\abs{\lambda_{j_0}}
\abs{\lambda_{k_0}}}{\abs{\lambda_{k_0}}+\abs{\lambda_{s_0}}}\Bigr)\frac{\pr^3\phi}{\pr\ol z_{k_0}\pr\ol z_{j_0}\pr z_{s_0}}(0).
\end{split}
\end{equation}
From \eqref{s2-e31}, we obtain 
\begin{equation} \label{s2-e37} 
\begin{split}
\frac{\pr^3\psi_3}{\pr z_{j_0}\pr z_{k_0}\pr\ol z_{s_0}}&=\frac{i}{\abs{\lambda_{j_0}}+\abs{\lambda_{k_0}}+\abs{\lambda_{s_0}}}
\Bigr(-\abs{\lambda_{j_0}}-\abs{\lambda_{k_0}}-\abs{\lambda_{s_0}} \\
&\quad +\frac{2\abs{\lambda_{j_0}}\abs{\lambda_{k_0}}}{\abs{\lambda_{j_0}}+\abs{\lambda_{s_0}}}+\frac{2\abs{\lambda_{j_0}}
\abs{\lambda_{k_0}}}{\abs{\lambda_{k_0}}+\abs{\lambda_{s_0}}}\Bigr)\frac{\pr^3\phi}{\pr z_{j_0}\pr z_{k_0}\pr\ol z_{s_0}}(0).
\end{split}
\end{equation} 
We can repeat the method above several times to determine all the terms $\frac{\pr^3\psi_3}{\pr\ol z^\alpha\pr z^\beta}$, $(\alpha'',\beta')=0$. The computation is straight forward. We omit the process. We state our result 

\begin{prop} \label{s2-p2} 
Under the assumptions and notations before, we have
\begin{equation} \label{s2-e38} 
\begin{split}
\psi(z,0)&=i\sum^n_{j=1}\abs{\lambda_j}\abs{z_j}^2\\
&+i\sum_{\abs{\alpha}+\abs{\beta}=3,(\alpha'',\beta')\neq0}\frac{<\lambda'',\alpha''>+<\lambda',\beta'>}{<\abs{\lambda''},\alpha''>+<\abs{\lambda'},\beta'>}\frac{\pr^3\phi}{\pr\ol z^\alpha\pr z^\beta}(0)\frac{\ol z^\alpha z^\beta}{\alpha!\beta!}\\
&+\frac{i}{2}\sum_{q+1\leq j,k\leq n,1\leq s\leq q}\frac{1}{\abs{\lambda_{j}}+\abs{\lambda_{k}}+\abs{\lambda_{s}}}
\Bigr(-\abs{\lambda_{j}}-\abs{\lambda_{k}}-\abs{\lambda_{s}} \\
&\quad +\frac{2\abs{\lambda_{j}}\abs{\lambda_{k}}}{\abs{\lambda_{j}}+\abs{\lambda_{s}}}+\frac{2\abs{\lambda_{j}}
\abs{\lambda_{k}}}{\abs{\lambda_{k}}+\abs{\lambda_{s}}}\Bigr)\frac{\pr^3\phi}{\pr z_{j}\pr z_{k}\pr\ol z_{s}}(0)z_jz_k\ol z_s \\
&+\frac{i}{2}\sum_{q+1\leq j\leq n,1\leq k,s\leq q}\frac{1}{\abs{\lambda_{j}}+\abs{\lambda_{k}}+\abs{\lambda_{s}}}
\Bigr(\abs{\lambda_{j}}+\abs{\lambda_{k}}+\abs{\lambda_{s}} \\
&\quad-\frac{2\abs{\lambda_{k}}\abs{\lambda_{s}}}{\abs{\lambda_{j}}+\abs{\lambda_{k}}}-\frac{2\abs{\lambda_{k}}
\abs{\lambda_{s}}}{\abs{\lambda_{j}}+\abs{\lambda_{s}}}\Bigr)\frac{\pr^3\phi}{\pr z_{j}\pr\ol z_{k}\pr\ol z_{s}}(0)z_j\ol z_k\ol z_s \\
&-\frac{i}{3}\sum_{q+1\leq j,k,s\leq n}\frac{\pr^3\phi}{\pr z_j\pr z_k\pr z_s}(0)z_jz_kz_s\\
&+\frac{i}{3}\sum_{1\leq j,k,s\leq q}\frac{\pr^3\phi}{\pr\ol z_j\pr\ol z_k\pr\ol z_s}(0)\ol z_j\ol z_k\ol z_s+O(\abs{z}^4),
\end{split}
\end{equation}
in some neighborhood of $0$.
\end{prop} 

Now, to complete the proof of Theorem~\ref{s2-t1}, we only need to compute the terms $\frac{\pr^4\psi}{\pr\ol z_j\pr z_j\pr\ol z_k\pr z_k}(0,0)$, $1\leq j,k\leq n$. We go through the steps in the computation of the case: $q+1\leq j\leq n$, $1\leq k\leq q$. The computations of the cases: $1\leq j,k\leq q$ and 
$q+1\leq j,k\leq n$, are similar but simpler and is therefore omitted. Now, we compute $\frac{\pr^4\psi}{\pr\ol z_{j_0}\pr z_{j_0}\pr\ol z_{k_0}\pr z_{k_0}}(0,0)$, $q+1\leq j_0\leq n$, $1\leq k_0\leq q$. Take $w=0$ in \eqref{s2-e11}, we obtain 
\begin{equation} \label{s2-e39} 
\begin{split} 
T\psi_4(z,0)&=-\sum^n_{j=1}\bigr(i\frac{\pr\psi_3}{\pr\ol z_j}(z, 0)+\frac{\pr\phi_3}{\pr\ol z_j}(z)\bigr)\bigr(-i\frac{\pr\psi_3}{\pr z_j}(z, 0)+\frac{\pr\phi_3}{\pr z_j}(z)\bigr)\\
&\quad-\sum^q_{j=1}2\lambda_jz_j\frac{\pr\phi_4}{\pr z_j}(z)
-\sum^n_{j=q+1}2\lambda_j\ol z_j\frac{\pr\phi_4}{\pr\ol z_j}(z) .
\end{split}
\end{equation} 
We write 
\begin{subequations}\label{s2-e40} 
\begin{gather} 
\begin{split}
i\frac{\pr\psi_3}{\pr\ol z_t}(z, 0)+\frac{\pr\phi_3}{\pr\ol z_t}& =\alpha^1_t\ol z_{j_0}\ol z_{k_0}+\alpha^2_t\ol z_{j_0}z_{k_0}+\alpha^3_tz_{j_0}\ol z_{k_0}\\
&\quad+\alpha^4_tz_{j_0}z_{k_0}+\alpha^5_t\ol z_{j_0}z_{j_0}+\alpha^6_t\ol z_{k_0}z_{k_0}+h(z), 
\end{split} \\
\begin{split}
-i\frac{\pr\psi_3}{\pr z_t}(z, 0)+\frac{\pr\phi_3}{\pr z_t}& =
\beta^1_t\ol z_{j_0}\ol z_{k_0}+\beta^2_t\ol z_{j_0}z_{k_0}+\beta^3_tz_{j_0}\ol z_{k_0}\\
&\quad+\beta^4_tz_{j_0}z_{k_0}+\beta^5_t\ol z_{j_0}z_{j_0}+\beta^6_t\ol z_{k_0}z_{k_0}+r(z),
\end{split}
\end{gather} 
\end{subequations}
where $\frac{\pr^2h}{\pr z_{j_0}\pr z_{k_0}}=\frac{\pr^2h}{\pr\ol z_{j_0}\pr z_{k_0}}=\frac{\pr^2h}{\pr z_{j_0}\pr\ol z_{k_0}}=\frac{\pr^2h}
{\pr\ol z_{j_0}\pr \ol z_{k_0}}=\frac{\pr^2h}{\pr z_{j_0}\pr\ol z_{j_0}}=\frac{\pr^2h}{\pr z_{k_0}\pr\ol z_{k_0}}=0$ and  $\frac{\pr^2r}{\pr z_{j_0}\pr z_{k_0}}=\frac{\pr^2r}{\pr\ol z_{j_0}\pr z_{k_0}}=\frac{\pr^2r}{\pr z_{j_0}\pr\ol z_{k_0}}=\frac{\pr^2r}{\pr\ol z_{j_0}\pr \ol z_{k_0}}=\frac{\pr^2r}{\pr z_{j_0}\pr\ol z_{j_0}}=\frac{\pr^2r}{\pr z_{k_0}\pr \ol z_{k_0}}=0$.
From Proposition~\ref{s2-p2}, it is straight forward to see that if $1\leq t\leq q$, then
\begin{subequations}\label{s2-e41}
\begin{gather} 
\alpha^1_t=0 , \\
\alpha^2_t=\frac{2\abs{\lambda_{k_0}}}{\abs{\lambda_{j_{0}}}+\abs{\lambda_{k_0}}}\frac{\pr^3\phi_3}{\pr\ol z_t\pr\ol z_{j_0}\pr z_{k_0}}, \\ 
\alpha^3_t=\frac{2\abs{\lambda_t}\abs{\lambda_{k_0}}}{\abs{\lambda_{t}}+\abs{\lambda_{j_0}}+\abs{\lambda_{k_0}}}\Bigr( 
\frac{1}{\abs{\lambda_{j_0}}+\abs{\lambda_{t}}}+\frac{1}{\abs{\lambda_{j_0}}+\abs{\lambda_{k_0}}}\Bigr)\frac{\pr^3\phi_3}{\pr\ol z_t\pr z_{j_0}\pr\ol z_{k_0}}, \\
\alpha^4_t=2\frac{\pr^3\phi_3}{\pr\ol z_t\pr z_{j_0}\pr z_{k_0}},\\
\alpha^5_t=0, \\
\alpha^6_t=2\frac{\pr^3\phi}{\pr\ol z_t\pr\ol z_{k_0}\pr z_{k_0}}(0),
\end{gather}
\end{subequations}   
and
\begin{subequations}\label{s2-e42} 
\begin{gather}
\beta^1_t=\frac{2\abs{\lambda_{j_0}}}{\abs{\lambda_t}+\abs{\lambda_{j_0}}}
\frac{\pr^3\phi_3}{\pr z_t\pr\ol z_{j_0}\pr\ol z_{k_0}}, \\
\beta^2_t=\frac{2\abs{\lambda_{j_0}}}{\abs{\lambda_t}+\abs{\lambda_{j_0}}
+\abs{\lambda_{k_0}}}\frac{\pr^3\phi_3}{\pr z_t\pr\ol z_{j_0}\pr z_{k_0}},\\ 
\beta^3_t=0 , \\
\beta^4_t=0 ,\\
\beta^5_t=\frac{2\abs{\lambda_{j_0}}}{\abs{\lambda_{j_0}}+\abs{\lambda_t}}\frac{\pr^3\phi}{\pr z_t\pr\ol z_{j_0}\pr z_{j_0}}(0),\\
\beta^6_t=0. 
\end{gather}
\end{subequations} 
If $q+1\leq t\leq n$, we have 
\begin{subequations}\label{s2-e43}
\begin{gather} 
\alpha^1_t=0 ,\\
\alpha^2_t=\frac{2\abs{\lambda_{k_0}}}{\abs{\lambda_t}+\abs{\lambda_{j_{0}}}+
\abs{\lambda_{k_0}}}\frac{\pr^3\phi_3}{\pr\ol z_t\pr\ol z_{j_0}\pr z_{k_0}}, \\ 
\alpha^3_t=0, \\
\alpha^4_t=\frac{2\abs{\lambda_{k_0}}}{\abs{\lambda_t}+\abs{\lambda_{k_0}}}
\frac{\pr^3\phi_3}{\pr\ol z_t\pr z_{j_0}\pr z_{k_0}}, \\
\alpha^5_t=0,\\
\alpha^6_t=\frac{2\abs{\lambda_{k_0}}}{\abs{\lambda_t}+\abs{\lambda_{k_0}}}\frac{\pr^3\phi}{\pr\ol z_t\pr\ol z_{k_0}\pr z_{k_0}}(0),
\end{gather}
\end{subequations}  
and
\begin{subequations}\label{s2-e44} 
\begin{gather}
\beta^1_t=2\frac{\pr^3\phi_3}{\pr z_t\pr\ol z_{j_0}\pr\ol z_{k_0}}, \\
\beta^2_t=\frac{2\abs{\lambda_{j_0}}}{\abs{\lambda_{j_0}}+\abs{\lambda_{k_0}}}
\frac{\pr^3\phi_3}{\pr z_t\pr\ol z_{j_0}\pr z_{k_0}},\\ 
\beta^3_t=\frac{2\abs{\lambda_{j_0}}\abs{\lambda_t}}{\abs{\lambda_{j_0}}
+\abs{\lambda_t}+\abs{\lambda_{k_0}}}\Bigr(\frac{1}{\abs{\lambda_{j_0}}
+\abs{\lambda_{k_0}}}+\frac{1}{\abs{\lambda_{t}}+\abs{\lambda_{k_0}}}\Bigr)
\frac{\pr^3\phi}{\pr z_{j_0}\pr z_t\pr\ol z_{k_0}}(0), \\
\beta^4_t=0 ,\\
\beta^5_t=2\frac{\pr^3\phi}{\pr z_t\pr\ol z_{j_0}\pr z_{j_0}}(0),\\
\beta^6_t=0. 
\end{gather}
\end{subequations} 

Now, using \eqref{s2-e40} in \eqref{s2-e39}, we can check that 
\begin{equation} \label{s2-e45} 
\begin{split} 
&2i\bigr(\abs{\lambda_{j_0}}+\abs{\lambda_{k_0}}\bigr)\frac{\pr^4\psi_4}{\pr\ol z_{j_0}\pr z_{j_0}\pr\ol z_{k_0}\pr z_{k_0}}\abs{z_{j_0}}^2\abs{z_{k_0}}^2\\
&=-\sum^n_{t=1}\bigr(\alpha^1_t\beta^4_t+\alpha^2_t\beta^3_t+\alpha^3_t\beta^2_t
+\alpha^4_t\beta^1_t+\alpha^5_t\beta^6_t+\alpha^6_t\beta^5_t\bigr)\abs{z_{j_0}}^2\abs{z_{k_0}}^2\\
&\quad-2\bigr(\lambda_{j_0}+\lambda_{k_0}\bigr)\frac{\pr^4\phi_4}{\pr\ol z_{j_0}\pr z_{j_0}\pr\ol z_{k_0}\pr z_{k_0}}\abs{z_{j_0}}^2\abs{z_{k_0}}^2+F(z),
\end{split}
\end{equation} 
where $\frac{\pr^4F}{\pr\ol z_{j_0}\pr z_{j_0}\pr\ol z_{k_0}\pr z_{k_0}}=0$. 
We recall that $T$ is given by \eqref{s2-e12}.
From \eqref{s2-e45}, we have 
\begin{equation} \label{s2-e46} 
\begin{split}
&\frac{\pr^4\psi_4}{\pr\ol z_{j_0}\pr z_{j_0}\pr\ol z_{k_0}\pr z_{k_0}}\\
&=\frac{i}{2\bigr(\abs{\lambda_{j_0}}+\abs{\lambda_{k_0}}\bigr)}
\sum^n_{t=1}\bigr(\alpha^1_t\beta^4_t+\alpha^2_t\beta^3_t+\alpha^3_t\beta^2_t
+\alpha^4_t\beta^1_t+\alpha^5_t\beta^6_t+\alpha^6_t\beta^5_t\bigr) \\
&\quad+i\frac{\lambda_{j_0}+\lambda_{k_0}}{\abs{\lambda_{j_0}}+\abs{\lambda_{k_0}}}\frac{\pr^4\phi_4}{\pr\ol z_{j_0}\pr z_{j_0}\pr \ol z_{k_0}\pr z_{k_0}}. 
\end{split}
\end{equation}
From \eqref{s2-e41}--\eqref{s2-e44} and \eqref{s2-e46}, it is straight forward to see 
that
\begin{equation} \label{s2-e47} 
\begin{split}
\frac{\pr^4\psi_4}{\pr\ol z_{j_0}\pr z_{j_0}\pr\ol z_{k_0}\pr z_{k_0}}
&=\frac{\pr^4\psi}{\pr\ol z_{j_0}\pr z_{j_0}\pr\ol z_{k_0}\pr z_{k_0}}(0,0) \\
&=\sum^q_{t=1}\Bigr(\frac{2i\abs{\lambda_{j_0}}}{\bigr(\abs{\lambda_{j_0}}+\abs{\lambda_{k_0}}\bigr)\bigr(\abs{\lambda_t}+
\abs{\lambda_{j_0}}\bigr)}\abs{\frac{\pr^3\phi}{\pr\ol z_t\pr z_{j_0}\pr z_{k_0}}(0)}^2\\
&\quad +\frac{2i\abs{\lambda_{t}}\abs{\lambda_{k_0}}\abs{\lambda_{j_0}}}{\bigr(\abs{\lambda_{j_0}}+\abs{\lambda_{k_0}}\bigr)
\bigr(\abs{\lambda_{t}}+\abs{\lambda_{k_0}}+\abs{\lambda_{j_0}}\bigr)^2} \\
&\quad\times\bigr(\frac{1}{\abs{\lambda_{j_0}}+\abs{\lambda_t}}+\frac{1}{\abs{\lambda_{j_0}}+\abs{\lambda_{k_0}}}\bigr)\abs{\frac{\pr^3\phi}
{\pr\ol z_t\pr z_{j_0}\pr\ol z_{k_0}}(0)}^2 \\
&\quad+\frac{2i\abs{\lambda_{j_0}}}{(\abs{\lambda_{j_0}}+\abs{\lambda_t})(\abs{\lambda_{j_0}}+\abs{\lambda_{k_0}})}\frac{\pr^3\phi}{\pr z_t\pr\ol z_{j_0}\pr z_{j_0}}(0)\frac{\pr^3\phi}{\pr\ol z_t\pr\ol z_{k_0}\pr z_{k_0}}(0)\Bigr)\\
&\quad +\sum^n_{t=q+1}\Bigr(\frac{2i\abs{\lambda_{k_0}}}{\bigr(\abs{\lambda_{j_0}}+\abs{\lambda_{k_0}}\bigr)
\bigr(\abs{\lambda_t}+\abs{\lambda_{k_0}}\bigr)}\abs{\frac{\pr^3\phi}{\pr\ol z_t\pr z_{j_0}\pr z_{k_0}}(0)}^2\\
&\quad +\frac{2i\abs{\lambda_{t}}\abs{\lambda_{k_0}}\abs{\lambda_{j_0}}}{\bigr(\abs{\lambda_{j_0}}+\abs{\lambda_{k_0}}\bigr)
\bigr(\abs{\lambda_{t}}+\abs{\lambda_{k_0}}+\abs{\lambda_{j_0}}\bigr)^2} \\
&\quad\times\bigr(\frac{1}{\abs{\lambda_{j_0}}+\abs{\lambda_{k_0}}}+\frac{1}{\abs{\lambda_t}+\abs{\lambda_{k_0}}}\bigr)\abs{\frac{\pr^3\phi}{\pr\ol z_t\pr \ol z_{j_0}\pr z_{k_0}}(0)}^2 \\ 
&\quad+\frac{2i\abs{\lambda_{k_0}}}{(\abs{\lambda_{k_0}}+\abs{\lambda_t})(\abs{\lambda_{j_0}}+\abs{\lambda_{k_0}})}\frac{\pr^3\phi}{\pr z_t\pr\ol z_{j_0}\pr z_{j_0}}(0)\frac{\pr^3\phi}{\pr\ol z_t\pr\ol z_{k_0}\pr z_{k_0}}(0)\Bigr)\\
&\quad +i\frac{\lambda_{j_0}+\lambda_{k_0}}{\abs{\lambda_{j_0}}+\abs{\lambda_{k_0}}}\frac{\pr^4\phi}{\pr\ol z_{j_0}\pr z_{j_0}\pr \ol z_{k_0}\pr z_{k_0}}(0).
\end{split}
\end{equation} 
Similarly, we can repeat the procedure above with minor change and get 
\begin{equation} \label{s2-e48} 
\begin{split}
\frac{\pr^4\psi}{\pr\ol z_{j_0}\pr z_{j_0}\pr\ol z_{k_0}\pr z_{k_0}}(0,0)
&=\sum^q_{t=1}\Bigr(\frac{2i}{\bigr(\abs{\lambda_t}+\abs{\lambda_{j_0}}+
\abs{\lambda_{k_0}}\bigr)^2}\bigr(\abs{\lambda_t}+\abs{\lambda_{j_0}}+
\abs{\lambda_{k_0}}\\
&-\frac{\abs{\lambda_{j_0}}\abs{\lambda_{k_0}}}{\abs{\lambda_{j_0}}+\abs{\lambda_t}}-\frac{\abs{\lambda_{j_0}}\abs{\lambda_{k_0}}}{\abs{\lambda_{k_0}}+\abs{\lambda_t}}\bigr)\abs{\frac{\pr^3\phi}{\pr z_t\pr \ol z_{j_0}\pr\ol z_{k_0}}(0)}^2\Bigr)\\
&+i\frac{\pr^4\phi}{\pr\ol z_{j_0}\pr z_{j_0}\pr\ol z_{k_0}\pr z_{k_0}}(0),
\end{split}
\end{equation} 
where $q+1\leq j_0,k_0\leq n$, and 
\begin{equation} \label{s2-e49} 
\begin{split}
\frac{\pr^4\psi}{\pr\ol z_{j_0}\pr z_{j_0}\pr\ol z_{k_0}\pr z_{k_0}}(0,0)
&=\sum^n_{t=q+1}\Bigr(\frac{2i}{\bigr(\abs{\lambda_t}+\abs{\lambda_{j_0}}+
\abs{\lambda_{k_0}}\bigr)^2}\bigr(\abs{\lambda_t}+\abs{\lambda_{j_0}}+
\abs{\lambda_{k_0}}\\
&-\frac{\abs{\lambda_{j_0}}\abs{\lambda_{k_0}}}{\abs{\lambda_{k_0}}+\abs{\lambda_t}}-\frac{\abs{\lambda_{j_0}}
\abs{\lambda_{k_0}}}{\abs{\lambda_{j_0}}+\abs{\lambda_t}}\bigr)\abs{\frac{\pr^3\phi}{\pr\ol z_t\pr z_{j_0}\pr z_{k_0}}(0)}^2\Bigr)\\
&-i\frac{\pr^4\phi}{\pr\ol z_{j_0}\pr z_{j_0}\pr\ol z_{k_0}\pr z_{k_0}}(0),
\end{split} 
\end{equation} 
where $1\leq j_0,k_0\leq q$. 

From \eqref{s2-e47}, \eqref{s2-e48}, \eqref{s2-e49} and Proposition~\ref{s2-p2}, 
Theorem~\ref{s2-t1} follows. 

\section{The transport equations for $\Box^{(q)}_k$} 

\subsection{The transport equations} 

In this section, we will write down the transport equations for $\Box^{(q)}_k$ and we will solve the first transport equation at $z=w$ in some sence. The main reference for this section is~\cite{BS05}. We first derive representations for $\ddbar$, $\ol\pr^*_k$ in spaces without exponential weights, by using the following unitary 
identifications: 
\begin{equation} \label{s3-e1} 
\left\{\begin{aligned}  
L^2_q(\Complex^n)&\leftrightarrow L^2_{q,k}(\Complex^n) \\
u&\leftrightarrow \hat u=e^{k\phi}u.
\end{aligned}
\right. 
\end{equation} 
Using \eqref{s3-e1}, we get 
\begin{equation} \label{s3-e2} 
\ddbar\hat u=e^{k\phi}\ddbar_su,
\end{equation}
where
\begin{equation} \label{s3-e3}
\ddbar_s=\sum^n_{j=1}\Bigr(d\ol z_j^\wedge\circ(\frac{\pr}{\pr\ol z_j}+k\frac{\pr\phi}{\pr\ol z_j})\Bigr).
\end{equation}
Now the formal adjoint of $\ddbar_s$ for the inner product $(\ |\ )$ given by \eqref{s1-e2} is 
\begin{equation} \label{s3-e4}
\ol{\pr}^*_s=\sum^n_{j=1}\Bigr(d\ol z_j^{\wedge,*}\circ(-\frac{\pr}{\pr z_j}+k\frac{\pr\phi}{\pr z_j})\Bigr),
\end{equation}
where in view of the unitarity of the relation \eqref{s3-e1}, 
\begin{equation} \label{s3-e5}
\ol{\pr}^*\hat u=e^{k\phi}\ol{\pr}^*_su.
\end{equation}
We can identify the Kodaira Laplacian with 
\begin{equation} \label{s3-e6}
\hat\Box^{(q)}_k=\ddbar_s\ol{\pr}^*_s+\ol{\pr}^*_s\ddbar_s.
\end{equation}
Put 
\begin{equation} \label{s3-e7}
\hat\Pi^{(q)}_k:L^2_q(\Complex^n)\To{\rm Ker\,}\hat\Box^{(q)}_k
\end{equation}
be the orthogonal projection with respect to $(\ |\ )$ and let $\hat\Pi^{(q)}_k(z,w)$ be the distribution kernel of $\hat\Pi^{(q)}_k$. From \eqref{s3-e2} and \eqref{s3-e5}, we have 
\begin{equation} \label{s3-e8}
\hat\Pi^{(q)}_k(z,w)=e^{-k\phi(z)}\Pi^{(q)}_k(z,w)e^{k\phi(w)}.
\end{equation}
In view of Theorem~\ref{s1-t1}, we see that
\begin{equation} \label{s3-e9}
\hat\Pi^{(q)}_k(z,w)=e^{ik\psi(z,w)}b(z,w,k)+S(z,w,k), 
\end{equation}
where $b(z. w, k)\sim\sum^\infty_ob_j(z,w)k^{n-j}$ in $C^\infty(\Complex^n\times\Complex^n;\, \mathscr L(\Lambda^{0,q}T^*_w(\Complex^n),
\Lambda^{0,q}T^*_z(\Complex^n)))$ and $S(z,w,k)$ is negligible.

From \eqref{s3-e3} and \eqref{s3-e9}, we have
\begin{equation} \label{s3-e10}
e^{-ik\psi}\ddbar_s\hat\Pi^{(q)}_k(z,w)=\sum^n_{j=1}k(i\frac{\pr\psi}{\pr\ol z_j}+\frac{\pr\phi}{\pr\ol z_j})d\ol z_j^{\wedge}b+\sum^n_{j=1}d\ol z_j^\wedge\frac{\pr b}{\pr\ol z_j}.
\end{equation} 
From this and \eqref{s3-e4}, we can compute 
\begin{equation} \label{s3-e11}
\begin{split}
e^{-ik\psi}\ol{\pr}^*_s\ddbar_s\hat\Pi^{(q)}_k(z,w)=&k^2\sum^n_{j,t=1}(-i\frac{\pr\psi}{\pr z_t}+\frac{\pr\phi}{\pr z_t})
(i\frac{\pr\psi}{\pr\ol z_j}+\frac{\pr\phi}{\pr\ol z_j})d\ol z_t^{\wedge,*}d\ol z_j^\wedge b \\
&-k\sum^n_{j,t=1}\frac{\pr}{\pr z_t}(i\frac{\pr\psi}{\pr\ol z_j}+\frac{\pr\phi}{\pr\ol z_j})d\ol z_t^{\wedge,*}d\ol z_j^\wedge b \\
&-k\sum^n_{j,t=1}(i\frac{\pr\psi}{\pr\ol z_j}+\frac{\pr\phi}{\pr\ol z_j})d\ol z_t^{\wedge,*}d\ol z_j^\wedge\frac{\pr b}{\pr z_t} \\
&+k\sum^n_{j,t=1}(-i\frac{\pr\psi}{\pr z_t}+\frac{\pr\phi}{\pr z_t})d\ol z_t^{\wedge,*}d\ol z_j^\wedge\frac{\pr b}{\pr\ol z_j}\\
&-\sum^n_{j,t=1}d\ol z_t^{\wedge,*}d\ol z_j^\wedge\frac{\pr^2b}{\pr z_t\pr\ol z_j}.
\end{split}
\end{equation}
Similarly, we have 
\begin{equation} \label{s3-e12}
\begin{split}
e^{-ik\psi}\ddbar_s\ol{\pr}^*_s\hat\Pi^{(q)}_k(z,w)=&k^2\sum^n_{j,t=1}(i\frac{\pr\psi}{\pr\ol z_j}+\frac{\pr\phi}{\pr\ol z_j})(-i\frac{\pr\psi}{\pr z_t}+\frac{\pr\phi}{\pr z_t})
d\ol z_j^\wedge d\ol z_t^{\wedge,*} b \\ 
&+k\sum^n_{j,t=1}\frac{\pr}{\pr\ol z_j}(-i\frac{\pr\psi}{\pr z_t}+\frac{\pr\phi}{\pr z_t})d\ol z_j^\wedge d\ol z_t^{\wedge,*} b \\
&+k\sum^n_{j,t=1}(-i\frac{\pr\psi}{\pr z_t}+\frac{\pr\phi}{\pr z_t})d\ol z_j^\wedge d\ol z_t^{\wedge,*}\frac{\pr b}{\pr\ol z_j} \\
&-k\sum^n_{j,t=1}(i\frac{\pr\psi}{\pr\ol z_j}+\frac{\pr\phi}{\pr\ol z_j})d\ol z_j^\wedge d\ol z_t^{\wedge,*}\frac{\pr b}{\pr z_t}\\
&-\sum^n_{j,t=1}d\ol z_j^\wedge d\ol z_t^{\wedge,*}\frac{\pr^2b}{\pr z_t\pr\ol z_j}.
\end{split} 
\end{equation}
Combining \eqref{s3-e11}, \eqref{s3-e12} with $d\ol z_j^\wedge\circ d\ol z_k^{\wedge,*}+d\ol z_k^{\wedge,*}\circ d\ol z_j^\wedge=\delta_{j,k}$, we get 
\begin{equation} \label{s3-e13}
\begin{split}
e^{-ik\psi}\hat\Box^{(q)}_k\hat\Pi^{(q)}_k(z,w)=&k^2\sum^n_{j=1}(i\frac{\pr\psi}{\pr\ol z_j}+
\frac{\pr\phi}{\pr\ol z_j})(-i\frac{\pr\psi}{\pr z_j}+\frac{\pr\phi}{\pr z_j})b  \\
&+k\sum^n_{j,t=1}\frac{\pr}{\pr\ol z_j}(-i\frac{\pr\psi}{\pr z_t}+\frac{\pr\phi}{\pr z_t})d\ol z_j^\wedge d\ol z_t^{\wedge,*} b \\
&-k\sum^n_{j,t=1}\frac{\pr}{\pr z_t}(i\frac{\pr\psi}{\pr\ol z_j}+\frac{\pr\phi}{\pr\ol z_j})d\ol z_t^{\wedge,*} d\ol z_j^\wedge b \\
&+k\sum^n_{j=1}(-i\frac{\pr\psi}{\pr z_j}+\frac{\pr\phi}{\pr z_j})\frac{\pr b}{\pr\ol z_j} \\
&-k\sum^n_{j=1}(i\frac{\pr\psi}{\pr\ol z_j}+\frac{\pr\phi}{\pr\ol z_j})\frac{\pr b}{\pr z_j}\\
&-\sum^n_{j=1}\frac{\pr^2b}{\pr z_j\pr\ol z_j}.
\end{split} 
\end{equation}
We regroup \eqref{s3-e13} according to the degree of $k$ and notice that the leading term in \eqref{s3-e13} vanishes to infinite order on $z=w$ and $\hat\Box^{(q)}_k\hat\Pi^{(q)}_k=0$, we obtain the following

\begin{prop} \label{s3-p1}
We have 
\begin{equation} \label{s3-e14}
\begin{split}
&\sum^n_{j,t=1}\frac{\pr}{\pr\ol z_j}(-i\frac{\pr\psi}{\pr z_t}+\frac{\pr\phi}{\pr z_t})d\ol z_j^\wedge d\ol z_t^{\wedge,*} b_0 
-\sum^n_{j,t=1}\frac{\pr}{\pr z_t}(i\frac{\pr\psi}{\pr\ol z_j}+\frac{\pr\phi}{\pr\ol z_j})d\ol z_t^{\wedge,*} d\ol z_j^\wedge b_0 \\
&+\sum^n_{j=1}(-i\frac{\pr\psi}{\pr z_j}+\frac{\pr\phi}{\pr z_j})\frac{\pr b_0}{\pr\ol z_j} 
-\sum^n_{j=1}(i\frac{\pr\psi}{\pr\ol z_j}+\frac{\pr\phi}{\pr\ol z_j})\frac{\pr b_0}{\pr z_j}
\end{split}
\end{equation}
vanishes to infinite order on $z=w$ and 
\begin{equation} \label{s3-e15}
\begin{split}
&\sum^n_{j,t=1}\frac{\pr}{\pr\ol z_j}(-i\frac{\pr\psi}{\pr z_t}+\frac{\pr\phi}{\pr z_t})d\ol z_j^\wedge d\ol z_t^{\wedge,*} b_1 -\sum^n_{j,t=1}\frac{\pr}{\pr z_t}(i\frac{\pr\psi}{\pr\ol z_j}+\frac{\pr\phi}{\pr\ol z_j})d\ol z_t^{\wedge,*} d\ol z_j^\wedge b_1 \\
&+\sum^n_{j=1}(-i\frac{\pr\psi}{\pr z_j}+\frac{\pr\phi}{\pr z_j})\frac{\pr b_1}{\pr\ol z_j} 
-\sum^n_{j=1}(i\frac{\pr\psi}{\pr\ol z_j}+\frac{\pr\phi}{\pr\ol z_j})\frac{\pr b_1}{\pr z_j}-\sum^n_{j=1}\frac{\pr^2b_0}{\pr z_j\pr\ol z_j}
\end{split}
\end{equation} 
vanishes to infinite order on $z=w$.
\end{prop}

\subsection{The first order of the Taylor expansion of $b_0(z,w)$ at $z=w$} 

Now, as in section 2, we assume that $\phi(z)=\sum^n_{j=1}\lambda_j\abs{z_j}^2+O(\abs{z}^3)$ near $0$ and that $\lambda_j<0$, $j=1,\ldots,q$ and $\lambda_j>0$, $j=q+1,\ldots,n$. We work in some neighborhood of $(0,0)$. Put 
\begin{equation} \label{s3-e16}
b_0(z,w)=b ^0_0(z,w)+b^1_0(z,w)+b^2_0(z,w)+\cdots,
\end{equation}
where $b^j_0$ is a homogeneous polynomialof degree $j$ in $(z,w)$. We recall that 
\begin{equation*} 
b^0_0=\pi^{-n}\abs{\lambda_1}\abs{\lambda_2}\cdots\abs{\lambda_n}\prod^q_{j=1}d\ol z_j^\wedge d\ol z_j^{\wedge,*}.
\end{equation*} 
(See Theorem~\ref{s1-t2}.)
For $1\leq s\leq n$, put $\chi_1(s)=1$ if $1\leq s\leq q$ and $\chi_1(s)=0$ if $q+1\leq s\leq n$ and put $\chi_2(s)=1-\chi_1(s)$. We recall that for any operator $T\in\mathscr L(\Lambda^{p,q}T^*(\Complex^n), \Lambda^{p,q}T^*(\Complex^n))$, the trace of $T$ is given by  \eqref{s1-e-1}. The first goal of this section is to prove the following 

\begin{thm} \label{s3-t1}  
Under the notations above, we have 
\begin{equation} \label{s3-e17}
\begin{split}
b^1_0(z,0)=&\sum_{q+1\leq j\leq n,1\leq k\leq q,1\leq s\leq n}\frac{1}{\abs{\lambda_j}+\abs{\lambda_k}+\abs{\lambda_s}\chi_1(s)}\\
&\times\frac{\pr^3\phi}{\pr\ol z_j\pr z_k\pr z_s}(0)z_sd\ol z_k^{\wedge,*}d\ol z_j^\wedge b^0_0\\
&+\sum_{q+1\leq j\leq n,1\leq k\leq q,1\leq s\leq n}\frac{1}{\abs{\lambda_j}+\abs{\lambda_k}+\abs{\lambda_s}\chi_2(s)}\\
&\times\frac{\pr^3\phi}{\pr\ol z_j\pr z_k\pr\ol z_s}(0)\ol z_sd\ol z_k^{\wedge,*}d\ol z_j^\wedge b^0_0\\ 
&+\sum_{q+1\leq j\leq n, 1\leq k\leq q, q+1\leq s\leq n}\frac{\abs{\lambda_s}}{(\abs{\lambda_j}+\abs{\lambda_k})(\abs{\lambda_j}+\abs{\lambda_k}+\abs{\lambda_s})}\\
&\times\frac{\pr^3\phi}{\pr z_j\pr\ol z_k\pr z_s}(0)z_sb^0_0d\ol z_j^{\wedge,*}d\ol z_k^\wedge\\
&+\sum_{q+1\leq j\leq n, 1\leq k\leq q, 1\leq s\leq q}\frac{\abs{\lambda_s}}{(\abs{\lambda_j}+\abs{\lambda_k})(\abs{\lambda_j}+\abs{\lambda_k}+\abs{\lambda_s})}\\
&\times\frac{\pr^3\phi}{\pr z_j\pr\ol z_k\pr\ol z_s}(0)\ol z_sb^0_0d\ol z_j^{\wedge,*}d\ol z_k^\wedge \\
&-\sum_{q+1\leq s\leq n,1\leq j\leq q}\frac{1}{\abs{\lambda_j}+\abs{\lambda_s}}\frac{\pr^3\phi}{\pr\ol z_j\pr z_j\pr\ol z_s}(0)\ol z_sb^0_0 \\
&-\sum_{1\leq s,j\leq q}\frac{1}{\abs{\lambda_j}}\frac{\pr^3\phi}{\pr\ol z_j\pr z_j\pr\ol z_s}(0)\ol z_sb^0_0\\
&+\sum_{1\leq s\leq q,q+1\leq j\leq n}\frac{\abs{\lambda_s}}{\abs{\lambda_j}(\abs{\lambda_j}+\abs{\lambda_s})}\frac{\pr^3\phi}{\pr\ol z_j\pr z_j\pr\ol z_s}(0)\ol z_sb^0_0  \\
&+\sum_{1\leq s\leq q,q+1\leq j\leq n}\frac{1}{\abs{\lambda_j}+\abs{\lambda_s}}\frac{\pr^3\phi}{\pr\ol z_j\pr z_j\pr z_s}(0) z_sb^0_0\\
&-\sum_{q+1\leq s\leq n,1\leq j\leq q}\frac{\abs{\lambda_s}}{\abs{\lambda_j}(\abs{\lambda_j}+\abs{\lambda_s})}\frac{\pr^3\phi}{\pr\ol z_j\pr z_j\pr z_s}(0) z_sb^0_0 \\
&+\sum_{q+1\leq s,j\leq n}\frac{1}{\abs{\lambda_j}}\frac{\pr^3\phi}{\pr\ol z_j\pr z_j\pr z_s}(0) z_sb^0_0.
\end{split}
\end{equation} 
\end{thm} 

As in section 2, we write 
\begin{equation} \label{s3-e18} 
\begin{split}
\psi(z,w)&=i\sum^n_{j=1}\abs{\lambda_j}\abs{z_j-w_j}^2+i\sum^n_{j=1}\lambda_j(\ol z_jw_j-\ol w_jz_j) \\
&\quad+\psi_3(z, w)+\psi_4(z, w)+\cdots,
\end{split}
\end{equation} 
and
\begin{equation} \label{s3-e19} 
\phi(z)=\sum^n_{j=1}\lambda_j\abs{z_j}^2+\phi_3(z)+\phi_4(z)+\cdots,
\end{equation} 
where $\psi_j(z, w)$ is a homogeneous polynomial of degree $j$ in $(z, w)$, $j=3,4,\ldots$, $\phi_j(z)$ is a homogeneous polynomial of degree $j$  in $z$, $j=3,4,\ldots$. Now, using \eqref{s3-e16}, \eqref{s3-e18} and \eqref{s3-e19} in \eqref{s3-e14}, we get 
\begin{equation} \label{s3-e20}
\begin{split}
&\Bigr(\sum^n_{j=q+1}2\abs{\lambda_j}d\ol z_j^\wedge d\ol z_j^{\wedge,*}+\sum^n_{j,t=1}
\frac{\pr}{\pr\ol z_j}(-i\frac{\pr\psi_3(z,w)}{\pr z_t}+\frac{\pr\phi_3(z)}{\pr z_t})d\ol z_j^\wedge d\ol z_t^{\wedge,*}\\
&+\sum^n_{j,t=1}\frac{\pr}{\pr\ol z_j}(-i\frac{\pr\psi_4(z,w)}{\pr z_t}+\frac{\pr\phi_4(z)}{\pr z_t})d\ol z_j^\wedge d\ol z_t^{\wedge,*}\Bigr)
(b^0_0+b^1_0(z,w)+b^2_0(z,w))\\
&+\Bigr(\sum^q_{j=1}2\abs{\lambda_j}d\ol z_j^{\wedge,*} d\ol z_j^\wedge-
\sum^n_{j,t=1}\frac{\pr}{\pr z_t}(i\frac{\pr\psi_3(z,w)}{\pr\ol z_j}+\frac{\pr\phi_3(z)}{\pr\ol z_j})d\ol z_t^{\wedge,*}d\ol z_j^\wedge \\
&-\sum^n_{j,t=1}
\frac{\pr}{\pr z_t}(i\frac{\pr\psi_4(z,w)}{\pr\ol z_j}+\frac{\pr\phi_4(z)}{\pr\ol z_j})d\ol z_t^{\wedge,*} d\ol z_j^\wedge\Bigr)
(b^0_0+b^1_0(z,w)+b^2_0(z,w))\\
&+\sum^n_{j=q+1}2\abs{\lambda_j}(\ol z_j-\ol w_j)(\frac{\pr b^1_0(z,w)}{\pr\ol z_j}+\frac{\pr b^2_0(z,w)}{\pr\ol z_j})+\sum^n_{j=1}(-i\frac{\pr\psi_3}{\pr z_j}+\frac{\pr\phi_3}{\pr z_j})\frac{\pr b^1_0(z,w)}{\pr\ol z_j}\\
&+\sum^q_{j=1}2\abs{\lambda_j}(z_j-w_j)(\frac{\pr b^1_0(z,w)}{\pr z_j}+\frac{\pr b^2_0(z,w)}{\pr z_j})-\sum^n_{j=1}(i\frac{\pr\psi_3}{\pr\ol z_j}+
\frac{\pr\phi_3}{\pr\ol z_j})\frac{\pr b^1_0(z,w)}{\pr z_j}\\
&=O(\abs{(z,w)}^3).
\end{split}
\end{equation} 
It is straight forward to see that the order $1$ and $2$ terms in \eqref{s3-e20} are the following 
\begin{equation} \label{s3-e21}
\begin{split}
Lb^1_0(z,w)=&-\Bigr(\sum^n_{j,t=1}
\frac{\pr}{\pr\ol z_j}(-i\frac{\pr\psi_3(z,w)}{\pr z_t}+\frac{\pr\phi_3(z)}{\pr z_t})d\ol z_j^\wedge d\ol z_t^{\wedge,*}\\
&-\sum^n_{j,t=1}
\frac{\pr}{\pr z_t}(i\frac{\pr\psi_3(z,w)}{\pr\ol z_j}+\frac{\pr\phi_3(z)}{\pr\ol z_j})d\ol z_t^{\wedge,*} d\ol z_j^\wedge \Bigr)b^0_0\\
&+\sum^n_{j=q+1}2\abs{\lambda_j}\ol w_j\frac{\pr b^1_0(z,w)}{\pr\ol z_j}+\sum^q_{j=1}2\abs{\lambda_j}w_j\frac{\pr b^1_0(z,w)}{\pr z_j}
\end{split}
\end{equation} 
and 
\begin{equation} \label{s3-e22}
\begin{split}
Lb^2_0(z,w)=&-\Bigr(\sum^n_{j,t=1}
\frac{\pr}{\pr\ol z_j}(-i\frac{\pr\psi_4(z,w)}{\pr z_t}+\frac{\pr\phi_4(z)}{\pr z_t})d\ol z_j^\wedge d\ol z_t^{\wedge,*}\\
&-\sum^n_{j,t=1}
\frac{\pr}{\pr z_t}(i\frac{\pr\psi_4(z,w)}{\pr\ol z_j}+\frac{\pr\phi_4(z)}{\pr\ol z_j})d\ol z_t^{\wedge,*} d\ol z_j^\wedge\Bigr)b^0_0\\
&-\Bigr(\sum^n_{j,t=1}
\frac{\pr}{\pr\ol z_j}(-i\frac{\pr\psi_3(z,w)}{\pr z_t}+\frac{\pr\phi_3(z)}{\pr z_t})d\ol z_j^\wedge d\ol z_t^{\wedge,*}\\
&-\sum^n_{j,t=1}
\frac{\pr}{\pr z_t}(i\frac{\pr\psi_3(z,w)}{\pr\ol z_j}+\frac{\pr\phi_3(z)}{\pr\ol z_j})d\ol z_t^{\wedge,*}d\ol z_j^\wedge\Bigr)b^1_0(z,w) \\
&-\sum^n_{j=1}(-i\frac{\pr\psi_3(z,w)}{\pr z_j}+\frac{\pr\phi_3(z)}{\pr z_j})\frac{\pr b^1_0(z,w)}{\pr\ol z_j}\\
&+\sum^n_{j=1}(i\frac{\pr\psi_3(z,w)}{\pr\ol z_j}+\frac{\pr\phi_3(z)}{\pr\ol z_j})\frac{\pr b^1_0(z,w)}{\pr z_j}\\
&+\sum^n_{j=q+1}2\abs{\lambda_j}\ol w_j\frac{\pr b^2_0(z,w)}{\pr\ol z_j}+
\sum^q_{j=1}2\abs{\lambda_j}w_j\frac{\pr b^2_0(z,w)}{\pr z_j},
\end{split}
\end{equation} 
where 
\begin{equation} \label{s3-e23}
L=\sum^n_{j=q+1}2\abs{\lambda_j}d\ol z_j^\wedge d\ol z_j^{\wedge,*}+
\sum^q_{j=1}2\abs{\lambda_j}d\ol z_j^{\wedge,*} d\ol z_j^\wedge
+\sum^n_{j=q+1}2\abs{\lambda_j}\ol z_j\frac{\pr}{\pr\ol z_j}+\sum^q_{j=1}2\abs{\lambda_j}z_j\frac{\pr}{\pr z_j}.
\end{equation}
We rewrite first term of the right side of \eqref{s3-e21}:
\begin{equation} \label{s3-e24}
\begin{split} 
&-\Bigr(\sum^n_{j,t=1}
\frac{\pr}{\pr\ol z_j}(-i\frac{\pr\psi_3(z,w)}{\pr z_t}+\frac{\pr\phi_3(z)}{\pr z_t})d\ol z_j^\wedge d\ol z_t^{\wedge,*}\\
&-\sum^n_{j,t=1}
\frac{\pr}{\pr z_t}(i\frac{\pr\psi_3(z,w)}{\pr\ol z_j}+\frac{\pr\phi_3(z)}{\pr\ol z_j})d\ol z_t^{\wedge,*} d\ol z_j^\wedge \Bigr)b^0_0\\
&=-\sum^n_{j=1}\frac{\pr}{\pr\ol z_j}(-i\frac{\pr\psi_3(z,w)}{\pr z_j}+\frac{\pr\phi_3(z)}{\pr z_j})d\ol z_j^\wedge d\ol z_j^{\wedge,*}b^0_0\\
&+\sum^n_{j=1}\frac{\pr}{\pr z_j}(i\frac{\pr\psi_3(z,w)}{\pr\ol z_j}+\frac{\pr\phi_3(z)}{\pr\ol z_j})d\ol z_j^{\wedge,*} d\ol z_j^\wedge b^0_0\\
&+\sum^n_{j,t=1,j\neq t}i\frac{\pr^2\psi_3(z,w)}{\pr\ol z_j\pr z_t}(d\ol z_j^\wedge d\ol z_t^{\wedge,*}+d\ol z_t^{\wedge,*}d\ol z_j^\wedge)b^0_0 \\
&-\sum^n_{j,t=1.j\neq t}\frac{\pr^2\phi_3(z)}{\pr\ol z_j\pr z_t}(d\ol z_j^\wedge d\ol z_t^{\wedge,*}-d\ol z_t^{\wedge,*}d\ol z_j^\wedge)b^0_0.
\end{split}
\end{equation} 
Note that 
$d\ol z_j^\wedge d\ol z_t^{\wedge,*}+d\ol z_t^{\wedge,*}d\ol z_j^\wedge=\delta_{j,t}$. From the form of $b^0_0$ (see Theorem~\ref{s1-t2}), we can check that 
$d\ol z_j^\wedge d\ol z_j^{\wedge,*}b^0_0= 0$ if $q+1\leq j\leq n$, 
$d\ol z_j^\wedge d\ol z_j^{\wedge,*}b^0_0=b^0_0$ if $1\leq j\leq q$,
$d\ol z_j^{\wedge,*}d\ol z_j^\wedge b^0_0= 0$ if $1\leq j\leq q$, 
$d\ol z_j^{\wedge,*}d\ol z_j^\wedge b^0_0=b^0_0$ if $q+1\leq j\leq n$ and 
when $j\neq t$, $d\ol z_t^{\wedge,*}d\ol z_j^\wedge b^0_0\neq0$ if and only if 
$q+1\leq j\leq n$ and $1\leq t\leq q$. From this observaton, \eqref{s3-e24} becomes 
\begin{equation} \label{s3-e25} 
\begin{split} 
&-\Bigr(\sum^n_{j,t=1}
\frac{\pr}{\pr\ol z_j}(-i\frac{\pr\psi_3(z,w)}{\pr z_t}+\frac{\pr\phi_3(z)}{\pr z_t})d\ol z_j^\wedge d\ol z_t^{\wedge,*}\\
&-\sum^n_{j,t=1}
\frac{\pr}{\pr z_t}(i\frac{\pr\psi_3(z,w)}{\pr\ol z_j}+\frac{\pr\phi_3(z)}{\pr\ol z_j})d\ol z_t^{\wedge,*} d\ol z_j^\wedge \Bigr)b^0_0\\
&=-\sum^q_{j=1}\frac{\pr}{\pr\ol z_j}(-i\frac{\pr\psi_3(z,w)}{\pr z_j}+\frac{\pr\phi_3(z)}{\pr z_j})b^0_0\\
&+\sum^n_{j=q+1}\frac{\pr}{\pr z_j}(i\frac{\pr\psi_3(z,w)}{\pr\ol z_j}+\frac{\pr\phi_3(z)}{\pr\ol z_j})b^0_0\\
&+2\sum^n_{q+1\leq j\leq n,1\leq t\leq q}\frac{\pr^2\phi_3(z)}{\pr\ol z_j\pr z_t}
d\ol z_t^{\wedge,*}d\ol z_j^\wedge b^0_0.
\end{split}
\end{equation} 

From \eqref{s2-e30}, it is straight forward to see that 
\begin{equation} \label{s3-e26} 
\frac{\pr}{\pr\ol z_j}(-i\frac{\pr\psi_3(z,w)}
{\pr z_j}+\frac{\pr\phi_3(z)}{\pr z_j})=\sum_{q+1\leq s\leq n}\frac{2\abs{\lambda_s}}{\abs{\lambda_j}+\abs{\lambda_s}}\frac{\pr^3\phi}{\pr\ol z_j\pr z_j\pr\ol z_s}(0)(\ol z_s-\ol w_s)
\end{equation} 
where $1\leq j\leq q$ and 
\begin{equation} \label{s3-e27} 
\frac{\pr}{\pr z_j}(i\frac{\pr\psi_3(z,w)}
{\pr\ol z_j}+\frac{\pr\phi_3(z)}{\pr\ol z_j})=\sum_{1\leq s\leq q}\frac{2\abs{\lambda_s}}{\abs{\lambda_j}+\abs{\lambda_s}}\frac{\pr^3\phi}{\pr\ol z_j\pr z_j\pr z_s}(0)(z_s-w_s),
\end{equation} 
where $q+1\leq j\leq n$. From \eqref{s3-e25}, \eqref{s3-e26} and \eqref{s3-e27}, \eqref{s3-e21} becomes 
\begin{equation} \label{s3-e28} 
\begin{split}
Lb^1_0(z,w)=&-\sum_{q+1\leq s\leq n,1\leq j\leq q}\frac{2\abs{\lambda_s}}{\abs{\lambda_j}+\abs{\lambda_s}}\frac{\pr^3\phi}{\pr\ol z_j\pr z_j\pr\ol z_s}(0)(\ol z_s-\ol w_s)b^0_0\\
&+\sum_{1\leq s\leq q,q+1\leq j\leq n}\frac{2\abs{\lambda_s}}{\abs{\lambda_j}+\abs{\lambda_s}}\frac{\pr^3\phi}{\pr\ol z_j\pr z_j\pr z_s}(0)(z_s-w_s)b^0_0\\
&+2\sum^n_{q+1\leq j\leq n,1\leq t\leq q,1\leq s\leq n}
\frac{\pr^3\phi}{\pr\ol z_j\pr z_t\pr z_s}(0)z_s
d\ol z_t^{\wedge,*}d\ol z_j^\wedge b^0_0\\
&+2\sum^n_{q+1\leq j\leq n,1\leq t\leq q,1\leq s\leq n}\frac{\pr^3\phi}{\pr\ol z_j\pr z_t\pr\ol z_s}(0)\ol z_s
d\ol z_t^{\wedge,*}d\ol z_j^\wedge b^0_0\\
&+\sum^n_{j=q+1}2\abs{\lambda_j}\ol w_j\frac{\pr b^1_0(z,w)}{\pr\ol z_j}+\sum^q_{j=1}2\abs{\lambda_j}w_j\frac{\pr b^1_0(z,w)}{\pr z_j}.
\end{split}
\end{equation} 
Now, we write \eqref{s3-e28} according to the degree of homogenity in $w$, we get 

\begin{equation} \label{s3-e29} 
\begin{split} 
Lb^1_0(z,0)=&-\sum_{q+1\leq s\leq n,1\leq j\leq q}\frac{2\abs{\lambda_s}}{\abs{\lambda_j}+\abs{\lambda_s}}\frac{\pr^3\phi}{\pr\ol z_j\pr z_j\pr\ol z_s}(0)\ol z_sb^0_0\\
&+\sum_{1\leq s\leq q,q+1\leq j\leq n}\frac{2\abs{\lambda_s}}{\abs{\lambda_j}+\abs{\lambda_s}}\frac{\pr^3\phi}{\pr\ol z_j\pr z_j\pr z_s}(0)z_sb^0_0\\
&+2\sum^n_{q+1\leq j\leq n,1\leq t\leq q,1\leq s\leq n}
\frac{\pr^3\phi}{\pr\ol z_j\pr z_t\pr z_s}(0)z_s
d\ol z_t^{\wedge,*}d\ol z_j^\wedge b^0_0\\
&+2\sum^n_{q+1\leq j\leq n,1\leq t\leq q,1\leq s\leq n}\frac{\pr^3\phi}{\pr\ol z_j\pr z_t\pr\ol z_s}(0)\ol z_s
d\ol z_t^{\wedge,*}d\ol z_j^\wedge b^0_0
\end{split}
\end{equation} 
and 
\begin{equation} \label{s3-e30} 
\begin{split}
Lb^1_0(0,w)=&\sum_{q+1\leq s\leq n,1\leq j\leq q}\frac{2\abs{\lambda_s}}{\abs{\lambda_j}+\abs{\lambda_s}}\frac{\pr^3\phi}{\pr\ol z_j\pr z_j\pr\ol z_s}(0)\ol w_sb^0_0\\
&-\sum_{1\leq s\leq q,q+1\leq j\leq n}\frac{2\abs{\lambda_s}}{\abs{\lambda_j}+\abs{\lambda_s}}\frac{\pr^3\phi}{\pr\ol z_j\pr z_j\pr z_s}(0)w_sb^0_0\\
&+\sum^n_{j=q+1}2\abs{\lambda_j}\ol w_j\frac{\pr b^1_0(z,0)}{\pr\ol z_j}+\sum^q_{j=1}2\abs{\lambda_j}w_j\frac{\pr b^1_0(z,0)}{\pr z_j}. 
\end{split}
\end{equation}

We pause and introduce some notations. For multi-index $J$, we write $\abs{J}=q$, $J\nearrow$, if $J=(j_1,\ldots,j_q)$, $1\leq j_1<j_2<\cdots<j_q\leq n$. Set
$d\ol z^J=d\ol z_{j_1}\wedge\cdots\wedge d\ol z_{j_q}$, $J=(j_1,\ldots,j_q)$. Then, $d\ol z^J$, $\abs{J}=q$, $J\nearrow$, is an orthonormal basis of $\Lambda^{0,q}T^*_0(\Complex^n)$, where $0$ is the origin in $\Complex^n$. Let $M_{d\ol z^J,d\ol z^K}$, $\abs{J}=\abs{K}=q$, $J,K\nearrow$, be the $\Complex$-linear operator: 
\begin{align} \label{s3-e31} 
M_{d\ol z^J,d\ol z^K}:\Lambda^{0,q}T^*_0(\Complex^n)&\To\Lambda^{0,q}T^*_0(\Complex^n)\nonumber \\
d\ol z^J&\To d\ol z^K\nonumber \\
d\ol z^I&\To 0\ \ \mbox{if $I\neq J$}.
\end{align}
It is clear that $M_{d\ol z^J,d\ol z^K}$, $\abs{J}=\abs{K}=q$, $J,K\nearrow$, is a basis of the vector space $\mathscr L(\Lambda^{0,q}T^*_0(\Complex^n),\Lambda^{0,q}T^*_0(\Complex^n))$. For $m\in\mathbb N\bigcup\set{0}$, put 
\begin{equation*} 
\begin{split}
&P^m(\mathscr L(\Lambda^{0,q}T^*_0(\Complex^n),\Lambda^{0,q}T^*_0(\Complex^n)))\\
&\quad=\set{\sum_{\abs{\alpha}+\abs{\beta}+\abs{\gamma}+\abs{\delta}=m}A_{\alpha,\beta}\ol z^\alpha z^\beta\ol w^\gamma w^\delta;\, A_{\alpha,\beta}\in\mathscr L(\Lambda^{0,q}T^*_0(\Complex^n),\Lambda^{0,q}T^*_0(\Complex^n))}. 
\end{split}
\end{equation*} 
For multi-index $J$, $\abs{J}=q$, we define 
\begin{equation} \label{s3-e32} 
F(J)=2\sum_{j\in J,q+1\leq j\leq n}\abs{\lambda_j}+2\sum_{j\notin J,1\leq j\leq q}\abs{\lambda_j}.
\end{equation} 
Put $I_0=(1,\ldots,q)$. Note that $F(J)\neq 0$ if and only if $J\neq I_0$. We have the following

\begin{lem} \label{s3-l1} 
We use the same notations as in the discussion before Theorem~\ref{s2-t1} and before.
If we consider $L$ as the operator (We recall that $L$ is given by \eqref{s3-e23}.)
\begin{equation*} 
L:P^m(\mathscr L(\Lambda^{0,q}T^*_0(\Complex^n),\Lambda^{0,q}T^*_0(\Complex^n)))\To P^m(\mathscr L(\Lambda^{0,q}T^*_0(\Complex^n),\Lambda^{0,q}T^*_0(\Complex^n))).
\end{equation*} 
Then, 
\begin{equation} \label{s3-e33} 
\begin{split}
&{\rm Ker\,}L\\
&=\set{\sum_{\abs{\alpha}+\abs{\beta}+\abs{\gamma}+\abs{\delta}=m, (\alpha'',\beta')=0,\abs{J}=q,J\nearrow}c^{\alpha,\beta,\gamma,\delta}_{J,I_0}M_{d\ol zJ,d\ol z^{I_0}}\ol z^\alpha z^\beta\ol w^\gamma w^\delta;\, c^{\alpha,\beta,\gamma,\delta}_{J,I_0}\in\Complex}.
\end{split}
\end{equation} 
Moreover, for $A\in P^m(\mathscr L(\Lambda^{0,q}T^*_0(\Complex^n),\Lambda^{0,q}T^*_0(\Complex^n)))$, we write 
\begin{equation*}
A=\sum_{\abs{\alpha}+\abs{\beta}+\abs{\gamma}+\abs{\delta}=m, \abs{J}=\abs{K}=q,J,K\nearrow}c^{\alpha,\beta,\gamma,\delta}_{J,K}
M_{d\ol z^J,d\ol z^K}\ol z^\alpha z^\beta\ol w^\gamma w^\delta,
\end{equation*}
$c^{\alpha,\beta,\gamma,\delta}_{J,K}\in\Complex$. If 
$c^{\alpha,\beta,\gamma,\delta}_{J,I_0}=0$ when $(\alpha'',\beta')=0$.
Then, we have $LB=A$, where 
\begin{equation} \label{s3-e34} 
\begin{split}
B=&\sum_{\abs{\alpha}+\abs{\beta}+\abs{\gamma}+\abs{\delta}=m, \abs{J}=\abs{K}=q,J,K\nearrow,K\neq I_0\ \mbox{if $(\alpha'',\beta')=0$}}c^{\alpha,\beta,\gamma,\delta}_{J,K}M_{d\ol z^J,d\ol z^K}\\
&\times\frac{1}{F(K)
+2<\abs{\lambda''},\alpha''>+2<\abs{\lambda'},\beta'>}\ol z^\alpha z^\beta\ol w^\gamma w^\delta+u(z), 
\end{split}
\end{equation} 
where $u\in{\rm Ker\,}L$.
\end{lem} 

\begin{proof} 
We recall that 
\begin{equation*}
L=\sum^n_{j=q+1}2\abs{\lambda_j}d\ol z_j^\wedge d\ol z_j^{\wedge,*}+
\sum^q_{j=1}2\abs{\lambda_j}d\ol z_j^{\wedge,*} d\ol z_j^\wedge
+\sum^n_{j=q+1}2\abs{\lambda_j}\ol z_j\frac{\pr}{\pr\ol z_j}+\sum^q_{j=1}2\abs{\lambda_j}z_j\frac{\pr}{\pr z_j}.
\end{equation*} 
For $M_{d\ol z^J,d\ol z^K}$, $\abs{J}=\abs{K}=q$, $J,K\nearrow$, we have 
\begin{align*} 
&\bigr(L(M_{d\ol z^J,d\ol z^K}\ol z^\alpha z^\beta\ol w^\gamma w^\delta)\bigr)d\ol z^J\\
&=\Bigr(\sum^n_{j=q+1}2\abs{\lambda_j}d\ol z_j^\wedge d\ol z_j^{\wedge,*}
+\sum^q_{j=1}2\abs{\lambda_j}d\ol z_j^{\wedge,*} d\ol z_j^\wedge\\
&\quad+\sum^n_{j=q+1}2\abs{\lambda_j}\ol z_j\frac{\pr}{\pr\ol z_j}+\sum^q_{j=1}2\abs{\lambda_j}z_j\frac{\pr}{\pr z_j}\Bigr)d\ol z^K\ol z^\alpha z^\beta\\
&=\Bigr(2\sum_{j\in K,q+1\leq j\leq n}\abs{\lambda_j}+2\sum_{j\notin K,1\leq j\leq q}\abs{\lambda_j}\\
&+2(<\abs{\lambda''},\alpha''>+<\abs{\lambda'},\beta'>)\Bigr)d\ol z^K\ol z^\alpha z^\beta\ol w^\gamma w^\delta\\
&=\Bigr(F(K)+2(<\abs{\lambda''},\alpha''>+<\abs{\lambda'},\beta'>)\Bigr)d\ol z^K\ol z^\alpha z^\beta\ol w^\gamma w^\delta
\end{align*} 
and 
\begin{equation*}
\bigr(L(M_{d\ol z^J,d\ol z^K}\ol z^\alpha z^\beta\ol w^\gamma w^\delta)\bigr)d\ol z^I=0
\end{equation*} 
if $I\neq J$. Thus, 
\begin{equation} \label{s3-e35}
\begin{split} 
&L(M_{d\ol z^J,d\ol z^K})\ol z^\alpha z^\beta\ol w^\gamma w^\delta\\
&=\Bigr(F(K)+2(<\abs{\lambda''},\alpha''>+<\abs{\lambda'},\beta'>)\Bigr)M_{d\ol z^J,d\ol z^K}\ol z^\alpha z^\beta\ol w^\gamma w^\delta. 
\end{split}
\end{equation} 
From \eqref{s3-e35}, the lemma follows.
\end{proof} 

It is not difficult to see that $b^0_0=\abs{\lambda_1}\abs{\lambda_2}\cdots\abs{\lambda_n}\pi^{-n}M_{d\ol z^{I_0},d\ol z^{I_0}}$ and 
\begin{equation} \label{s3-e35-1}
d\ol z_j^\wedge d\ol z_t^{\wedge,*}b^0_0=
\abs{\lambda_1}\abs{\lambda_2}\cdots\abs{\lambda_n}\pi^{-n}M_{d\ol z^{I_0},d\ol z_j^\wedge d\ol z_t^{\wedge,*}d\ol z^{I_0}},
\end{equation}
where $q+1\leq j\leq n$, $1\leq t\leq q$. Combining this with \eqref{s3-e29} and Lemma~\ref{s3-l1}, 
we get the following 

\begin{prop} \label{s3-p2}
We have that 
\begin{equation} \label{s3-e36} 
\begin{split} 
b^1_0(z,0)=&\sum_{q+1\leq j\leq n,1\leq t\leq q,1\leq s\leq n}\frac{1}{\abs{\lambda_j}+\abs{\lambda_t}+\abs{\lambda_s}\chi_1(s)}\frac{\pr^3\phi}{\pr\ol z_j\pr z_t\pr z_s}(0)z_sd\ol z_t^{\wedge,*}d\ol z_j^\wedge b^0_0\\
&+\sum_{q+1\leq j\leq n,1\leq t\leq q,1\leq s\leq n}\frac{1}{\abs{\lambda_j}+\abs{\lambda_t}+\abs{\lambda_s}\chi_2(s)}\frac{\pr^3\phi}{\pr\ol z_j\pr z_t\pr\ol z_s}(0)\ol z_sd\ol z_t^{\wedge,*}d\ol z_j^\wedge b^0_0\\
&-\sum_{q+1\leq s\leq n,1\leq j\leq q}\frac{1}{\abs{\lambda_j}+\abs{\lambda_s}}\frac{\pr^3\phi}{\pr\ol z_j\pr z_j\pr\ol z_s}(0)\ol z_sb^0_0 \\
&+\sum_{1\leq s\leq q,q+1\leq j\leq n}\frac{1}{\abs{\lambda_j}+\abs{\lambda_s}}\frac{\pr^3\phi}{\pr\ol z_j\pr z_j\pr z_s}(0) z_sb^0_0 \\
&+u(z),
\end{split}
\end{equation} 
where $u(z)\in{\rm Ker\,}L$.
\end{prop}

Now, we compute $b^1_0(0,w)$. From \eqref{s3-e36}, we can compute the last two terms of the right side of \eqref{s3-e30}: 
\begin{equation} \label{s3-e37} 
\begin{split}
&\sum^n_{s=q+1}2\abs{\lambda_s}\ol w_s\frac{\pr b^1_0(z,0)}{\pr\ol z_s}+\sum^q_{s=1}2\abs{\lambda_s}w_s\frac{\pr b^1_0(z,0)}{\pr z_s}\\ 
&=\sum_{q+1\leq j\leq n,1\leq t\leq q,1\leq s\leq q}\frac{2\abs{\lambda_s}}{\abs{\lambda_j}+\abs{\lambda_t}+\abs{\lambda_s}}\frac{\pr^3\phi}{\pr\ol z_j\pr z_t\pr z_s}(0)w_sd\ol z_t^{\wedge,*}d\ol z_j^\wedge b^0_0\\
&+\sum_{q+1\leq j\leq n,1\leq t\leq q,q+1\leq s\leq n}\frac{2\abs{\lambda_s}}{\abs{\lambda_j}+\abs{\lambda_t}+\abs{\lambda_s}}\frac{\pr^3\phi}{\pr\ol z_j\pr z_t\pr\ol z_s}(0)\ol w_sd\ol z_t^{\wedge,*}d\ol z_j^\wedge b^0_0\\
&-\sum_{q+1\leq s\leq n,1\leq j\leq q}\frac{2\abs{\lambda_s}}{\abs{\lambda_j}+\abs{\lambda_s}}\frac{\pr^3\phi}{\pr\ol z_j\pr z_j\pr\ol z_s}(0)\ol w_sb^0_0\\
&+\sum_{1\leq s\leq q,q+1\leq j\leq n}\frac{2\abs{\lambda_s}}{\abs{\lambda_j}+\abs{\lambda_s}}\frac{\pr^3\phi}{\pr\ol z_j\pr z_j\pr z_s}(0)w_sb^0_0.
\end{split}
\end{equation} 
Combining this with \eqref{s3-e30}, we obtain 
\begin{equation} \label{s3-e38} 
\begin{split}
&Lb^1_0(0,w)=\sum_{q+1\leq j\leq n,1\leq t\leq q,1\leq s\leq q}\frac{2\abs{\lambda_s}}{\abs{\lambda_j}+\abs{\lambda_t}+\abs{\lambda_s}}\frac{\pr^3\phi}{\pr\ol z_j\pr z_t\pr z_s}(0)w_sd\ol z_t^{\wedge,*}d\ol z_j^\wedge b^0_0\\
&+\sum_{q+1\leq j\leq n,1\leq t\leq q,q+1\leq s\leq n}\frac{2\abs{\lambda_s}}{\abs{\lambda_j}+\abs{\lambda_t}+\abs{\lambda_s}}\frac{\pr^3\phi}{\pr\ol z_j\pr z_t\pr\ol z_s}(0)\ol w_sd\ol z_t^{\wedge,*}d\ol z_j^\wedge b^0_0.
\end{split}
\end{equation} 
From this, \eqref{s3-e35-1} and Lemma~\ref{s3-l1}, we get 

\begin{prop} \label{s3-p3}
We have that 
\begin{equation} \label{s3-e39} 
\begin{split} 
&b^1_0(0,w)=\sum_{q+1\leq j\leq n,1\leq t\leq q,1\leq s\leq q}\frac{\abs{\lambda_s}}{
(\abs{\lambda_j}+\abs{\lambda_t})(\abs{\lambda_j}+\abs{\lambda_t}+\abs{\lambda_s})}\\
&\times\frac{\pr^3\phi}{\pr\ol z_j\pr z_t\pr z_s}(0)w_sd\ol z_t^{\wedge,*}d\ol z_j^\wedge b^0_0\\
&+\sum_{q+1\leq j\leq n,1\leq t\leq q,q+1\leq s\leq n}\frac{\abs{\lambda_s}}{(\abs{\lambda_j}+\abs{\lambda_t})(\abs{\lambda_j}+\abs{\lambda_t}+\abs{\lambda_s})}\\
&\times\frac{\pr^3\phi}{\pr\ol z_j\pr z_t\pr\ol z_s}(0)\ol w_sd\ol z_t^{\wedge,*}d\ol z_j^\wedge b^0_0+v(w),
\end{split}
\end{equation} 
where $v(w)\in{\rm Ker\,}L$.
\end{prop} 

In view of Proposition~\ref{s3-p2}, we know that to prove Theorem~\ref{s3-t1}, we only need to compute $u(z)$, where $u(z)$ is as in \eqref{s3-e36}. Now, we compute $u(z)$. Note that $u(z)\in{\rm Ker\,}L$. From Lemma~\ref{s3-l1}, we may write 
\begin{equation} \label{s3-e40} 
\begin{split}
u(z)&=\sum_{q+1\leq s\leq n,\abs{J}=q,J\nearrow,J\neq I_0}c^s_{J,I_0}M_{d\ol z^J,d\ol z^{I_0}}z_s\\
&+\sum_{1\leq s\leq q,\abs{J}=q,J\nearrow,J\neq I_0}c^s_{J,I_0}M_{d\ol z^J,d\ol z^{I_0}}\ol z_s\\
&+\sum_{q+1\leq s\leq n}z_sc^sM_{d\ol z^{I_0},d\ol z^{I_0}}+\sum_{1\leq s\leq q}\ol z_sc^sM_{d\ol z^{I_0},d\ol z^{I_0}}, 
\end{split}
\end{equation} 
where $c^s_{J,I_0}, c^s\in\Complex$, for all $s=1,\ldots,n$, $\abs{J}=q$, $J\nearrow$, $J\neq I_0$. Let $u^*(z)$ be the adjoint of $u(z)$ with respect to $(\ |\ )$ in the space $\mathscr L(\Lambda^{0,q}T^*(\Complex^n),\Lambda^{0,q}T^*(\Complex^n))$. We can check that 
\begin{equation} \label{s3-e41} 
\begin{split}
u^*(z)&=\sum_{q+1\leq s\leq n,\abs{J}=q,J\nearrow,J\neq I_0}\ol{c^s_{J,I_0}}M_{d\ol z^{I_0},d\ol z^J}\ol z_s\\
&+\sum_{1\leq s\leq q,\abs{J}=q,J\nearrow,J\neq I_0}\ol{c^s_{J,I_0}}M_{d\ol z^{I_0},d\ol z^J} z_s\\
&+\sum_{1\leq s\leq q}\ol{c^s}z_sM_{d\ol z^{I_0},d\ol z^{I_0}}+\sum_{q+1\leq s\leq n}\ol {c^s}\ol z_sM_{d\ol z^{I_0},d\ol z^{I_0}}.
\end{split}
\end{equation}
We notice that the Bergman projection $\Pi^{(q)}_k$ is self-adjoint. From this observation, we deduce that 
\begin{equation} \label{s3-e42}
(b^1_0(w, 0))^*=b^1_0(0, w),
\end{equation} 
where $(b^1_0(w, 0))^*$ is the adjoint of $b^1_0(w, 0)$ with respect to the inner product $(\ |\ )$ in the space $\mathscr L(\Lambda^{0,q}T^*(\Complex^n),\Lambda^{0,q}T^*(\Complex^n))$. From \eqref{s3-e36} and \eqref{s3-e41}  and recall that ${\rm Ker\,}L$ is given by \eqref{s3-e33}, we deduce that 
\begin{equation} \label{s3-e43} 
\begin{split} 
(b^1_0(w, 0))^*&=\sum_{q+1\leq s\leq n,\abs{J}=q,J\nearrow,J\neq I_0}\ol{c^s_{J,I_0}}M_{d\ol z^{I_0},d\ol z^J}\ol w_s\\
&+\sum_{1\leq s\leq q,\abs{J}=q,J\nearrow,J\neq I_0}\ol{c^s_{J,I_0}}M_{d\ol z^{I_0},d\ol z^J} w_s+r(w),
\end{split}
\end{equation} 
where $r(w)\in{\rm Ker\,}L$. From \eqref{s3-e39} and \eqref{s3-e35-1}, we have 
\begin{equation} \label{s3-e44} 
\begin{split} 
&b^1_0(0,w)=\abs{\lambda_1}\cdots\abs{\lambda_n}\pi^{-n}\times\\
&\Bigr(\sum_{q+1\leq j\leq n,1\leq t\leq q,1\leq s\leq q}\frac{\abs{\lambda_s}}{
(\abs{\lambda_j}+\abs{\lambda_t})(\abs{\lambda_j}+\abs{\lambda_t}+\abs{\lambda_s})}\\
&\times\frac{\pr^3\phi}{\pr\ol z_j\pr z_t\pr z_s}(0)w_sM_{d\ol z^{I_0},d\ol z_t^{\wedge,*}d\ol z_j^\wedge d\ol z^{I_0}}\\
&+\sum_{q+1\leq j\leq n,1\leq t\leq q,q+1\leq s\leq n}\frac{\abs{\lambda_s}}{(\abs{\lambda_j}+\abs{\lambda_t})(\abs{\lambda_j}+\abs{\lambda_t}+\abs{\lambda_s})}\\
&\times\frac{\pr^3\phi}{\pr\ol z_j\pr z_t\pr\ol z_s}(0)\ol w_sM_{d\ol z^{I_0},d\ol z_t^{\wedge,*}d\ol z_j^\wedge d\ol z^{I_0}}\Bigr)+v(w),
\end{split}
\end{equation} 
where $v(w)\in{\rm Ker\,}L$. From \eqref{s3-e42}, \eqref{s3-e43} and \eqref{s3-e44}, we get 
\begin{align}
&c^s_{J,I_0}=\abs{\lambda_1}\cdots\abs{\lambda_n}\pi^{-n}\times\frac{\abs{\lambda_s}}{
(\abs{\lambda_j}+\abs{\lambda_t})(\abs{\lambda_j}+\abs{\lambda_t}+\abs{\lambda_s})}\times\frac{\pr^3\phi}{\pr z_j\pr\ol z_t\pr\ol z_s}(0)\nonumber\\
&\mbox{if $J=d\ol z_t^{\wedge,*}d\ol z_j^\wedge d\ol z^{I_0}$, $q+1\leq j\leq n$, $1\leq t\leq q$, and $s=1,\ldots,q$},\label{s3-e45}\\
&c^s_{J,I_0}=\abs{\lambda_1}\cdots\abs{\lambda_n}\pi^{-n}\times\frac{\abs{\lambda_s}}{
(\abs{\lambda_j}+\abs{\lambda_t})(\abs{\lambda_j}+\abs{\lambda_t}+\abs{\lambda_s})}\times\frac{\pr^3\phi}{\pr z_j\pr\ol z_t\pr z_s}(0)\nonumber\\
&\mbox{if $J=d\ol z_t^{\wedge,*}d\ol z_j^\wedge d\ol z^{I_0}$, $q+1\leq j\leq n$, $1\leq t\leq q$, and $s=q+1,\ldots,n$},\label{s3-e46}\\
&c^s_{J,I_0}=0\ \ \mbox{otherwise}.\label{s3-e46}
\end{align}
Combining above with \eqref{s3-e40} and \eqref{s3-e35-1}, we obtain 
\begin{equation} \label{s3-e47}
\begin{split} 
b^1_0(z,0)=&\sum_{q+1\leq j\leq n,1\leq t\leq q,1\leq s\leq n}\frac{1}{\abs{\lambda_j}+\abs{\lambda_t}+\abs{\lambda_s}\chi_1(s)}\frac{\pr^3\phi}{\pr\ol z_j\pr z_t\pr z_s}(0)z_sd\ol z_t^{\wedge,*}d\ol z_j^\wedge b^0_0\\
&+\sum_{q+1\leq j\leq n,1\leq t\leq q,1\leq s\leq n}\frac{1}{\abs{\lambda_j}+\abs{\lambda_t}+\abs{\lambda_s}\chi_2(s)}\frac{\pr^3\phi}{\pr\ol z_j\pr z_t\pr\ol z_s}(0)\ol z_sd\ol z_t^{\wedge,*}d\ol z_j^\wedge b^0_0\\
&-\sum_{q+1\leq s\leq n,1\leq j\leq q}\frac{1}{\abs{\lambda_j}+\abs{\lambda_s}}\frac{\pr^3\phi}{\pr\ol z_j\pr z_j\pr\ol z_s}(0)\ol z_sb^0_0 \\
&+\sum_{1\leq s\leq q,q+1\leq j\leq n}\frac{1}{\abs{\lambda_j}+\abs{\lambda_s}}\frac{\pr^3\phi}{\pr\ol z_j\pr z_j\pr z_s}(0) z_sb^0_0 \\
&+\sum_{q+1\leq j\leq n, 1\leq t\leq q, q+1\leq s\leq n}\frac{\abs{\lambda_s}}{(\abs{\lambda_j}+\abs{\lambda_t})(\abs{\lambda_j}+\abs{\lambda_t}+\abs{\lambda_s})}\\
&\times\frac{\pr^3\phi}{\pr z_j\pr\ol z_t\pr z_s}(0)z_sb^0_0d\ol z_j^{\wedge,*}d\ol z_t^\wedge\\
&+\sum_{q+1\leq j\leq n, 1\leq t\leq q, 1\leq s\leq q}\frac{\abs{\lambda_s}}{(\abs{\lambda_j}+\abs{\lambda_t})(\abs{\lambda_j}+\abs{\lambda_t}+\abs{\lambda_s})}\\
&\times\frac{\pr^3\phi}{\pr z_j\pr\ol z_t\pr\ol z_s}(0)\ol z_sb^0_0d\ol z_j^{\wedge,*}d\ol z_t^\wedge \\
&+\sum_{q+1\leq s\leq n}z_sc^sM_{d\ol z^{I_0},d\ol z^{I_0}}+\sum_{1\leq s\leq q}\ol z_sc^sM_{d\ol z^{I_0},d\ol z^{I_0}}.
\end{split} 
\end{equation} 
Now, to complete the proof of Theorem~\ref{s3-t1}, we only need to know $c^s$, $s=1,\ldots,n$.  From \eqref{s3-e47}, we know that 
\begin{equation} \label{s3-e48} 
\begin{split}
&{\rm Tr\,}b^1_0(z,0)\\
&=-\abs{\lambda_1}\cdots\abs{\lambda_n}\pi^{-n}\Bigr(\sum_{q+1\leq s\leq n,1\leq j\leq q}\frac{1}{\abs{\lambda_j}+\abs{\lambda_s}}\frac{\pr^3\phi}{\pr\ol z_j\pr z_j\pr\ol z_s}(0)\ol z_s\Bigr) \\
&+\abs{\lambda_1}\cdots\abs{\lambda_n}\pi^{-n}\Bigr(\sum_{1\leq s\leq q,q+1\leq j\leq n}\frac{1}{\abs{\lambda_j}+\abs{\lambda_s}}\frac{\pr^3\phi}{\pr\ol z_j\pr z_j\pr z_s}(0) z_s\Bigr)\\ 
&+\sum_{q+1\leq s\leq n}z_sc^s+\sum_{1\leq s\leq q}\ol z_sc^s.
\end{split}
\end{equation} 
From Theorem~\ref{s1-t2}, we know that ${\rm Tr\,}b_0(z, z)=\pi^{-n}(-1)^q{\rm det\,}\left(\frac{\pr^2\phi}{\pr\ol z_j\pr z_k}(z)\right)^n_{j,k=1}$. From this, we can compute 
\begin{align} \label{s3-e49} 
&\frac{\pr}{\pr z_s}{\rm Tr\,}b_0(z, z)|_{z=0}=\pi^{-n}(-1)^q\frac{\pr}{\pr z_s}{\rm det\,}\left(\frac{\pr^2\phi}{\pr\ol z_j\pr z_k}(z)\right)^n_{j,k=1}|_{z=0}\nonumber\\
&=\pi^{-n}\abs{\lambda_1}\cdots\abs{\lambda_n}\Bigr(-\sum_{1\leq j\leq q}\frac{1}{\abs{\lambda_j}}\frac{\pr^3\phi}{\pr\ol z_j\pr z_j\pr z_s}(0)+\sum_{q+1\leq j\leq n}\frac{1}{\abs{\lambda_j}}\frac{\pr^3\phi}{\pr\ol z_j\pr z_j\pr z_s}(0)\Bigr).
\end{align} 
From \eqref{s3-e42}, we know that $\ol{{\rm Tr\,}b^1_0(z,0)}={\rm Tr\,}b^1_0(0,z)$. Note that $b^1_0(z,z)=b^1_0(z,0)+b^1_0(0,z)$. 
From this, we see that
\begin{align} \label{s3-e50} 
\frac{\pr}{\pr z_s}{\rm Tr\,}b_0(z, z)|_{z=0}&=
\frac{\pr}{\pr z_s}{\rm Tr\,}b^1_0(z,0)+\frac{\pr}{\pr z_s}{\rm Tr\,}b^1_0(0,z)\nonumber\\
&=\frac{\pr}{\pr z_s}{\rm Tr\,}b^1_0(z,0)+\ol{\frac{\pr}{\pr\ol z_s}{\rm Tr\,}b^1_0(z,0)}.
\end{align} 
Thus, if $q+1\leq s\leq n$, from \eqref{s3-e48} and \eqref{s3-e50}, we can check that 
\begin{equation} \label{s3-e51} 
\begin{split}
&\frac{\pr}{\pr z_s}{\rm Tr\,}b_0(z, z)|_{z=0}\\
&=-\abs{\lambda_1}\cdots\abs{\lambda_n}\pi^{-n}\Bigr(\sum_{1\leq j\leq q}\frac{1}{\abs{\lambda_j}+\abs{\lambda_s}}\frac{\pr^3\phi}{\pr\ol z_j\pr z_j\pr z_s}(0)\Bigr)+c^s. 
\end{split}
\end{equation} 
From \eqref{s3-e49} and \eqref{s3-e51}, we can compute 
\begin{equation} \label{s3-e52} 
\begin{split}
&c_s=\pi^{-n}\abs{\lambda_1}\cdots\abs{\lambda_n}\Bigr(-\sum_{1\leq j\leq q}\frac{1}{\abs{\lambda_j}}\frac{\pr^3\phi}{\pr\ol z_j\pr z_j\pr z_s}(0)+\sum_{q+1\leq j\leq n}\frac{1}{\abs{\lambda_j}}\frac{\pr^3\phi}{\pr\ol z_j\pr z_j\pr z_s}(0)\Bigr)\\
&+\pi^{-n}\abs{\lambda_1}\cdots\abs{\lambda_n}\Bigr(\sum_{1\leq j\leq q}\frac{1}{\abs{\lambda_j}+\abs{\lambda_s}}\frac{\pr^3\phi}{\pr\ol z_j\pr z_j\pr z_s}(0)\Bigr)\\
&=\pi^{-n}\abs{\lambda_1}\cdots\abs{\lambda_n}\\
&\times\Bigr(-\sum_{1\leq j\leq q}
\frac{\abs{\lambda_s}}{\abs{\lambda_j}(\abs{\lambda_j}+\abs{\lambda_s})}\frac{\pr^3\phi}{\pr\ol z_j\pr z_j\pr z_s}(0)
+\sum_{q+1\leq j\leq n}\frac{1}{\abs{\lambda_j}}\frac{\pr^3\phi}{\pr\ol z_j\pr z_j\pr z_s}(0)\Bigr),
\end{split}
\end{equation} 
$q+1\leq s\leq n$. Similarly, we can repeat the procedure above and get
\begin{equation} \label{s3-e53} 
\begin{split}
&c_s=\pi^{-n}\abs{\lambda_1}\cdots\abs{\lambda_n}\\
&\times\Bigr(-\sum_{1\leq j\leq q}\frac{1}{\abs{\lambda_j}}\frac{\pr^3\phi}{\pr\ol z_j\pr z_j\pr\ol z_s}(0)+\sum_{q+1\leq j\leq n}
\frac{\abs{\lambda_s}}{\abs{\lambda_j}(\abs{\lambda_j}+\abs{\lambda_s})}\frac{\pr^3\phi}{\pr\ol z_j\pr z_j\pr\ol z_s}(0)\Bigr),
\end{split}
\end{equation}
if $1\leq s\leq q$. Combining \eqref{s3-e53}, \eqref{s3-e52} with \eqref{s3-e47}, Theorem~\ref{s3-t1} follows.

\subsection{The second order of the Taylor expansion of $b_0(z,w)$ at $z=w$ and the $b_1$ term} 

The second goal of this section is to prove the following  

\begin{thm} \label{s3-t2}
Put 
\begin{equation} \label{s3-e54}
{\rm Tr\,} b^1_0(z,0)=\pi^{-n}\abs{\lambda_1}\cdots\abs{\lambda_n}\Bigr(\sum^n_{s=1}(a_sz_s+b_s\ol z_s)\Bigr).
\end{equation}
Then for $b^2_0(z,w)$ in \eqref{s3-e16}, we have 
\begin{equation} \label{s3-e54-1}
\begin{split}
b^2_0(z,0)=&\sum_{q+1\leq j\leq n,1\leq k\leq q, 1\leq s\leq n}\frac{1}{\abs{\lambda_s}(\abs{\lambda_j}+\abs{\lambda_k}+\abs{\lambda_s}\chi_1(s))}\\
&\times\abs{\frac{\pr^3\phi}{\pr\ol z_j\pr z_k\pr z_s}(0)}^2\abs{z_s}^2b^0_0\\
&+\sum_{q+1\leq j\leq n,1\leq k\leq q, 1\leq s\leq n}\frac{1}{\abs{\lambda_s}(\abs{\lambda_j}+\abs{\lambda_k}+\abs{\lambda_s}\chi_2(s))}\\
&\times\abs{\frac{\pr^3\phi}{\pr\ol z_j\pr z_k\pr\ol z_s}(0)}^2\abs{z_s}^2b^0_0\\
&+\sum_{q+1\leq j\leq n, 1\leq k\leq q, q+1\leq s\leq n}\frac{\abs{\lambda_s}}{(\abs{\lambda_j}+\abs{\lambda_k})(\abs{\lambda_j}+\abs{\lambda_k}+\abs{\lambda_s})^2}\\
&\times\abs{\frac{\pr^3\phi}{\pr\ol z_j\pr z_k\pr\ol z_s}(0)}^2\abs{z_s}^2d\ol z_k^{\wedge,*}d\ol z_j^\wedge b^0_0d\ol z_j^{\wedge,*}d\ol z_k^\wedge\\
&+\sum_{q+1\leq j\leq n, 1\leq k\leq q, 1\leq s\leq q}\frac{\abs{\lambda_s}}{(\abs{\lambda_j}+\abs{\lambda_k})(\abs{\lambda_j}+\abs{\lambda_k}+\abs{\lambda_s})^2}\\
&\times\abs{\frac{\pr^3\phi}{\pr\ol z_j\pr z_k\pr z_s}(0)}^2\abs{z_s}^2d\ol z_k^{\wedge,*}d\ol z_j^\wedge b^0_0d\ol z_j^{\wedge,*}d\ol z_k^\wedge \\ 
&-\sum_{q+1\leq u\leq n,1\leq s\leq n}\frac{1}{\abs{\lambda_u}+\abs{\lambda_s}\chi_1(s)}b_s\frac{\pr^3\phi}{\pr\ol z_u\pr z_u\pr z_s}(0)\abs{z_u}^2b^0_0\\
&+\sum_{1\leq u\leq q,1\leq s\leq n}\frac{1}{\abs{\lambda_u}
+\abs{\lambda_s}\chi_2(s)}a_s\frac{\pr^3\phi}{\pr\ol z_u\pr z_u\pr \ol z_s}(0)\abs{z_u}^2b^0_0\\
&-\sum_{1\leq j\leq q,q+1\leq s\leq n}\frac{1}{\abs{\lambda_j}+\abs{\lambda_s}}
a_s\frac{\pr^3\phi}{\pr\ol z_j\pr z_j\pr \ol z_s}(0)\abs{z_s}^2b^0_0\\
&+\sum_{q+1\leq j\leq n,1\leq s\leq q}\frac{1}{\abs{\lambda_j}+\abs{\lambda_s}}
b_s\frac{\pr^3\phi}{\pr\ol z_j\pr z_j\pr z_s}(0)\abs{z_s}^2b^0_0\\
&-\sum_{1\leq j\leq q,1\leq k\leq n}\frac{\pr^4(-i\psi(z,0)+\phi)}{\pr\ol z_j\pr z_j\pr\ol z_k\pr z_k}(0)\frac{\abs{z_k}^2}{2\abs{\lambda_k}}b^0_0\\
&+\sum_{q+1\leq j\leq n,1\leq k\leq n}\frac{\pr^4(i\psi(z,0)+\phi)}{\pr\ol z_j\pr z_j\pr\ol z_k\pr z_k}(0)\frac{\abs{z_k}^2}{2\abs{\lambda_k}}b^0_0
+r(z)+h(z)
\end{split}
\end{equation}
where ${\rm Tr\,}h=0$ and $\frac{\pr^2r}{\pr\ol z_j\pr z_j}=0$, $j=1,\ldots,n$. 

Furthermore, we have 
\begin{equation} \label{s3-e54-2}
\begin{split}
b_1(0,0)=&\sum_{q+1\leq j\leq n, 1\leq k\leq q, 1\leq s\leq q}\frac{\abs{\lambda_s}}{2(\abs{\lambda_j}+\abs{\lambda_k})^2
(\abs{\lambda_j}+\abs{\lambda_k}+\abs{\lambda_s})^2}\\
&\times\abs{\frac{\pr^3\phi}{\pr\ol z_j\pr z_k\pr z_s}(0)}^2d\ol z_k^{\wedge,*}d\ol z_j^\wedge b^0_0d\ol z_j^{\wedge,*}d\ol z_k^\wedge\\
&+\sum_{q+1\leq j\leq n, 1\leq k\leq q,q+1\leq s\leq n}\frac{\abs{\lambda_s}}{2(\abs{\lambda_j}+\abs{\lambda_k})^2
(\abs{\lambda_j}+\abs{\lambda_k}+\abs{\lambda_s})^2}\\
&\times\abs{\frac{\pr^3\phi}{\pr\ol z_j\pr z_k\pr\ol z_s}(0)}^2d\ol z_k^{\wedge,*}d\ol z_j^\wedge b^0_0d\ol z_j^{\wedge,*}d\ol z_k^\wedge\\
&+cb^0_0+R,
\end{split}
\end{equation} 
where $c\in\Complex$, ${\rm Tr\,}R=0$.
\end{thm}

Take $w=0$ in \eqref{s3-e22}, we get 
\begin{equation} \label{s3-e55}
\begin{split}
Lb^2_0(z,0)=&-\Bigr(\sum^n_{j,t=1}
\frac{\pr}{\pr\ol z_j}(-i\frac{\pr\psi_4(z,0)}{\pr z_t}+\frac{\pr\phi_4(z)}{\pr z_t})d\ol z_j^\wedge d\ol z_t^{\wedge,*}\\
&-\sum^n_{j,t=1}
\frac{\pr}{\pr z_t}(i\frac{\pr\psi_4(z,0)}{\pr\ol z_j}+\frac{\pr\phi_4(z)}{\pr\ol z_j})d\ol z_t^{\wedge,*} d\ol z_j^\wedge\Bigr)b^0_0\\
&-\Bigr(\sum^n_{j,t=1}
\frac{\pr}{\pr\ol z_j}(-i\frac{\pr\psi_3(z,0)}{\pr z_t}+\frac{\pr\phi_3(z)}{\pr z_t})d\ol z_j^\wedge d\ol z_t^{\wedge,*}\\
&-\sum^n_{j,t=1}
\frac{\pr}{\pr z_t}(i\frac{\pr\psi_3(z,0)}{\pr\ol z_j}+\frac{\pr\phi_3(z)}{\pr\ol z_j})d\ol z_t^{\wedge,*}d\ol z_j^\wedge\Bigr)b^1_0(z,0) \\
&-\sum^n_{j=1}(-i\frac{\pr\psi_3(z,0)}{\pr z_j}+\frac{\pr\phi_3(z)}{\pr z_j})\frac{\pr b^1_0(z,0)}{\pr\ol z_j}\\
&+\sum^n_{j=1}(i\frac{\pr\psi_3(z,0)}{\pr\ol z_j}+\frac{\pr\phi_3(z)}{\pr\ol z_j})\frac{\pr b^1_0(z,0)}{\pr z_j}\\
&=I+II+III,
\end{split}
\end{equation} 
where 
\begin{equation*} 
\begin{split}
&I=-\Bigr(\sum^n_{j,t=1}
\frac{\pr}{\pr\ol z_j}(-i\frac{\pr\psi_4(z,0)}{\pr z_t}+\frac{\pr\phi_4(z)}{\pr z_t})d\ol z_j^\wedge d\ol z_t^{\wedge,*}\\
&-\sum^n_{j,t=1}
\frac{\pr}{\pr z_t}(i\frac{\pr\psi_4(z,0)}{\pr\ol z_j}+\frac{\pr\phi_4(z)}{\pr\ol z_j})d\ol z_t^{\wedge,*} d\ol z_j^\wedge\Bigr)b^0_0, 
\end{split}
\end{equation*} 
\begin{equation*} 
\begin{split}
&II=-\Bigr(\sum^n_{j,t=1}
\frac{\pr}{\pr\ol z_j}(-i\frac{\pr\psi_3(z,0)}{\pr z_t}+\frac{\pr\phi_3(z)}{\pr z_t})d\ol z_j^\wedge d\ol z_t^{\wedge,*}\\
&-\sum^n_{j,t=1}
\frac{\pr}{\pr z_t}(i\frac{\pr\psi_3(z,0)}{\pr\ol z_j}+\frac{\pr\phi_3(z)}{\pr\ol z_j})d\ol z_t^{\wedge,*}d\ol z_j^\wedge\Bigr)b^1_0(z,0) 
\end{split} 
\end{equation*}
and 
\begin{equation*} 
\begin{split}
&III=-\sum^n_{j=1}(-i\frac{\pr\psi_3(z,0)}{\pr z_j}+\frac{\pr\phi_3(z)}{\pr z_j})\frac{\pr b^1_0(z,0)}{\pr\ol z_j}\\
&+\sum^n_{j=1}(i\frac{\pr\psi_3(z,0)}{\pr\ol z_j}+\frac{\pr\phi_3(z)}{\pr\ol z_j})\frac{\pr b^1_0(z,0)}{\pr z_j}
\end{split} 
\end{equation*}

First, we deal with the sum $I$. From the form of $b^0_0$, we can check that 
\begin{equation} \label{s3-e56}
\begin{split}
&-\sum^n_{j,t=1}
\frac{\pr}{\pr\ol z_j}(-i\frac{\pr\psi_4(z,0)}{\pr z_t}+\frac{\pr\phi_4(z)}{\pr z_t})d\ol z_j^\wedge d\ol z_t^{\wedge,*}\\
&+\sum^n_{j,t=1}
\frac{\pr}{\pr z_t}(i\frac{\pr\psi_4(z,0)}{\pr\ol z_j}+\frac{\pr\phi_4(z)}{\pr\ol z_j})d\ol z_t^{\wedge,*} d\ol z_j^\wedge b^0_0\\
&=-\sum_{1\leq j\leq q,1\leq k\leq n}\frac{\pr^4(-i\psi(z,0)+\phi)}{\pr\ol z_j\pr z_j\pr\ol z_k\pr z_k}(0)\abs{z_k}^2b^0_0\\
&+\sum_{q+1\leq j\leq n,1\leq k\leq n}\frac{\pr^4(i\psi(z,0)+\phi)}{\pr\ol z_j\pr z_j\pr\ol z_k\pr z_k}(0)\abs{z_k}^2b^0_0\\
&+h_1(z)+r_1(z),
\end{split}
\end{equation}
where ${\rm Tr\,}h_1=0$ and $\frac{\pr^2r_1}{\pr\ol z_j\pr z_j}=0$, $j=1,\ldots,n$. 

Now, we deal with the sum $II$. From\eqref{s3-e17}, we may write 
\begin{equation} \label{s3-e57}
\begin{split}
&b^1_0(z,0)=\sum_{1\leq s\leq n}(a_sz_s+b_s\ol z_s)b^0_0\\
&+\sum_{q+1\leq j\leq n,1\leq t\leq q,1\leq s\leq n}
(\alpha^s_{j,t}z_s+\beta^s_{j,t}\ol z_s)d\ol z_t^{\wedge,*}d\ol z_j^\wedge b^0_0\\
&+\sum_{q+1\leq j\leq n,1\leq t\leq q,1\leq s\leq n}
(\gamma^ s_{j,t}z_s+\delta^s_{j,t}\ol z_s)b^0_0d\ol z_j^{\wedge,*}d\ol z_t^\wedge,
\end{split}
\end{equation} 
where $a_s,b_s, \alpha^s_{j,t}, \beta^s_{j,t}, \gamma^s_{j,t}, \delta^s_{j,t}\in\Complex$, $s=1,\ldots,n$, $j=q+1,\ldots,n$, $t=1,\ldots,q$. We note that $a_s$, $b_s$, $s=1,\ldots,n$, are the same as in \eqref{s3-e54}.  For $j=q+1,\ldots,n$, $t=1,\ldots,q$ 
and $p, q=1,\ldots,n$, we can check that 
\begin{equation} \label{s3-e58}
{\rm Tr\,}(d\ol z_p^{\wedge,*}d\ol z_q^\wedge d\ol z_t^{\wedge,*}d\ol z_j^\wedge b^0_0)\neq0\ \ \mbox{if and only if 
$q=t$ and $p=j$}. 
\end{equation}
We notice that $d\ol z_j^{\wedge,*}d\ol z_t^\wedge d\ol z_t^{\wedge,*}d\ol z_j^\wedge b^0_0=b^0_0$. Similarly, 
For $j=q+1,\ldots,n$, $t=1,\ldots,q$ 
and $p, q=1,\ldots,n$, we can check that 
\begin{equation} \label{s3-e59} 
{\rm Tr\,}(d\ol z_p^{\wedge,*}d\ol z_q^\wedge b^0_0d\ol z_j^{\wedge,*}d\ol z_t^\wedge b^0_0)\neq0\ \ \mbox{if and only if  $p=t$ and $q=j$}. 
\end{equation} 
From \eqref{s3-e57}, \eqref{s3-e58} and \eqref{s3-e59}, it is straight forward to see that 
\begin{equation} \label{s3-e60}
\begin{split} 
&-\Bigr(\sum^n_{j,t=1}
\frac{\pr}{\pr\ol z_j}(-i\frac{\pr\psi_3(z,0)}{\pr z_t}+\frac{\pr\phi_3(z)}{\pr z_t})d\ol z_j^\wedge d\ol z_t^{\wedge,*}\\
&-\sum^n_{j,t=1}
\frac{\pr}{\pr z_t}(i\frac{\pr\psi_3(z,0)}{\pr\ol z_j}+\frac{\pr\phi_3(z)}{\pr\ol z_j})d\ol z_t^{\wedge,*}d\ol z_j^\wedge\Bigr)b^1_0(z,0)\\
&=-\sum_{1\leq j\leq q,1\leq s\leq n}(-i\frac{\pr^3\psi_3}{\pr z_s\pr z_j\pr\ol z_j}+\frac{\pr^3\phi_3}{\pr z_s\pr z_j\pr\ol z_j})\abs{z_s}^2b_sb^0_0\\
&-\sum_{1\leq j\leq q,1\leq s\leq n}(-i\frac{\pr^3\psi_3}{\pr\ol z_s\pr z_j\pr\ol z_j}+\frac{\pr^3\phi_3}{\pr\ol z_s\pr z_j\pr\ol z_j})\abs{z_s}^2a_sb^0_0\\
&+\sum_{q+1\leq j\leq n,1\leq s\leq n}(i\frac{\pr^3\psi_3}{\pr z_s\pr z_j\pr\ol z_j}+\frac{\pr^3\phi_3}{\pr z_s\pr z_j\pr\ol z_j})\abs{z_s}^2b_sb^0_0\\
&+\sum_{q+1\leq j\leq n,1\leq s\leq n}(i\frac{\pr^3\psi_3}{\pr\ol z_s\pr z_j\pr\ol z_j}+\frac{\pr^3\phi_3}{\pr\ol z_s\pr z_j\pr\ol z_j})\abs{z_s}^2a_sb^0_0\\
&+2\sum_{q+1\leq j\leq n,1\leq t\leq q,1\leq s\leq n}\frac{\pr^3\phi_3}{\pr z_s\pr\ol z_t\pr z_j}\beta^s_{j,t}\abs{z_s}^2b^0_0\\
&+2\sum_{q+1\leq j\leq n,1\leq t\leq q,1\leq s\leq n}\frac{\pr^3\phi_3}{\pr\ol z_s\pr\ol z_t\pr z_j}\alpha^s_{j,t}\abs{z_s}^2b^0_0\\
&+2\sum_{q+1\leq j\leq n,1\leq t\leq q,1\leq s\leq n}\frac{\pr^3\phi_3}{\pr z_s\pr z_t\pr\ol z_j}\delta^s_{j,t}\abs{z_s}^2
d\ol z_t^{\wedge,*}d\ol z_j^\wedge b^0_0d\ol z_j^{\wedge,*}d\ol z_t^\wedge\\
&+2\sum_{q+1\leq j\leq n,1\leq t\leq q,1\leq s\leq n}\frac{\pr^3\phi_3}{\pr\ol z_s\pr z_t\pr\ol z_j}\gamma^s_{j,t}\abs{z_s}^2
d\ol z_t^{\wedge,*}d\ol z_j^\wedge b^0_0d\ol z_j^{\wedge,*}d\ol z_t^\wedge\\
&+h_2(z)+r_2(z),
\end{split}
\end{equation} 
where ${\rm Tr\,}h_2=0$ and $\frac{\pr^2r_2}{\pr\ol z_j\pr z_j}=0$, $j=1,\ldots,n$. 

From \eqref{s2-e2} and \eqref{s3-e17}, it is straight forward to see that \eqref{s3-e60} becomes: 
\begin{equation} \label{s3-e61}
\begin{split} 
&-\Bigr(\sum^n_{j,t=1}
\frac{\pr}{\pr\ol z_j}(-i\frac{\pr\psi_3(z,0)}{\pr z_t}+\frac{\pr\phi_3(z)}{\pr z_t})d\ol z_j^\wedge d\ol z_t^{\wedge,*}\\
&-\sum^n_{j,t=1}
\frac{\pr}{\pr z_t}(i\frac{\pr\psi_3(z,0)}{\pr\ol z_j}+\frac{\pr\phi_3(z)}{\pr\ol z_j})d\ol z_t^{\wedge,*}d\ol z_j^\wedge\Bigr)b^1_0(z,0)\\
&=-\sum_{1\leq j\leq q,q+1\leq s\leq n}\frac{2\abs{\lambda_s}}{\abs{\lambda_j}+\abs{\lambda_s}}\frac{\pr^3\phi}{\pr\ol z_s\pr z_j\pr\ol z_j}(0)\abs{z_s}^2a_sb^0_0\\
&+\sum_{q+1\leq j\leq n,1\leq s\leq q}\frac{2\abs{\lambda_s}}{\abs{\lambda_j}+\abs{\lambda_s}}\frac{\pr^3\phi}{\pr z_s\pr z_j\pr\ol z_j}(0)\abs{z_s}^2b_sb^0_0\\ 
&+\sum_{q+1\leq j\leq n,1\leq t\leq q,1\leq s\leq n}\frac{2}{\abs{\lambda_j}+\abs{\lambda_t}+\abs{\lambda_s}\chi_2(s)}\abs{\frac{\pr^3\phi}{\pr\ol z_s\pr z_t\pr\ol z_j}(0)}^2\abs{z_s}^2b^0_0\\
&+\sum_{q+1\leq j\leq n,1\leq t\leq q,1\leq s\leq n}\frac{2}{\abs{\lambda_j}+\abs{\lambda_t}+\abs{\lambda_s}\chi_1(s)}\abs{\frac{\pr^3\phi}{\pr z_s\pr z_t\pr\ol z_j}(0)}^2\abs{z_s}^2b^0_0\\
&+\sum_{q+1\leq j\leq n,1\leq s,t\leq q}\frac{2\abs{\lambda_s}}{(\abs{\lambda_j}+\abs{\lambda_t})(\abs{\lambda_j}+\abs{\lambda_t}+\abs{\lambda_s})}\\
&\times\abs{\frac{\pr^3\phi}{\pr z_s\pr z_t\pr\ol z_j}(0)}^2\abs{z_s}^2
d\ol z_t^{\wedge,*}d\ol z_j^\wedge b^0_0d\ol z_j^{\wedge,*}d\ol z_t^\wedge\\
&+\sum_{q+1\leq j,s\leq n,1\leq t\leq q}\frac{2\abs{\lambda_s}}{(\abs{\lambda_j}+\abs{\lambda_t})(\abs{\lambda_j}+\abs{\lambda_t}+\abs{\lambda_s})}\\
&\times\abs{\frac{\pr^3\phi}{\pr\ol z_s\pr z_t\pr\ol z_j}(0)}^2\abs{z_s}^2
d\ol z_t^{\wedge,*}d\ol z_j^\wedge b^0_0d\ol z_j^{\wedge,*}d\ol z_t^\wedge\\
&+h_2(z)+r_2(z).
\end{split}
\end{equation} 

Now, we deal with the sum $III$.
From \eqref{s2-e2}, \eqref{s3-e17} and \eqref{s3-e57},  we can check that 
\begin{equation} \label{s3-e62}
\begin{split}
&-\sum^n_{j=1}(-i\frac{\pr\psi_3(z,0)}{\pr z_j}+\frac{\pr\phi_3(z)}{\pr z_j})\frac{\pr b^1_0(z,0)}{\pr\ol z_j}\\
&+\sum^n_{j=1}(i\frac{\pr\psi_3(z,0)}{\pr\ol z_j}+\frac{\pr\phi_3(z)}{\pr\ol z_j})\frac{\pr b^1_0(z,0)}{\pr z_j}\\
&=-\sum_{q+1\leq u\leq n,1\leq s\leq n}\frac{2\abs{\lambda_u}}{\abs{\lambda_u}+\abs{\lambda_s}\chi_1(s)}b_s\frac{\pr^3\phi}{\pr\ol z_u\pr z_u\pr z_s}(0)\abs{z_u}^2\\
&+\sum_{1\leq u\leq q,1\leq s\leq n}\frac{2\abs{\lambda_u}}{\abs{\lambda_u}+\abs{\lambda_s}\chi_2(s)}a_s\frac{\pr^3\phi}{\pr\ol z_u\pr z_u\pr\ol z_s}(0)\abs{z_u}^2\\
&+h_3(z)+r_3(z),
\end{split}
\end{equation} 
where ${\rm Tr\,}h_3=0$ and $\frac{\pr^2r_3}{\pr\ol z_j\pr z_j}=0$, $j=1,\ldots,n$. 

From \eqref{s3-e62}, \eqref{s3-e61} and \eqref{s3-e56}, \eqref{s3-e55} becomes 
\begin{equation} \label{s3-e63}
\begin{split}
&Lb^2_0(z,0)=-\sum_{1\leq j\leq q,1\leq k\leq n}\frac{\pr^4(-i\psi(z,0)+\phi)}{\pr\ol z_j\pr z_j\pr\ol z_k\pr z_k}(0)\abs{z_k}^2b^0_0\\
&+\sum_{q+1\leq j\leq n,1\leq k\leq n}\frac{\pr^4(i\psi(z,0)+\phi)}{\pr\ol z_j\pr z_j\pr\ol z_k\pr z_k}(0)\abs{z_k}^2b^0_0\\
&-\sum_{1\leq j\leq q,q+1\leq s\leq n}\frac{2\abs{\lambda_s}}{\abs{\lambda_j}+\abs{\lambda_s}}\frac{\pr^3\phi}{\pr\ol z_s\pr z_j\pr\ol z_j}(0)\abs{z_s}^2a_sb^0_0\\
&+\sum_{q+1\leq j\leq n,1\leq s\leq q}\frac{2\abs{\lambda_s}}{\abs{\lambda_j}+\abs{\lambda_s}}\frac{\pr^3\phi}{\pr z_s\pr z_j\pr\ol z_j}(0)\abs{z_s}^2b_sb^0_0\\ 
&+\sum_{q+1\leq j\leq n,1\leq t\leq q,1\leq s\leq n}\frac{2}{\abs{\lambda_j}+\abs{\lambda_t}+\abs{\lambda_s}\chi_2(s)}\abs{\frac{\pr^3\phi}{\pr\ol z_s\pr z_t\pr\ol z_j}(0)}^2\abs{z_s}^2b^0_0\\
&+\sum_{q+1\leq j\leq n,1\leq t\leq q,1\leq s\leq n}\frac{2}{\abs{\lambda_j}+\abs{\lambda_t}+\abs{\lambda_s}\chi_1(s)}\abs{\frac{\pr^3\phi}{\pr z_s\pr z_t\pr\ol z_j}(0)}^2\abs{z_s}^2b^0_0\\
&+\sum_{q+1\leq j\leq n,1\leq s,t\leq q}\frac{2\abs{\lambda_s}}{(\abs{\lambda_j}+\abs{\lambda_t})(\abs{\lambda_j}+\abs{\lambda_t}+\abs{\lambda_s})}\\
&\times\abs{\frac{\pr^3\phi}{\pr z_s\pr z_t\pr\ol z_j}(0)}^2\abs{z_s}^2
d\ol z_t^{\wedge,*}d\ol z_j^\wedge b^0_0d\ol z_j^{\wedge,*}d\ol z_t^\wedge\\
&+\sum_{q+1\leq j,s\leq n,1\leq t\leq q}\frac{2\abs{\lambda_s}}{(\abs{\lambda_j}+\abs{\lambda_t})(\abs{\lambda_j}+\abs{\lambda_t}+\abs{\lambda_s})}\\
&\times\abs{\frac{\pr^3\phi}{\pr\ol z_s\pr z_t\pr\ol z_j}(0)}^2\abs{z_s}^2
d\ol z_t^{\wedge,*}d\ol z_j^\wedge b^0_0d\ol z_j^{\wedge,*}d\ol z_t^\wedge\\
&-\sum_{q+1\leq u\leq n,1\leq s\leq n}\frac{2\abs{\lambda_u}}{\abs{\lambda_u}+\abs{\lambda_s}\chi_1(s)}b_s\frac{\pr^3\phi}{\pr\ol z_u\pr z_u\pr z_s}(0)\abs{z_u}^2\\
&+\sum_{1\leq u\leq q,1\leq s\leq n}\frac{2\abs{\lambda_u}}{\abs{\lambda_u}+\abs{\lambda_s}\chi_2(s)}a_s\frac{\pr^3\phi}{\pr\ol z_u\pr z_u\pr\ol z_s}(0)\abs{z_u}^2\\
&+h_4(z)+r_4(z),
\end{split}
\end{equation}
where ${\rm Tr\,}h_4=0$ and $\frac{\pr^2r_4}{\pr\ol z_j\pr z_j}=0$, $j=1,\ldots,n$. We notice that $d\ol z_t^{\wedge,*}d\ol z_j^\wedge b^0_0d\ol z_j^{\wedge,*}d\ol z_t^\wedge=M_{d\ol z_t^{\wedge,*}d\ol z_j^\wedge d\ol z^{I_0},d\ol z_t^{\wedge,*}d\ol z_j^\wedge d\ol z^{I_0}}$. From this observation and Lemma~\ref{s3-l1}, we get \eqref{s3-e54-1}. 

Now, to complete the proof of Theorem~\ref{s3-t2}, we only need to prove \eqref{s3-e54-2}. From \eqref{s3-e15} and \eqref{s2-e2}, we see that 
\begin{equation} \label{s3-e64} 
\Bigr(\sum^q_{j=1}2\abs{\lambda_j}d\ol z^\wedge_jd\ol z_j^{\wedge,*}+
\sum^n_{j=q+1}2\abs{\lambda_j}d\ol z^{\wedge,*}_jd\ol z_j^\wedge\Bigr)b_1(0,0)=\sum^n_{j=1}\frac{\pr^2b_0}{\pr z_j\pr\ol z_j}(0,0).
\end{equation} 
From \eqref{s3-e54-1}, we see that 
\begin{equation} \label{s3-e65} 
\begin{split}
&\sum^n_{j=1}\frac{\pr^2b_0}{\pr z_j\pr\ol z_j}(0,0)=\sum_{q+1\leq j\leq n, 1\leq k\leq q, 1\leq s\leq q}\frac{\abs{\lambda_s}}{(\abs{\lambda_j}+\abs{\lambda_k})
(\abs{\lambda_j}+\abs{\lambda_k}+\abs{\lambda_s})^2}\\
&\times\abs{\frac{\pr^3\phi}{\pr\ol z_j\pr z_k\pr z_s}(0)}^2d\ol z_k^{\wedge,*}d\ol z_j^\wedge b^0_0d\ol z_j^{\wedge,*}d\ol z_k^\wedge\\
&+\sum_{q+1\leq j\leq n, 1\leq k\leq q,q+1\leq s\leq n}\frac{\abs{\lambda_s}}{(\abs{\lambda_j}+\abs{\lambda_k})
(\abs{\lambda_j}+\abs{\lambda_k}+\abs{\lambda_s})^2}\\
&\times\abs{\frac{\pr^3\phi}{\pr\ol z_j\pr z_k\pr\ol z_s}(0)}^2d\ol z_k^{\wedge,*}d\ol z_j^\wedge b^0_0d\ol z_j^{\wedge,*}d\ol z_k^\wedge+\alpha b^0_0+f,
\end{split}
\end{equation} 
where $\alpha\in\Complex$ and ${\rm Tr\,}f=0$. Since we can solve \eqref{s3-e64}, we conclude that $\alpha=0$. Thus, \eqref{s3-e64} becomes: 
\begin{equation} \label{s3-e66} 
\begin{split} 
&\Bigr(\sum^q_{j=1}2\abs{\lambda_j}d\ol z^\wedge_jd\ol z_j^{\wedge,*}+
\sum^n_{j=q+1}2\abs{\lambda_j}d\ol z^{\wedge,*}_jd\ol z_j^\wedge\Bigr)b_1(0,0)\\
&=\sum_{q+1\leq j\leq n, 1\leq k\leq q, 1\leq s\leq q}\frac{\abs{\lambda_s}}{(\abs{\lambda_j}+\abs{\lambda_k})
(\abs{\lambda_j}+\abs{\lambda_k}+\abs{\lambda_s})^2}\\
&\times\abs{\frac{\pr^3\phi}{\pr\ol z_j\pr z_k\pr z_s}(0)}^2d\ol z_k^{\wedge,*}d\ol z_j^\wedge b^0_0d\ol z_j^{\wedge,*}d\ol z_k^\wedge\\
&+\sum_{q+1\leq j\leq n, 1\leq k\leq q,q+1\leq s\leq n}\frac{\abs{\lambda_s}}{(\abs{\lambda_j}+\abs{\lambda_k})
(\abs{\lambda_j}+\abs{\lambda_k}+\abs{\lambda_s})^2}\\
&\times\abs{\frac{\pr^3\phi}{\pr\ol z_j\pr z_k\pr\ol z_s}(0)}^2d\ol z_k^{\wedge,*}d\ol z_j^\wedge b^0_0d\ol z_j^{\wedge,*}d\ol z_k^\wedge+f.
\end{split}
\end{equation} 
Again, from Lemma~\ref{s3-l1}, we get \eqref{s3-e54-2}. Theorem~\ref{s3-t2} follows. 

\section{The trace of the $b_1$ term} 

As before, in this section, we assume that $\phi(z)=\sum^n_{j=1}\lambda_j\abs{z_j}^2+O(\abs{z}^3)$ near $0$ and that 
$\lambda_j<0$, $j=1,\ldots,q$, and $\lambda_j>0$, $j=q+1,\ldots,n$. We work in some neighborhood of $(0,0)$ and we shall use 
the same notaions as before. In view of \eqref{s3-e54-2}, we know that 
\begin{equation} \label{s4-e1}
\begin{split}
b_1(0,0)=&\sum_{q+1\leq j\leq n, 1\leq k\leq q, 1\leq s\leq q}\frac{\abs{\lambda_s}}{2(\abs{\lambda_j}+\abs{\lambda_k})^2
(\abs{\lambda_j}+\abs{\lambda_k}+\abs{\lambda_s})^2}\\
&\times\abs{\frac{\pr^3\phi}{\pr\ol z_j\pr z_k\pr z_s}(0)}^2d\ol z_k^{\wedge,*}d\ol z_j^\wedge b^0_0d\ol z_j^{\wedge,*}d\ol z_k^\wedge\\
&+\sum_{q+1\leq j\leq n, 1\leq k\leq q,q+1\leq s\leq n}\frac{\abs{\lambda_s}}{2(\abs{\lambda_j}+\abs{\lambda_k})^2
(\abs{\lambda_j}+\abs{\lambda_k}+\abs{\lambda_s})^2}\\
&\times\abs{\frac{\pr^3\phi}{\pr\ol z_j\pr z_k\pr\ol z_s}(0)}^2d\ol z_k^{\wedge,*}d\ol z_j^\wedge b^0_0d\ol z_j^{\wedge,*}d\ol z_k^\wedge\\
&+cb^0_0+R,
\end{split}
\end{equation} 
where $c\in\Complex$, ${\rm Tr\,}R=0$. The main goal of this section is to determine the constant $c$. We notice that the projection $\hat\Pi^{(q)}_k$ has the following property: 
\begin{equation} \label{s4-e2} 
\hat\Pi^{(q)}_k\circ\hat\Pi^{(q)}_k=\hat\Pi^{(q)}_k.
\end{equation} 
We recall that $\hat\Pi^{(q)}_k$ is given by \eqref{s3-e7}. From \eqref{s3-e9}, we have 
\begin{equation} \label{s4-e3} 
\hat\Pi^{(q)}_k\circ\hat\Pi^{(q)}_k(u,w)=\int e^{ik\psi(z,u)}b(u,z,k)e^{ik\psi(z,w)}b(z,w,k)dm(z)+F(u,w,k),
\end{equation} 
where $F(u,w,k)$ is negligible. Take $u=w=0$ in \eqref{s4-e3} and from \eqref{s4-e2} and \eqref{s3-e9}, we get 
\begin{equation} \label{s4-e4} 
\int e^{ik(\psi(0,z)+\psi(z,0))}b(0,z,k)b(z,0,k)dm(z)\sim k^nb_0(0,0)+k^{n-1}b_1(0,0)+\cdots.
\end{equation} 

We use $z_j=x_{2j-1}+ix_{2j}$, $j=1,\ldots,n$. We write $\pr_{x_j}$ to denote the operator $\frac{\pr}{\pr x_j}$, $j=1,\ldots,2n$. 
For  multi-indix $\alpha=(\alpha_1,\ldots,\alpha_{2n})$, $\alpha_j\in\mathbb N\bigcup\set{0}$, $j=1,\ldots,2n$. We write $\abs{\alpha}=N$ if $\sum^{2n}_{j=1}\alpha_j=N$ and we write $\pr^\alpha_x$ to denote the operator $\pr^{\alpha_1}_{x_1}\cdots\pr^{\alpha_{2n}}_{x_{2n}}$. We recall the stationary phase formula of H\"{o}rmander (see Theorem~7.7.5 in H\"{o}rmander ~\cite{Hor03}) 

\begin{thm} \label{s4-t1} 
Let $K\subset\Complex^n$ be a compact set, $X$ an open neighborhood of $K$ and $N$ a positive integer. If $u\in C^\infty_0(K)$, 
$f\in C^\infty(X)$ and ${\rm Im\,}f\geq0$ in $X$, ${\rm Im\,}f(x_0)=0$, $f'(x_0)=0$, ${\rm det\,}f''(x_0)\neq0$, $f'\neq0$ in $K\setminus\set{x_0}$ then 
\begin{equation} \label{s4-e5} 
\begin{split}
&\abs{\int u(z)e^{ikf(z)}dm-2^ne^{ikf(x_0)}{\rm det\,}\left(\frac{kf''(x_0)}{2\pi i}\right)^{-\frac{1}{2}}\sum_{j<N}k^{-j}L_ju} \\
&\quad\leq Ck^{-N}\sum_{\abs{\alpha}\leq 2N}\sup{\abs{\pr^\alpha_x u}},\ \ k>0,
\end{split}
\end{equation} 
where $C$ is bounded when $f$ stays in a bounded set in $C^\infty(X)$ and $\frac{\abs{x-x_0}}{\abs{f'(x)}}$ has a uniform bounded and 
\begin{equation} \label{s4-e6} 
L_ju=\sum_{\nu-\mu=j}\sum_{2\nu\geq 3\mu}i^{-j}2^{-\nu}<f''(x_0)^{-1}D,D>^\nu\frac{(g^\mu_{x_0}u)(x_0)}{\nu!\mu!}.
\end{equation} 
Here 
\begin{equation} \label{s4-e7} 
g_{x_0}(x)=f(x)-f(x_0)-\frac{1}{2}<f''(x_0)(x-x_0),x-x_0>
\end{equation} 
and $D=\left(\begin{array}[c]{ccc}
  &-i\pr_{x_1}  \\
  &\vdots\\
  &-i\pr_{x_{2n}}
\end{array}\right)$.
\end{thm} 

Now, we apply \eqref{s4-e5} to the left side of \eqref{s4-e4}. From \eqref{s2-e2}, we know that 
\begin{equation} \label{s4-e8} 
\psi(z,0)+\psi(0,z)=2i{\rm\, Im}\psi(z,0)=2i\sum^n_{j=1}\abs{\lambda_j}\abs{z_j}^2+O(\abs{z}^3).
\end{equation}
Since ${\rm Im\,}(\psi(0,z)+\psi(z,0))>0$ when $z\neq0$, we may assume that $b(0,z,k)$ has compact supports in some small neighborhood $K$ of $0\in\Complex^n$. From \eqref{s4-e5}, we have 
\begin{equation} \label{s4-e9}
\begin{split}
&\int e^{ik(\psi(0,z)+\psi(z,0))}b(0,z,k)b(z,0,k)dm(z)\\
&=\abs{\lambda_1}^{-1}\cdots\abs{\lambda_n}^{-1}\pi^nk^{-n}
\Bigr(b(0,0,k)^2+k^{-1}L_1(b(0,z,k)b(z,0,k))|_{z=0}+O(k^{2n-2})\Bigr).
\end{split}
\end{equation} 
We can check that
\begin{equation} \label{s4-e9-1} 
b(0,0,k)^2=k^{2n}b_0(0,0)^2+k^{2n-1}(b_0(0,0)b_1(0,0)+b_1(0,0)b_0(0,0))+O(k^{2n-2}).
\end{equation} 
The computation of the term $L_1(b(0,z,k)b(z,0,k))|_{z=0}$ is straight forward. We omit the process. We state our result 

\begin{prop} \label{s4-p1}
Under the notations above, we have 
\begin{equation} \label{s4-e10} 
\begin{split}
&k^{-2n}L_1(b(0,z,k)b(z,0,k))|_{z=0}+O(k^{-1})\\
&=\frac{1}{2}\sum^n_{j=1}\frac{1}{\abs{\lambda_j}}\Bigr(\frac{\pr^2b_0}{\pr z_j\pr\ol z_j}(z,0)|_{z=0}b_0(0,0)+
\frac{\pr^2 b_0}{\pr z_j\pr\ol z_j}(0,z)|_{z=0}b_0(0,0)\\
&+(\frac{\pr b_0}{\pr\ol z_j})(0,z)|_{z=0}\frac{\pr b_0}{\pr z_j}(z,0)|_{z=0}+\frac{\pr b_0}{\pr z_j}(0,z)|_{z=0}
\frac{\pr b_0}{\pr\ol z_j}(z,0)|_{z=0}\Bigr)\\
&-\frac{1}{4}\sum^n_{j,t=1}\frac{1}{\abs{\lambda_j}
\abs{\lambda_t}}\frac{\pr^4{\rm Im\,}\psi(z,0)}{\pr\ol z_j\pr z_j\pr\ol z_t\pr z_t}|_{z=0}b_0(0,0)^2\\
&-\frac{1}{2}\sum^n_{j,s=1}\frac{1}{\abs{\lambda_j}\abs{\lambda_s}}\frac{\pr^3{\rm Im\,}\psi(z,0)}{\pr\ol z_j\pr z_j\pr\ol z_s}|_{z=0}\Bigr(\frac{\pr b_0}{\pr z_s}(0,z)|_{z=0}+\frac{\pr b_0}{\pr z_s}(z,0)|_{z=0}\Bigr)b_0(0,0)\\
&-\frac{1}{2}\sum^n_{j,s=1}\frac{1}{\abs{\lambda_j}\abs{\lambda_s}}\frac{\pr^3{\rm Im\,}\psi(z,0)}{\pr\ol z_j\pr z_j\pr z_s}|_{z=0}\Bigr(\frac{\pr b_0}{\pr\ol z_s}(0,z)|_{z=0}+\frac{\pr b_0}{\pr\ol z_s}(z,0)|_{z=0}\Bigr)b_0(0,0)\\
&+\frac{1}{4}\sum^n_{j,t,s=1}\frac{1}{\abs{\lambda_j}\abs{\lambda_t}\abs{\lambda_s}}\abs{\frac{\pr^3{\rm Im\,}\psi(z,0)}{\pr z_j\pr z_t\pr\ol z_s}|_{z=0}}^2b_0(0,0)^2\\
&+\frac{1}{2}\sum^n_{j,t,s=1}\frac{1}{\abs{\lambda_j}\abs{\lambda_t}\abs{\lambda_s}}\frac{\pr^3{\rm Im\,}\psi(z,0)}{\pr\ol z_j\pr z_j\pr\ol z_s}|_{z=0}\frac{\pr^3{\rm Im\,}\psi(z,0)}{\pr\ol z_t\pr z_t\pr z_s}|_{z=0}b_0(0,0)^2\\
&+\frac{1}{12}\sum^n_{j,k,s=1}\frac{1}{\abs{\lambda_j}\abs{\lambda_t}\abs{\lambda_s}}\abs{\frac{\pr^3{\rm Im\,}\psi(z,0)}{\pr z_j\pr z_t\pr z_s}|_{z=0}}^2b_0(0,0)^2.
\end{split}
\end{equation}
\end{prop} 

From \eqref{s4-e4}, \eqref{s4-e9}, \eqref{s4-e9-1} and \eqref{s4-e10}, we get the following 

\begin{thm} \label{s4-t2} 
We have
\begin{equation} \label{s4-e11}
b_1(0,0)=\abs{\lambda_1}^{-1}\cdots\abs{\lambda_n}^{-1}\pi^n\Bigr(b_1(0,0)b_0(0,0)+b_0(0,0)b_1(0,0)+C_0\Bigr),
\end{equation} 
where $C_0$ denote the right side of \eqref{s4-e10}.
\end{thm} 
For $A, B\in\mathscr L(\Lambda^{0,q}T^*_0(\Complex^n),\Lambda^{0,q}T^*_0(\Complex^n))$, we write $A\bot B$ if 
$(Au\ |\ Bu)=0$ for all $u\in\Lambda^{0,q}T^*_0(\Complex^n)$.
For $F\in\mathscr L(\Lambda^{0,q}T^*_0(\Complex^n),\Lambda^{0,q}T^*_0(\Complex^n))$, we werit $\hat F$ to denote the component of $A$ in the direction $d\ol z_1^\wedge d\ol z_1^{\wedge,*}\cdots d\ol z_q^\wedge d\ol z_q^{\wedge,*}$. More precisely, if 
\begin{equation} \label{s4-e12}
F=\alpha d\ol z_1^\wedge d\ol z_1^{\wedge,*}\cdots d\ol z_q^\wedge d\ol z_q^{\wedge,*}+G,\ \ \alpha\in\Complex,\ \ 
G\bot d\ol z_1^\wedge d\ol z_1^{\wedge,*}\cdots d\ol z_q^\wedge d\ol z_q^{\wedge,*},
\end{equation}
then $\hat F=\alpha d\ol z_1^\wedge d\ol z_1^{\wedge,*}\cdots d\ol z_q^\wedge d\ol z_q^{\wedge,*}$. Now, we are ready to compute the constant $c$. We recall that $c$ is given by \eqref{s4-e1}. From \eqref{s4-e1}, we know that 
\begin{equation} \label{s4-e13}
\hat b_1(0,0)=cb_0(0,0).
\end{equation} 
Note that $b_0(0,0)=b^0_0$. From \eqref{s4-e1}, \eqref{s4-e11} and \eqref{s4-e13}, we get 
\begin{equation} \label{s4-e15} 
-cb_0(0,0)=\abs{\lambda_1}^{-1}\cdots\abs{\lambda_n}^{-1}\pi^n\hat C_0.
\end{equation} 
Now, we compute $\hat C_0$. The computation is very complicate. We only give the outline of the computation. We write $C_0=I+II+III+IV+V$, where 
\begin{equation} \label{s4-e15-1}
\begin{split}
&I=\frac{1}{2}\sum^n_{j=1}\frac{1}{\abs{\lambda_j}}\Bigr(\frac{\pr^2b_0}{\pr z_j\pr\ol z_j}(z,0)|_{z=0}b_0(0,0)+
\frac{\pr^2 b_0}{\pr z_j\pr\ol z_j}(0,z)|_{z=0}b_0(0,0)\\
&+(\frac{\pr b_0}{\pr\ol z_j})(0,z)|_{z=0}\frac{\pr b_0}{\pr z_j}(z,0)|_{z=0}+\frac{\pr b_0}{\pr z_j}(0,z)|_{z=0}
\frac{\pr b_0}{\pr\ol z_j}(z,0)|_{z=0}\Bigr),
\end{split}
\end{equation} 
\begin{equation} \label{s4-e15-2}
II=-\frac{1}{4}\sum^n_{j,t=1}\frac{1}{\abs{\lambda_j}
\abs{\lambda_t}}\frac{\pr^4{\rm Im\,}\psi(z,0)}{\pr\ol z_j\pr z_j\pr\ol z_t\pr z_t}|_{z=0}b_0(0,0)^2,
\end{equation} 
\begin{equation} \label{s4-e15-3}
\begin{split}
&III=-\frac{1}{2}\sum^n_{j,s=1}\frac{1}{\abs{\lambda_j}\abs{\lambda_s}}\frac{\pr^3{\rm Im\,}\psi(z,0)}{\pr\ol z_j\pr z_j\pr\ol z_s}|_{z=0}\Bigr(\frac{\pr b_0}{\pr z_s}(0,z)|_{z=0}+\frac{\pr b_0}{\pr z_s}(z,0)|_{z=0}\Bigr)b_0(0,0)\\
&-\frac{1}{2}\sum^n_{j,t=1}\frac{1}{\abs{\lambda_j}\abs{\lambda_s}}\frac{\pr^3{\rm Im\,}\psi(z,0)}{\pr\ol z_j\pr z_j\pr z_s}|_{z=0}\Bigr(\frac{\pr b_0}{\pr\ol z_s}(0,z)|_{z=0}+\frac{\pr b_0}{\pr\ol z_s}(z,0)|_{z=0}\Bigr)b_0(0,0),
\end{split}
\end{equation} 
\begin{equation} \label{s4-e15-4}
IV=\frac{1}{4}\sum^n_{j,t,s=1}\frac{1}{\abs{\lambda_j}\abs{\lambda_t}\abs{\lambda_s}}\abs{\frac{\pr^3{\rm Im\,}\psi(z,0)}{\pr z_j\pr z_t\pr\ol z_s}|_{z=0}}^2b_0(0,0)^2
\end{equation} 
and 
\begin{equation} \label{s4-e15-5}
\begin{split}
&V=\frac{1}{2}\sum^n_{j,t,s=1}\frac{1}{\abs{\lambda_j}\abs{\lambda_t}\abs{\lambda_s}}\frac{\pr^3{\rm Im\,}\psi(z,0)}{\pr\ol z_j\pr z_j\pr\ol z_s}|_{z=0}\frac{\pr^3{\rm Im\,}\psi(z,0)}{\pr\ol z_t\pr z_t\pr z_s}|_{z=0}b_0(0,0)^2\\
&+\frac{1}{12}\sum^n_{j,t,s=1}\frac{1}{\abs{\lambda_j}\abs{\lambda_t}\abs{\lambda_s}}\abs{\frac{\pr^3{\rm Im\,}\psi(z,0)}{\pr z_j\pr z_t\pr z_s}|_{z=0}}^2b_0(0,0)^2.
\end{split}
\end{equation} 

Now, we deal with the sum $I$. From \eqref{s3-e17} and \eqref{s3-e54-1}, we may write 
\begin{equation} \label{s4-e16} 
\begin{split}
&b_0(z,0)=\sum_{1\leq s\leq n}(a_sz_s+b_s\ol z_s)b_0(0,0)\\
&+\sum_{q+1\leq j\leq n, 1\leq t\leq q,1\leq s\leq n}
(\alpha^s_{j,t}z_s+\beta^s_{j,t}\ol z_s)d\ol z_t^{\wedge,*}d\ol z_j^\wedge b_0(0,0)\\
&+\sum_{q+1\leq j\leq n,1\leq t\leq q,1\leq s\leq n}
(\gamma^ s_{j,t}z_s+\delta^s_{j,t}\ol z_s)b_0(0,0)d\ol z_j^{\wedge,*}d\ol z_t^\wedge \\
&+\sum_{1\leq s\leq n}d_s\abs{z_s}^2b_0(0,0)+\sum_{q+1\leq j\leq n,1\leq t\leq q,1\leq s\leq n}\sigma^s_{j,t}\abs{z_s}^2d\ol z_t^{\wedge,*}d\ol z_j^\wedge b_0(0,0)d\ol z_j^{\wedge,*}d\ol z_t^\wedge\\
&+r(z)+h(z),
\end{split}
\end{equation} 
where 
\[a_s, b_s, \alpha^s_{j,t}, \beta^s_{j,t}, \gamma^s_{j,t}, \delta^s_{j,t}, d_s, \sigma^s_{j,t}\in\Complex,\ 
s=1,\ldots,n, j=q+1,\ldots,n, t=1,\ldots,q,\] 
$r(z)=O(\abs{z}^2)$, $h(z)=O(\abs{z}^2)$, ${\rm Tr\,}h=0$ and $\frac{\pr^2r}{\pr\ol z_j\pr z_j}(0)=0$, $j=1,\ldots,n$. As in section $3$, for $u(z)\in\mathscr L(\Lambda^{0,q}T^*(\Complex^n),\Lambda^{0,q}T^*(\Complex^n))$, 
we write $u^*(z)$ to denote the adjoint of $u(z)$ with respect to $(\ |\ )$ in $\mathscr L(\Lambda^{0,q}T^*(\Complex^n),\Lambda^{0,q}T^*(\Complex^n))$. We notice that 
\[b_0(0,z,k)=b_0(z,0,k)^*.\] 
From this and \eqref{s4-e16}, we deduce that 
\begin{equation} \label{s4-e17} 
\begin{split}
&b_0(0,z)=\sum_{1\leq s\leq n}(\ol a_s\ol z_s+\ol b_sz_s)b_0(0,0)\\
&+\sum_{q+1\leq j\leq n, 1\leq t\leq q,1\leq s\leq n}
(\ol{\alpha^s_{j,t}}\ol z_s+\ol{\beta^s_{j,t}}z_s)b_0(0,0) d\ol z_j^{\wedge,*}d\ol z_t^\wedge\\
&+\sum_{q+1\leq j\leq n,1\leq t\leq q,1\leq s\leq n}
(\ol{\gamma^ s_{j,t}}\ol z_s+\ol{\delta^s_{j,t}}z_s)d\ol z_t^{\wedge,*}d\ol z_j^\wedge b_0(0,0)\\
&+\sum_{1\leq s\leq n}\ol{d_s}\abs{z_s}^2b_0(0,0)+\sum_{q+1\leq j\leq n,1\leq t\leq q,1\leq s\leq n}\ol{\sigma^s_{j,t}}\abs{z_s}^2d\ol z_t^{\wedge,*}d\ol z_j^\wedge b_0(0,0)d\ol z_j^{\wedge,*}d\ol z_t^\wedge\\
&+r^*(z)+h^*(z).
\end{split}
\end{equation} 
Note that we still have ${\rm Tr\,}h^*=0$ and $\frac{\pr^2r^*}{\pr\ol z_j\pr z_j}(0)=0$, $j=1,\ldots,n$.
From \eqref{s4-e16} and \eqref{s4-e17}, we can check that

\begin{equation} \label{s4-e18}
\begin{split}
&\frac{1}{2}\sum^n_{s=1}\frac{1}{\abs{\lambda_s}}\Bigr(\frac{\pr^2b_0}{\pr z_s\pr\ol z_s}(z,0)|_{z=0}b_0(0,0)+
\frac{\pr^2 b_0}{\pr z_s\pr\ol z_s}(0,z)|_{z=0}b_0(0,0) \\ 
&+(\frac{\pr b_0}{\pr\ol z_s})(0,z)|_{z=0}\frac{\pr b_0}{\pr z_s}(z,0)|_{z=0}+\frac{\pr b_0}{\pr z_s}(0,z)|_{z=0}
\frac{\pr b_0}{\pr\ol z_s}(z,0)|_{z=0}\Bigr)\\
&=\frac{1}{2}\sum^n_{s=1}\frac{1}{\abs{\lambda_s}}(d_s+\ol{d_s})b_0(0,0)^2
+\frac{1}{2}\sum^n_{s=1}\frac{1}{\abs{\lambda_s}}(\abs{a_s}^2+\abs{b_s}^2)\\
&+\frac{1}{2}\sum_{1\leq s\leq n,q+1\leq j\leq n,1\leq t\leq q}\frac{1}{\abs{\lambda_s}}(\abs{\alpha^s_{j,t}}^2+\abs{\beta^s_{j,t}}^2)b_0(0,0)^2\\
&+R_1,
\end{split}
\end{equation} 
where $R_1\bot d\ol z_1^\wedge d\ol z_1^{\wedge,*}\cdots d\ol z_q^\wedge d\ol z_q^{\wedge,*}$.
From \eqref{s3-e17}, \eqref{s3-e54-1} and the fourth order Taylor expansion of the phase $\psi$ (see \eqref{s2-e3} \eqref{s2-e4} 
and \eqref{s2-e5}), we can write down the explicit formula of the right side of \eqref{s4-e18}. The computation is very complicate but elemenrary and is therefore omitted. We state our result 

\begin{prop} \label{s4-p2} 
For $I$ in \eqref{s4-e15-1}, we have 
\begin{equation} \label{s4-e21} 
\begin{split} 
&I=\sum_{q+1\leq j\leq n,1\leq t,s\leq q}\frac{3\abs{\lambda_j}^2\abs{\lambda_t}-\abs{\lambda_t}^2\abs{\lambda_j}-2\abs{\lambda_t}^2
\abs{\lambda_s}+2\abs{\lambda_j}\abs{\lambda_t}\abs{\lambda_s}}{2(\abs{\lambda_j}+\abs{\lambda_s})^2\abs{\lambda_j}\abs{\lambda_t}^2(\abs{\lambda_j}+\abs{\lambda_t})}\\
&\times\abs{\frac{\pr^3\phi}{\pr\ol z_s\pr z_j\pr z_t}(0)}^2b_0(0,0)^2\\
&+\sum_{q+1\leq j,t\leq n,1\leq s\leq q}\alpha_{s,j,t}\abs{\frac{\pr^3\phi}{\pr\ol z_s\pr z_j\pr z_t}(0)}^2b_0(0,0)^2\\
&+\sum_{q+1\leq j,s\leq n,1\leq t\leq q}\frac{3\abs{\lambda_j}\abs{\lambda_t}^2-\abs{\lambda_t}\abs{\lambda_j}^2-2\abs{\lambda_j}^2
\abs{\lambda_s}+2\abs{\lambda_j}\abs{\lambda_t}\abs{\lambda_s}}{2(\abs{\lambda_t}+\abs{\lambda_s})^2\abs{\lambda_j}^2\abs{\lambda_t}(\abs{\lambda_j}+\abs{\lambda_t})}\\
&\times\abs{\frac{\pr^3\phi}{\pr\ol z_s\pr z_j\pr z_t}(0)}^2b_0(0,0)^2\\ 
&+\sum_{1\leq j,t\leq q,q+1\leq s\leq n}\alpha_{s,j,t}\abs{\frac{\pr^3\phi}{\pr\ol z_s\pr z_j\pr z_t}(0)}^2b^2_0(0,0)\\
&+\frac{1}{2}\sum_{q+1\leq j,t\leq n,1\leq s\leq q}\frac{-\abs{\lambda_s}^2+\abs{\lambda_j}\abs{\lambda_t}+\abs{\lambda_s}(\abs{\lambda_j}+\abs{\lambda_t})}{(\abs{\lambda_j}+\abs{\lambda_s})(\abs{\lambda_t}+\abs{\lambda_s})\abs{\lambda_j}\abs{\lambda_t}\abs{\lambda_s}}\\
&\times\frac{\pr^3\phi}{\pr\ol z_j\pr z_j\pr z_s}(0)\frac{\pr^3\phi}{\pr\ol z_t\pr z_t\pr\ol z_s}(0)b_0(0,0)^2\\
&+\frac{1}{2}\sum_{1\leq j,t\leq q,q+1\leq s\leq n}\frac{-\abs{\lambda_s}^2+\abs{\lambda_j}\abs{\lambda_t}+\abs{\lambda_s}(\abs{\lambda_j}+\abs{\lambda_t})}{(\abs{\lambda_j}+\abs{\lambda_s})(\abs{\lambda_t}+\abs{\lambda_s})\abs{\lambda_j}\abs{\lambda_t}\abs{\lambda_s}}\\
&\times\frac{\pr^3\phi}{\pr\ol z_j\pr z_j\pr z_s}(0)\frac{\pr^3\phi}{\pr\ol z_t\pr z_t\pr\ol z_s}(0)b^2_0(0,0)\\
&+\frac{1}{2}\sum_{q+1\leq j\leq n,1\leq t,s\leq q}\frac{\abs{\lambda_j}\abs{\lambda_t}\abs{\lambda_s}-\abs{\lambda_j}^2\abs{\lambda_t}-\abs{\lambda_j}\abs{\lambda_t}^2-
\abs{\lambda_s}\abs{\lambda_t}^2}{(\abs{\lambda_j}+\abs{\lambda_s})(\abs{\lambda_j}+\abs{\lambda_t})\abs{\lambda_j}\abs{\lambda_t}^2\abs{\lambda_s}}\\
&\times\Bigr(\frac{\pr^3\phi}{\pr\ol z_j\pr z_j\pr z_s}(0)\frac{\pr^3\phi}{\pr\ol z_t\pr z_t\pr\ol z_s}(0)+
\frac{\pr^3\phi}{\pr\ol z_j\pr z_j\pr\ol z_s}(0)\frac{\pr^3\phi}{\pr\ol z_t\pr z_t\pr z_s}(0)\Bigr)b_0(0,0)^2\\
&+\frac{1}{2}\sum_{q+1\leq j,s\leq n,1\leq t\leq q}\frac{\abs{\lambda_j}\abs{\lambda_t}\abs{\lambda_s}-\abs{\lambda_j}^2\abs{\lambda_t}-\abs{\lambda_j}\abs{\lambda_t}^2-
\abs{\lambda_s}\abs{\lambda_j}^2}{(\abs{\lambda_t}+\abs{\lambda_s})(\abs{\lambda_j}+\abs{\lambda_t})\abs{\lambda_j}^2\abs{\lambda_t}\abs{\lambda_s}}\\
&\times\Bigr(\frac{\pr^3\phi}{\pr\ol z_j\pr z_j\pr z_s}(0)\frac{\pr^3\phi}{\pr\ol z_t\pr z_t\pr\ol z_s}(0)+
\frac{\pr^3\phi}{\pr\ol z_j\pr z_j\pr\ol z_s}(0)\frac{\pr^3\phi}{\pr\ol z_t\pr z_t\pr z_s}(0)\Bigr)b_0(0,0)^2\\
&+\frac{1}{2}\sum_{q+1\leq j,t,s\leq n}\frac{1}{\abs{\lambda_j}\abs{\lambda_t}\abs{\lambda_s}}\frac{\pr^3\phi}{\pr\ol z_j\pr z_j\pr z_s}(0)\frac{\pr^3\phi}{\pr\ol z_t\pr z_t\pr\ol z_s}(0)b_0(0,0)^2\\
&+\frac{1}{2}\sum_{1\leq j,t,s\leq q}\frac{1}{\abs{\lambda_j}\abs{\lambda_t}\abs{\lambda_s}}\frac{\pr^3\phi}{\pr\ol z_j\pr z_j\pr z_s}(0)\frac{\pr^3\phi}{\pr\ol z_t\pr z_t\pr\ol z_s}(0)b_0(0,0)^2\\
&+\sum_{q+1\leq j\leq n,1\leq t\leq q}\frac{1}{\abs{\lambda_j}\abs{\lambda_t}}\frac{\abs{\lambda_j}-\abs{\lambda_t}}{\abs{\lambda_j}+\abs{\lambda_t}}\frac{\pr^4\phi}{\pr\ol z_j\pr z_j\pr\ol z_t\pr z_t}(0)b_0(0,0)^2\\
&+R_1,
\end{split}
\end{equation}
where $R_1\bot d\ol z_1^\wedge d\ol z_1^{\wedge,*}\cdots d\ol z_q^\wedge d\ol z_q^{\wedge,*}$ and 
\begin{equation} \label{s4-e22} 
\begin{split}
&\alpha_{s,j,t}=
\frac{1}{4}\frac{1}{(\abs{\lambda_j}+\abs{\lambda_t}+\abs{\lambda_s})^2}\bigr(\frac{1}{\abs{\lambda_j}}+\frac{1}{\abs{\lambda_t}}\bigr)\\
&+\frac{\abs{\lambda_t}^2\abs{\lambda_j}\abs{\lambda_s}-\abs{\lambda_j}^2\abs{\lambda_t}^2}{2\abs{\lambda_j}\abs{\lambda_t}\abs{\lambda_s}(\abs{\lambda_j}+\abs{\lambda_s})(\abs{\lambda_j}+\abs{\lambda_t}+\abs{\lambda_s})^2}\bigr(\frac{1}{\abs{\lambda_j}+\abs{\lambda_s}}+\frac{1}{\abs{\lambda_t}+\abs{\lambda_s}}\bigr)\\
&+\frac{\abs{\lambda_t}\abs{\lambda_j}^2\abs{\lambda_s}-\abs{\lambda_j}^2\abs{\lambda_t}^2}{2\abs{\lambda_j}\abs{\lambda_t}\abs{\lambda_s}(\abs{\lambda_t}+\abs{\lambda_s})(\abs{\lambda_j}+\abs{\lambda_t}+\abs{\lambda_s})^2}\bigr(\frac{1}{\abs{\lambda_j}+\abs{\lambda_s}}+\frac{1}{\abs{\lambda_t}+\abs{\lambda_s}}\bigr).
\end{split}
\end{equation} 
\end{prop}
 
Now, we deal with the sum $II$. From \eqref{s2-e3}, \eqref{s2-e4} and \eqref{s2-e5}, we can check that 
\begin{equation} \label{s4-e23}
\begin{split}
&-\frac{1}{4}\sum_{1\leq j,s\leq n}\frac{1}{\abs{\lambda_j}\abs{\lambda_s}}\frac{\pr^4{\rm Im\,}\psi(z,0)}{\pr\ol z_j\pr z_j\pr\ol z_s\pr z_s}(0)b_0(0,0)^2\\
&=\sum_{q+1\leq j,t\leq n,1\leq s\leq q}\beta_{s,j,t}\abs{\frac{\pr^3\phi}{\pr\ol z_s\pr z_j\pr z_t}(0)}^2\\
&+\sum_{1\leq j,t\leq q,q+1\leq s\leq n}\beta_{s,j,t}\abs{\frac{\pr^3\phi}{\pr\ol z_s\pr z_j\pr z_t}(0)}^2\\
&-\sum_{1\leq s,t\leq q,q+1\leq j\leq n}\frac{1}{\abs{\lambda_t}(\abs{\lambda_j}+\abs{\lambda_t})(\abs{\lambda_s}+\abs{\lambda_j})}
\abs{\frac{\pr^3\phi}{\pr\ol z_s\pr z_j\pr z_t}(0)}^2\\
&-\sum_{q+1\leq s,j\leq n,1\leq t\leq q}\frac{1}{\abs{\lambda_j}(\abs{\lambda_j}+\abs{\lambda_t})(\abs{\lambda_s}+\abs{\lambda_t})}
\abs{\frac{\pr^3\phi}{\pr\ol z_s\pr z_j\pr z_t}(0)}^2\\
&-\frac{1}{2}\sum_{1\leq s,t\leq q,q+1\leq j\leq n}\frac{1}{\abs{\lambda_t}(\abs{\lambda_j}+\abs{\lambda_s})(\abs{\lambda_j}+\abs{\lambda_t})}\\
&\times\Bigr(\frac{\pr^3\phi}{\pr z_s\pr\ol z_j\pr z_j}(0)
\frac{\pr^3\phi}{\pr\ol z_s\pr\ol z_t\pr z_t}(0)+\frac{\pr^3\phi}{\pr\ol z_s\pr\ol z_j\pr z_j}(0)\frac{\pr^3\phi}{\pr z_s\pr\ol z_t\pr z_t}(0)\Bigr)\\
&-\frac{1}{2}\sum_{q+1\leq s,j\leq n,1\leq t\leq q}\frac{1}{\abs{\lambda_j}(\abs{\lambda_t}+\abs{\lambda_s})(\abs{\lambda_j}+\abs{\lambda_t})}\\
&\times\Bigr(\frac{\pr^3\phi}{\pr z_s\pr\ol z_j\pr z_j}(0)
\frac{\pr^3\phi}{\pr\ol z_s\pr\ol z_t\pr z_t}(0)+\frac{\pr^3\phi}{\pr\ol z_s\pr\ol z_j\pr z_j}(0)\frac{\pr^3\phi}{\pr z_s\pr\ol z_t\pr z_t}(0)\Bigr)\\
&-\frac{1}{2}\sum_{q+1\leq j\leq n,1\leq t\leq q}\frac{\abs{\lambda_j}-\abs{\lambda_t}}{\abs{\lambda_j}\abs{\lambda_t}(\abs{\lambda_j}+\abs{\lambda_t})}\frac{\pr^4\phi}{\pr\ol z_j\pr z_j\pr\ol z_t\pr z_t}(0)\\
&+\frac{1}{4}\sum_{1\leq j,t\leq q}\frac{1}{\abs{\lambda_j}\abs{\lambda_t}}\frac{\pr^4\phi}{\pr\ol z_j\pr z_j\pr\ol z_t\pr z_t}(0)-\frac{1}{4}\sum_{q+1\leq j,t\leq n}\frac{1}{\abs{\lambda_j}\abs{\lambda_t}}\frac{\pr^4\phi}{\pr\ol z_j\pr z_j\pr\ol z_t\pr z_t}(0),
\end{split}
\end{equation}
where 
\begin{equation} \label{s4-e24}
\begin{split}
&\beta_{s,j,t}=
-\frac{1}{2(\abs{\lambda_s}+\abs{\lambda_j}+\abs{\lambda_t})^2\abs{\lambda_j}\abs{\lambda_t}}\\
&\times\Bigr(\abs{\lambda_s}+
\abs{\lambda_j}+\abs{\lambda_t}-\frac{\abs{\lambda_j}\abs{\lambda_t}}{\abs{\lambda_j}+\abs{\lambda_s}}-\frac{\abs{\lambda_j}\abs{\lambda_t}}{\abs{\lambda_t}+\abs{\lambda_s}}\Bigr)\\
&-\frac{1}{2}\frac{\abs{\lambda_t}}{(\abs{\lambda_j}+\abs{\lambda_s})(\abs{\lambda_j}+\abs{\lambda_s}+\abs{\lambda_t})^2}
\bigr(\frac{1}{\abs{\lambda_j}+\abs{\lambda_s}}+\frac{1}{\abs{\lambda_t}+\abs{\lambda_s}}\bigr)\\
&-\frac{1}{2}\frac{\abs{\lambda_j}}{(\abs{\lambda_t}+\abs{\lambda_s})(\abs{\lambda_j}+\abs{\lambda_s}+\abs{\lambda_t})^2}
\bigr(\frac{1}{\abs{\lambda_j}+\abs{\lambda_s}}+\frac{1}{\abs{\lambda_t}+\abs{\lambda_s}}\bigr).
\end{split}
\end{equation} 

Now, we deal wih the sum $III$, $IV$ and $V$.  First we need the following theorem which follows from Theorem~\ref{s2-t1} 

\begin{thm} \label{s4-t3} 
We have that 
\begin{align} \label{s4-e25} 
&\frac{\pr^3{\rm Im\,}\psi(z,0)}{\pr\ol z_s\pr z_j\pr z_t}(0)=\frac{\pr^3\phi}{\pr\ol z_s\pr z_j\pr z_t}(0),\ \ q+1\leq j,t,s\leq n,\nonumber \\
&\frac{\pr^3{\rm Im\,}\psi(z,0)}{\pr\ol z_s\pr z_j\pr z_t}(0)=\frac{\abs{\lambda_s}}{\abs{\lambda_s}+\abs{\lambda_t}}\frac{\pr^3\phi}{\pr\ol z_s\pr z_j\pr z_t}(0),\ \ q+1\leq j,s\leq n,1\leq t\leq q\nonumber\\
&\frac{\pr^3{\rm Im\,}\psi(z,0)}{\pr\ol z_s\pr z_j\pr z_t}(0)=\frac{1}{\abs{\lambda_j}+\abs{\lambda_t}+\abs{\lambda_s}}\Bigr(\abs{\lambda_s}-\frac{\abs{\lambda_j}\abs{\lambda_t}}{\abs{\lambda_j}+\abs{\lambda_s}}-\frac{\abs{\lambda_j}\abs{\lambda_t}}{\abs{\lambda_s}+\abs{\lambda_t}}\Bigr)\nonumber\\
&\times\frac{\pr^3\phi}{\pr\ol z_s\pr z_j\pr z_t}(0),\ \ q+1\leq s\leq n,1\leq j,t\leq q,\nonumber\\
&\frac{\pr^3{\rm Im\,}\psi(z,0)}{\pr\ol z_s\pr z_j\pr z_t}(0)=-\frac{\pr^3\phi}{\pr\ol z_s\pr z_j\pr z_t}(0),\ \ 1\leq j,t,s\leq q,\nonumber\\
&\frac{\pr^3{\rm Im\,}\psi(z,0)}{\pr\ol z_s\pr z_j\pr z_t}(0)=-\frac{\abs{\lambda_s}}{\abs{\lambda_j}+\abs{\lambda_s}}\frac{\pr^3\phi}{\pr\ol z_s\pr z_j\pr z_t}(0),\ \ 1\leq s,t\leq q,q+1\leq j\leq n,\nonumber\\
&\frac{\pr^3{\rm Im\,}\psi(z,0)}{\pr\ol z_s\pr z_j\pr z_t}(0)=\frac{1}{\abs{\lambda_j}+\abs{\lambda_t}+\abs{\lambda_s}}\Bigr(-\abs{\lambda_s}+\frac{\abs{\lambda_t}\abs{\lambda_j}}{\abs{\lambda_j}+\abs{\lambda_s}}+\frac{\abs{\lambda_t}\abs{\lambda_j}}{\abs{\lambda_s}+\abs{\lambda_t}}\Bigr)\nonumber\\
&\times\frac{\pr^3\phi}{\pr\ol z_s\pr z_j\pr z_t}(0),\ \ q+1\leq j,t\leq n,1\leq s\leq q,\nonumber\\
&\frac{\pr^3{\rm Im\,}\psi(z,0)}{\pr z_s\pr z_j\pr z_t}(0)=0,\ \ 1\leq s,j,t\leq n.
\end{align}
\end{thm} 

Note that 
\begin{equation*} 
\begin{split}
&\frac{\pr b_0(0,z)}{\pr z_s}|_{z=0}+\frac{\pr b_0(z,0)}{\pr z_s}|_{z=0}\\
&=\Bigr(-\sum^q_{j=1}\frac{1}{\abs{\lambda_j}}\frac{\pr^3\phi}{\pr\ol z_j\pr z_j\pr z_s}(0)+\sum^n_{j=q+1}\frac{1}{\abs{\lambda_j}}\frac{\pr^3\phi}{\pr\ol z_j\pr z_j\pr z_s}(0)\Bigr)b_0(0,0),
\end{split}
\end{equation*}
and 
\begin{equation*} 
\begin{split}
&\frac{\pr b_0(0,z)}{\pr\ol z_s}|_{z=0}+\frac{\pr b_0(z,0)}{\pr\ol z_s}|_{z=0}\\
&=\Bigr(-\sum^q_{j=1}\frac{1}{\abs{\lambda_j}}\frac{\pr^3\phi}{\pr\ol z_j\pr z_j\pr\ol z_s}(0)+\sum^n_{j=q+1}\frac{1}{\abs{\lambda_j}}\frac{\pr^3\phi}{\pr\ol z_j\pr\ol z_j\pr z_s}(0)\Bigr)b_0(0,0),
\end{split} 
\end{equation*} 
$s=1,\ldots,n$. From this and \eqref{s4-e25}, we can check that 
\begin{equation} \label{s4-e26} 
\begin{split}
&III=-\frac{1}{2}\sum^n_{j,s=1}\frac{1}{\abs{\lambda_j}\abs{\lambda_s}}\frac{\pr^3{\rm Im\,}\psi(z,0)}{\pr\ol z_j\pr z_j\pr z_s}|_{z=0}\Bigr(\frac{\pr b_0}{\pr\ol z_s}(0,z)|_{z=0}+\frac{\pr b_0}{\pr\ol z_s}(z,0)|_{z=0}\Bigr)b_0(0,0)\\
&=-\sum_{1\leq j,t, s\leq q}\frac{1}{\abs{\lambda_j}\abs{\lambda_s}\abs{\lambda_t}}
\frac{\pr^3\phi}{\pr\ol z_j\pr z_j\pr z_s}(0)\frac{\pr^3\phi}{\pr\ol z_t\pr z_t\pr\ol z_s}(0)b_0(0,0)^2 \\
&+\frac{1}{2}\sum_{q+1\leq j\leq n,1\leq t,s\leq q}\frac{2\abs{\lambda_j}+\abs{\lambda_s}}{\abs{\lambda_j}\abs{\lambda_t}\abs{\lambda_s}(\abs{\lambda_j}+\abs{\lambda_s})}\\
&\times\Bigr(\frac{\pr^3\phi}{\pr\ol z_j\pr z_j\pr z_s}(0)\frac{\pr^3\phi}{\pr\ol z_t\pr z_t\pr\ol z_s}(0)+\frac{\pr^3\phi}{\pr\ol z_j\pr z_j\pr\ol z_s}(0)\frac{\pr^3\phi}{\pr\ol z_t\pr z_t\pr z_s}(0)\Bigr)b_0(0,0)^2\\
&-\frac{1}{2}\sum_{q+1\leq j,t\leq n,1\leq s\leq q}\frac{\abs{\lambda_j}}{\abs{\lambda_j}\abs{\lambda_t}\abs{\lambda_s}(\abs{\lambda_j}+\abs{\lambda_s})}\\
&\times\Bigr(\frac{\pr^3\phi}{\pr\ol z_j\pr z_j\pr z_s}(0)
\frac{\pr^3\phi}{\pr\ol z_t\pr z_t\pr\ol z_s}(0)+\frac{\pr^3\phi}{\pr\ol z_j\pr z_j\pr\ol z_s}(0)\frac{\pr^3\phi}{\pr\ol z_t\pr z_t\pr z_s}(0)\Bigr)b_0(0,0)^2\\
&-\sum_{q+1\leq j,t, s\leq n}\frac{1}{\abs{\lambda_j}\abs{\lambda_s}\abs{\lambda_t}}
\frac{\pr^3\phi}{\pr\ol z_j\pr z_j\pr z_s}(0)\frac{\pr^3\phi}{\pr\ol z_t\pr z_t\pr\ol z_s}(0)b_0(0,0)^2 \\
&+\frac{1}{2}\sum_{q+1\leq j,s\leq n,1\leq t\leq q}\frac{2\abs{\lambda_t}+\abs{\lambda_s}}{\abs{\lambda_j}\abs{\lambda_t}\abs{\lambda_s}(\abs{\lambda_t}+\abs{\lambda_s})}\\
&\times\Bigr(\frac{\pr^3\phi}{\pr\ol z_j\pr z_j\pr z_s}(0)\frac{\pr^3\phi}{\pr\ol z_t\pr z_t\pr\ol z_s}(0)+\frac{\pr^3\phi}{\pr\ol z_j\pr z_j\pr\ol z_s}(0)\frac{\pr^3\phi}{\pr\ol z_t\pr z_t\pr z_s}(0)\Bigr)b_0(0,0)^2\\
&-\frac{1}{2}\sum_{q+1\leq s\leq n,1\leq j,t\leq q}\frac{\abs{\lambda_t}}{\abs{\lambda_j}\abs{\lambda_t}\abs{\lambda_s}(\abs{\lambda_t}+\abs{\lambda_s})}\\
&\times\Bigr(\frac{\pr^3\phi}{\pr\ol z_j\pr z_j\pr z_s}(0)
\frac{\pr^3\phi}{\pr\ol z_t\pr z_t\pr\ol z_s}(0)+\frac{\pr^3\phi}{\pr\ol z_j\pr z_j\pr\ol z_s}(0)\frac{\pr^3\phi}{\pr\ol z_t\pr z_t\pr z_s}(0)\Bigr)b_0(0,0)^2,
\end{split}
\end{equation} 

\begin{equation} \label{s4-e27} 
\begin{split} 
&IV=\frac{1}{4}\sum_{q+1\leq j,t,s\leq n}\frac{1}{\abs{\lambda_j}\abs{\lambda_t}\abs{\lambda_s}}\abs{\frac{\pr^3\phi}
{\pr\ol z_s\pr z_j\pr z_t}(0)}^2b_0(0,0)^2\\
&+\frac{1}{2}\sum_{q+1\leq j,s\leq n,1\leq t\leq q}\frac{\abs{\lambda_s}}{\abs{\lambda_j}\abs{\lambda_t}(\abs{\lambda_s}+\abs{\lambda_t})^2}\abs{\frac{\pr^3\phi}{\pr\ol z_s\pr z_j\pr z_t}(0)}^2b^2_0(0,0)\\
&+\frac{1}{4}\sum_{q+1\leq s\leq n,1\leq j,t\leq q}\gamma_{s,j,t}\abs{\frac{\pr^3\phi}{\pr\ol z_s\pr z_j\pr z_t}(0)}^2b_0(0,0)^2\\
&+\frac{1}{4}\sum_{1\leq j,t,s\leq q}\frac{1}{\abs{\lambda_j}\abs{\lambda_t}\abs{\lambda_s}}\abs{\frac{\pr^3\phi}
{\pr\ol z_s\pr z_j\pr z_t}(0)}^2b_0(0,0)^2\\
&+\frac{1}{2}\sum_{1\leq s,t\leq q,q+1\leq j\leq n}\frac{\abs{\lambda_s}}{\abs{\lambda_j}\abs{\lambda_t}(\abs{\lambda_s}+\abs{\lambda_j})^2}\abs{\frac{\pr^3\phi}{\pr\ol z_s\pr z_j\pr z_t}(0)}^2b_0(0,0)^2\\
&+\frac{1}{4}\sum_{1\leq s\leq q,q+1\leq j,t\leq n}\gamma_{s,j,t}\abs{\frac{\pr^3\phi}{\pr\ol z_s\pr z_j\pr z_t}(0)}^2b_0(0,0)^2,
\end{split}
\end{equation} 
where
\begin{equation} \label{s4-e27-1} 
\begin{split}
&\gamma_{s,j,t}=\frac{1} {\abs{\lambda_j}\abs{\lambda_t}\abs{\lambda_s}(\abs{\lambda_j}+\abs{\lambda_s}+\abs{\lambda_t})^2}\\
&\times\Bigr(\abs{\lambda_s}^2+\frac{\abs{\lambda_j}^2\abs{\lambda_t}^2}{(\abs{\lambda_s}+\abs{\lambda_j})^2}+\frac{\abs{\lambda_j}^2\abs{\lambda_t}^2}{(\abs{\lambda_s}+\abs{\lambda_t})^2}-\frac{2\abs{\lambda_j}\abs{\lambda_t}\abs{\lambda_s}}{\abs{\lambda_s}+\abs{\lambda_j}}\\
&-\frac{2\abs{\lambda_j}\abs{\lambda_t}\abs{\lambda_s}}{\abs{\lambda_s}+\abs{\lambda_t}}+
\frac{2\abs{\lambda_j}^2\abs{\lambda_t}^2}{(\abs{\lambda_s}+\abs{\lambda_j})(\abs{\lambda_s}+\abs{\lambda_t})}\Bigr). 
\end{split} 
\end{equation}
and 
\begin{equation} \label{s4-e28} 
\begin{split}
&V=\frac{1}{2}\sum^n_{j,k,s=1}\frac{1}{\abs{\lambda_j}\abs{\lambda_k}\abs{\lambda_s}}\frac{\pr^3{\rm Im\,}\psi(z,0)}{\pr\ol z_j\pr z_j\pr\ol z_s}|_{z=0}\frac{\pr^3{\rm Im\,}\psi(z,0)}{\pr\ol z_k\pr z_k\pr z_s}|_{z=0}b_0(0,0)^2\\
&+\frac{1}{12}\sum^n_{j,k,s=1}\frac{1}{\abs{\lambda_j}\abs{\lambda_k}\abs{\lambda_s}}\abs{\frac{\pr^3{\rm Im\,}\psi(z,0)}{\pr z_j\pr z_k\pr z_s}|_{z=0}}^2b_0(0,0)^2\\
&=\frac{1}{2}\sum_{q+1\leq j,t,s\leq n}\frac{1}{\abs{\lambda_j}\abs{\lambda_t}\abs{\lambda_s}}\frac{\pr^3\phi}{\pr\ol z_j\pr z_j\pr z_s}(0)\frac{\pr^3\phi}{\pr\ol z_t\pr z_t\pr\ol z_s}(0)b_0(0,0)^2\\
&-\frac{1}{2}\sum_{q+1\leq s,j\leq n,1\leq t\leq q}\frac{1}{(\abs{\lambda_t}+\abs{\lambda_s})\abs{\lambda_j}\abs{\lambda_s}}\\
&\times\Bigr(\frac{\pr^3\phi}{\pr\ol z_j\pr z_j\pr z_s}(0)\frac{\pr^3\phi}{\pr\ol z_t\pr z_t\pr\ol z_s}(0)+\frac{\pr^3\phi}{\pr\ol z_j\pr z_j\pr\ol z_s}(0)\frac{\pr^3\phi}{\pr\ol z_t\pr z_t\pr z_s}(0)\Bigr)b_0(0,0)^2\\
&+\frac{1}{2}\sum_{q+1\leq s\leq n,1\leq j,t\leq q}\frac{1}{(\abs{\lambda_t}+\abs{\lambda_s})(\abs{\lambda_j}+\abs{\lambda_s})\abs{\lambda_s}}\\
&\times\frac{\pr^3\phi}{\pr\ol z_j\pr z_j\pr z_s}(0)\frac{\pr^3\phi}{\pr\ol z_t\pr z_t\pr\ol z_s}(0)b_0(0,0)^2\\
&+\frac{1}{2}\sum_{1\leq j,t,s\leq q}\frac{1}{\abs{\lambda_j}\abs{\lambda_t}\abs{\lambda_s}}\frac{\pr^3\phi}{\pr\ol z_j\pr z_j\pr z_s}(0)\frac{\pr^3\phi}{\pr\ol z_t\pr z_t\pr\ol z_s}(0)b_0(0,0)^2\\
&-\frac{1}{2}\sum_{1\leq s,t\leq q,q+1\leq j\leq n}\frac{1}{(\abs{\lambda_j}+\abs{\lambda_s})\abs{\lambda_t}\abs{\lambda_s}}\\
&\times\Bigr(\frac{\pr^3\phi}{\pr\ol z_j\pr z_j\pr z_s}(0)\frac{\pr^3\phi}{\pr\ol z_t\pr z_t\pr\ol z_s}(0)+\frac{\pr^3\phi}{\pr\ol z_j\pr z_j\pr\ol z_s}(0)\frac{\pr^3\phi}{\pr\ol z_t\pr z_t\pr z_s}(0)\Bigr)b_0(0,0)^2\\
&+\frac{1}{2}\sum_{1\leq s\leq n,q+1\leq j,t\leq n}\frac{1}{(\abs{\lambda_t}+\abs{\lambda_s})(\abs{\lambda_j}+\abs{\lambda_s})\abs{\lambda_s}}\\
&\times\frac{\pr^3\phi}{\pr\ol z_j\pr z_j\pr z_s}(0)\frac{\pr^3\phi}{\pr\ol z_t\pr z_t\pr\ol z_s}(0)b_0(0,0)^2.
\end{split}
\end{equation} 

From \eqref{s4-e28}, \eqref{s4-e27}, \eqref{s4-e26}, \eqref{s4-e23}, \eqref{s4-e21}, \eqref{s4-e15} and some straight forward computations, we get the main result of this section 

\begin{thm} \label{s4-t4}  
For $c$ in \eqref{s4-e1}, we have
\begin{equation} \label{s4-e29}  
\begin{split}
&c=\frac{1}{4}\Bigr(-\sum_{1\leq j,t\leq q}\frac{1}{\abs{\lambda_j}\abs{\lambda_t}}\frac{\pr^4\phi}{\pr\ol z_j\pr z_j\pr\ol z_t\pr z_t}(0)+\sum_{q+1\leq j,t\leq n}\frac{1}{\abs{\lambda_j}\abs{\lambda_t}}\frac{\pr^4\phi}{\pr\ol z_j\pr z_j\pr\ol z_t\pr z_t}(0)\\
&-2\sum_{q+1\leq j\leq n,1\leq t\leq q}\frac{\abs{\lambda_j}-\abs{\lambda_t}}{\abs{\lambda_j}\abs{\lambda_t}(\abs{\lambda_j}+\abs{\lambda_t})}\frac{\pr^4\phi}{\pr\ol z_j\pr z_j\pr\ol z_t\pr z_t}(0)\Bigr)\\
&-\frac{1}{2}\sum_{q+1\leq j\leq n,1\leq t,s\leq q}\frac{\abs{\lambda_j}-\abs{\lambda_t}}{\abs{\lambda_j}\abs{\lambda_t}(\abs{\lambda_j}+\abs{\lambda_t})(\abs{\lambda_s}+\abs{\lambda_j})}\abs{\frac{\pr^3\phi}{\pr\ol z_s\pr z_j\pr z_t}(0)}^2\\
&+\frac{1}{2}\sum_{q+1\leq j,s\leq n,1\leq t\leq q}\frac{\abs{\lambda_j}-\abs{\lambda_t}}{\abs{\lambda_j}\abs{\lambda_t}(\abs{\lambda_j}+\abs{\lambda_t})(\abs{\lambda_s}+\abs{\lambda_t})}\abs{\frac{\pr^3\phi}{\pr\ol z_s\pr z_j\pr z_t}(0)}^2\\
&+\frac{1}{4}\sum_{q+1\leq j,t\leq n,1\leq s\leq q}\Bigr(\frac{\abs{\lambda_t}^2\abs{\lambda_j}^2}{\abs{\lambda_s}\abs{\lambda_j}\abs{\lambda_t}(\abs{\lambda_t}+\abs{\lambda_j}+\abs{\lambda_s})^2}\bigr(\frac{1}{\abs{\lambda_j}+\abs{\lambda_s}}+\frac{1}{\abs{\lambda_t}+\abs{\lambda_s}}\bigr)^2\\
&+\frac{1}{(\abs{\lambda_t}+\abs{\lambda_j}+\abs{\lambda_s})\abs{\lambda_j}\abs{\lambda_t}}\Bigr)\abs{\frac{\pr^3\phi}{\pr\ol z_s\pr z_j\pr z_t}(0)}^2\\
&+\frac{1}{4}\sum_{q+1\leq s\leq n,1\leq j,t\leq q}\Bigr(\frac{\abs{\lambda_t}^2\abs{\lambda_j}^2}{\abs{\lambda_s}\abs{\lambda_j}\abs{\lambda_t}(\abs{\lambda_t}+\abs{\lambda_j}+\abs{\lambda_s})^2}\bigr(\frac{1}{\abs{\lambda_j}+\abs{\lambda_s}}+\frac{1}{\abs{\lambda_t}+\abs{\lambda_s}}\bigr)^2\\
&+\frac{1}{(\abs{\lambda_t}+\abs{\lambda_j}+\abs{\lambda_s})\abs{\lambda_j}\abs{\lambda_t}}\Bigr)\abs{\frac{\pr^3\phi}{\pr\ol z_s\pr z_j\pr z_t}(0)}^2\\
&-\frac{1}{4}\sum_{1\leq j,t,s\leq q}\frac{1}{\abs{\lambda_j}\abs{\lambda_t}\abs{\lambda_s}}\abs{\frac{\pr^3\phi}
{\pr\ol z_s\pr z_j\pr z_t}(0)}^2\\
&-\frac{1}{4}\sum_{q+1\leq j,t,s\leq n}\frac{1}{\abs{\lambda_j}\abs{\lambda_t}\abs{\lambda_s}}\abs{\frac{\pr^3\phi}
{\pr\ol z_s\pr z_j\pr z_t}(0)}^2\\
&+\frac{1}{2}\sum_{q+1\leq j,t\leq n,1\leq s\leq q}\frac{\abs{\lambda_s}}{\abs{\lambda_j}\abs{\lambda_t}
(\abs{\lambda_j}+\abs{\lambda_s})(\abs{\lambda_t}+\abs{\lambda_s})}
\frac{\pr^3\phi}{\pr\ol z_j\pr z_j\pr z_s}(0)\frac{\pr^3\phi}{\pr\ol z_t\pr z_t\pr\ol z_s}(0)\\ 
&+\frac{1}{2}\sum_{1\leq j,t\leq q,q+1\leq s\leq n}\frac{\abs{\lambda_s}}{\abs{\lambda_j}\abs{\lambda_t}
(\abs{\lambda_j}+\abs{\lambda_s})(\abs{\lambda_t}+\abs{\lambda_s})}
\frac{\pr^3\phi}{\pr\ol z_j\pr z_j\pr z_s}(0)\frac{\pr^3\phi}{\pr\ol z_t\pr z_t\pr\ol z_s}(0)\\
&-\frac{1}{2}\sum_{1\leq s,t\leq q,q+1\leq j\leq n}\frac{1}{(\abs{\lambda_j}+\abs{\lambda_t})
(\abs{\lambda_j}+\abs{\lambda_s})\abs{\lambda_t}}\\
&\times\Bigr(\frac{\pr^3\phi}{\pr\ol z_j\pr z_j\pr z_s}(0)\frac{\pr^3\phi}{\pr\ol z_t\pr z_t\pr\ol z_s}(0)+\frac{\pr^3\phi}{\pr\ol z_j\pr z_j\pr z_s}(0)\frac{\pr^3\phi}{\pr\ol z_t\pr z_t\pr\ol z_s}(0)\Bigr)\\
&-\frac{1}{2}\sum_{1\leq t\leq q,q+1\leq j,s\leq n}\frac{1}{(\abs{\lambda_j}+\abs{\lambda_t})
(\abs{\lambda_t}+\abs{\lambda_s})\abs{\lambda_j}}\\
&\times\Bigr(\frac{\pr^3\phi}{\pr\ol z_j\pr z_j\pr z_s}(0)\frac{\pr^3\phi}{\pr\ol z_t\pr z_t\pr\ol z_s}(0)+\frac{\pr^3\phi}{\pr\ol z_j\pr z_j\pr z_s}(0)\frac{\pr^3\phi}{\pr\ol z_t\pr z_t\pr\ol z_s}(0)\Bigr).
\end{split}
\end{equation}
\end{thm} 

\section{The end of the proof of Theorem~\ref{t-main}}

In this section, we will use the same notations as in section $1$. Let 
\[F:\Lambda^{1,0}T(\Complex^n)\To\Lambda^{1,0}T(\Complex^n)\otimes\Lambda^{p_1,q_1}T^*(\Complex^n)\] 
and 
\[T: \Lambda^{1,0}T(\Complex^n))\To\Lambda^{1,0}T(\Complex^n)\otimes\Lambda^{p_2,q_2}T^*(\Complex^n)\]
be linear operaors, where $p_1, q_1, p_2, q_2\in\mathbb Z$, $p_1, q_1, p_2, q_2\geq0$. We write $F=\left(F_{j,k}\right)^n_{j,k=1}$, $F_{j,k}\in\Lambda^{p_1,q_1}T^*(\Complex^n)$, $1\leq j,k\leq n$, and $T=\left(T_{j,k}\right)^n_{j,k=1}$, $T_{j,k}\in\Lambda^{p_2,q_2}T^*(\Complex^n)$, $j, k=1,\ldots,n$, as in 
\eqref{s1-e12}.
$TF:\Lambda^{1,0}T(\Complex^n)\To\Lambda^{1,0}T(\Complex^n)\otimes\Lambda^{p_1+p_2,q_1+q_2}T^*(\Complex^n)$
is the linear operator defined by $TF\frac{\pr}{\pr z_k}=\sum^n_{s,j=1}\frac{\pr}{\pr z_j}\otimes(T_{j,s}\wedge F_{s,k})$, $k=1,\ldots,n$. We have $TF=\left(a_{j,k}\right)^n_{j,k=1}$. We can compute 
\begin{equation} \label{s5-e1}
a_{j,k}=\sum^n_{s=1}T_{j,s}\wedge F_{s,k}\in\Lambda^{p_1+p_2,q_1+q_2}T^*(\Complex^n),\ \ j,k=1,\ldots,n.
\end{equation} 

For $p\in\Complex^n$, we may assume that $p=0$ and $\phi(z)=\sum^n_{j=1}\lambda_j\abs{z_j}^2+O(\abs{z}^3)$ near $0$ and as before we suppose that $\lambda_j<0$, $j=1,\ldots,q$, and $\lambda_j>0$, $j=q+1,\ldots,n$. Let $M_\phi:C^\infty(\Complex^n;\, \Lambda^{1,0}T(\Complex^n))\To C^\infty(\Complex^n;\, \Lambda^{1,0}T(\Complex^n))$ be as in \eqref{s1-e15}. We recall that 
$M_\phi=\left(\frac{\pr^2\phi}{\pr\ol z_j\pr z_k}\right)^n_{j,k=1}$. We have 
\begin{equation} \label{s5-e1-1}
U_j=\frac{\pr}{\pr z_j}\ \ \mbox{at $0$},\ \ j=1,\ldots,n,
\end{equation} 
where $U_j$, $j=1,\ldots,n$, are given by \eqref{s1-e17}. Moreover, we have 
\begin{align}
&(\frac{\pr}{\pr z_j}\ |\ \frac{\pr}{\pr z_k})_{\abs{\phi}}=(\frac{\pr}{\pr\ol z_j}\ |\ \frac{\pr}{\pr\ol z_k})_{\abs{\phi}}=\delta_{j,k}\abs{\lambda_j}\ \ \mbox{at $0$},\ \ j,k=1,\ldots,n,\label{s5-e1-2}\\
&(d z_j\ |\ d z_k)_{\abs{\phi}}=(d\ol z_j\ |\ d\ol z_k)_{\abs{\phi}}=\delta_{j,k}\frac{1}{\abs{\lambda_j}}\ \ \mbox{at $0$},\ \ j,k=1,\ldots,n.\label{s5-e1-3}
\end{align}
Here $\delta_{j,k}=1$ if $j=k$, $\delta_{j,k}=0$ if $j\neq k$. For the definition of the Hermitian metrix $(\ |\ )_{\abs{\phi}}$, 
see the discussion after \eqref{s1-e17}.

Put 
$M^{-1}_\phi(z)=\left(b_{j,k}(z)\right)^n_{j,k=1}$. We claim that 
\begin{equation} \label{s5-e2} 
b_{j,k}(z)=\frac{\delta_{j,k}}{\lambda_j}-\frac{1}{\lambda_j\lambda_k}\sum^n_{s=1}\Bigr(\frac{\pr^3\phi}{\pr\ol z_j\pr z_k\pr z_s}(0)z_s+\frac{\pr^3\phi}{\pr\ol z_j\pr z_k\pr\ol z_s}(0)\ol z_s\Bigr)+O(\abs{z}^2),
\end{equation} 
$j, k=1,\ldots,n$, near $0\in\Complex^n$.  In fact, if we put $\Td B=\left(\Td b_{j,k}\right)^n_{j,k=1}$, $\Td b_{j,k}=\frac{\delta_{j,k}}{\lambda_j}-\frac{1}{\lambda_j\lambda_k}\sum^n_{s=1}\Bigr(\frac{\pr^3\phi}{\pr\ol z_j\pr z_k\pr z_s}(0)z_s+\frac{\pr^3\phi}{\pr\ol z_j\pr z_k\pr\ol z_s}(0)\ol z_s\Bigr)$. Then $M_\phi\Td B=C$, $C=\left(c_{j,k}\right)^n_{j,k=1}$, 
\begin{align} \label{s5-e3}
c_{j,k}&=\sum^n_{t=1}\frac{\pr^2\phi}{\pr\ol z_j\pr z_t}(z)\Td b_{t,k}(z)\nonumber\\
&=\sum^n_{t=1}\frac{\pr^2\phi}{\pr\ol z_j\pr z_t}(z)\Bigr(\frac{\delta_{t,k}}{\lambda_t}-\frac{1}{\lambda_t\lambda_k}\sum^n_{s=1}\bigr(\frac{\pr^3\phi}{\pr\ol z_t\pr z_k\pr z_s}(0)z_s+\frac{\pr^3\phi}{\pr\ol z_t\pr z_k\pr\ol z_s}(0)\ol z_s\bigr)\Bigr)\nonumber\\
&=\sum^n_{t=1}\Bigr(\lambda_t\delta_{j,t}+\sum^n_{s=1}\bigr(\frac{\pr^3\phi}{\pr\ol z_j\pr z_t\pr z_s}(0)z_s+\frac{\pr^3\phi}{\pr\ol z_j\pr z_t\pr\ol z_s}(0)\ol z_s\bigr)\Bigr)\nonumber \\
&\times\Bigr(\frac{\delta_{t,k}}{\lambda_t}-\frac{1}{\lambda_t\lambda_k}\sum^n_{s=1}\bigr(\frac{\pr^3\phi}{\pr\ol z_t\pr z_k\pr z_s}(0)z_s+\frac{\pr^3\phi}{\pr\ol z_t\pr z_k\pr\ol z_s}(0)\ol z_s\bigr)\Bigr)+O(\abs{z}^2)\nonumber\\
&=\delta_{j,k}+O(\abs{z}^2).
\end{align} 
Thus, $\Td B=M_\phi^{-1}(z)+O(\abs{z}^2)$. \eqref{s5-e2} follows. 

Let $Q:\Lambda^{1,0}T(\Complex^n)\To\Lambda^{1,0}T^*(\Complex^n)\otimes\Lambda^{1,0}T(\Complex^n)$ be as in \eqref{s1-e21}. 
We write $Q=\left(Q_{j,k}\right)^n_{j,k=1}$, $Q_{j,k}(z)=\sum^n_{s=1}\frac{\pr^3\phi}{\pr\ol z_j\pr z_k\pr z_s}(0)q_{j,k,s}dz_s$ 
at $0$, where 
\begin{equation} \label{s5-e4} 
\begin{split}
&q_{j,k,s}=\frac{\abs{\lambda_k}\delta_j(k)+\abs{\lambda_s}\delta_j(s)}{\abs{\lambda_j}+\abs{\lambda_k}\delta_j(k)+\abs{\lambda_s}\delta_j(s)}
-\delta_j(k)\delta_j(s)\times\\
&\frac{\abs{\lambda_k}^2\abs{\lambda_s}^2}{(\abs{\lambda_j}+\abs{\lambda_k}+\abs{\lambda_s})^2}\bigr(\frac{1}{\abs{\lambda_j}+\abs{\lambda_k}}+\frac{1}{\abs{\lambda_j}+\abs{\lambda_s}}\bigr)^2.
\end{split}
\end{equation} 
We recall that the definition of $\delta_k(j)$ is given by \eqref{s1-e20}. Put $\ddbar M_\phi^{-1}=\left(d_{j,k}\right)^n_{j,k=1}$. 
From \eqref{s5-e2}, we see that 
\begin{equation} \label{s5-e5} 
d_{j,k}=-\sum^n_{s=1}\frac{1}{\lambda_j\lambda_k}\frac{\pr^3\phi}{\pr\ol z_j\pr z_k\pr\ol z_s}(0)d\ol z_s
\end{equation} 
at $0$. Put $\ddbar M^{-1}_\phi Q=\left(f_{j,k}\right)^n_{j,k=1}$. From \eqref{s5-e1}, \eqref{s5-e5}, we see that 
\begin{equation} \label{s5-e6}
f_{k,j}=-\sum_{1\leq t,s,u\leq n}\frac{1}{\lambda_k\lambda_t}q_{t,j,u}\frac{\pr^3\phi}{\pr\ol z_k\pr z_t\pr\ol z_s}(0)\frac{\pr^3\phi}{\pr\ol z_t\pr z_j\pr z_u}(0)d\ol z_s\wedge dz_u
\end{equation} 
at $0$.
As in section $1$, put $e_j=\frac{1}{\sqrt{\abs{\lambda_j}}}U_j$, $j=1,\ldots,n$. From \eqref{s5-e6}, it is not difficult to see that
\begin{equation} \label{s5-e7} 
 < (\ddbar M^{-1}_\phi Qe_j\ |\ e_k)_{\abs{\phi}}, \ol e_j\wedge e_k>(0)=-\sum_{1\leq t\leq n}\frac{1}
 {\lambda_k\lambda_t\abs{\lambda_j}}q_{t,j,k}\abs{\frac{\pr^3\phi}{\pr\ol z_k\pr z_t\pr\ol z_j}(0)}^2,
\end{equation} 
$j,k=1,\ldots,n$. From \eqref{s5-e7} and \eqref{s5-e4}, it is straight forward to see that

\begin{prop} \label{s5-l1} 
We have that 
\begin{equation} \label{s5-e8} 
\begin{split}
&< (\ddbar M^{-1}_\phi Qe_j\ |\ e_k)_{\abs{\phi}}, \ol e_j\wedge e_k>(0)=\sum_{1\leq t\leq q}\frac{1}{\abs{\lambda_t}
\abs{\lambda_j}\abs{\lambda_k}}\abs{\frac{\pr^3\phi}{\pr\ol z_k\pr z_t\pr\ol z_j}(0)}^2\\
&\times\Bigr(\frac{\abs{\lambda_j}+\abs{\lambda_k}}{\abs{\lambda_j}+
\abs{\lambda_k}+\abs{\lambda_t}}-\frac{\abs{\lambda_j}^2\abs{\lambda_k}^2}{(\abs{\lambda_j}+\abs{\lambda_k}+\abs{\lambda_t})^2}\bigr(\frac{1}{\abs{\lambda_j}+\abs{\lambda_t}}+\frac{1}{\abs{\lambda_k}+\abs{\lambda_t}}\bigr)^2\Bigr),
\end{split} 
\end{equation} 
where $q+1\leq j,k\leq n$, 
\begin{equation} \label{s5-e9} 
\begin{split}
&< (\ddbar M^{-1}_\phi Qe_j\ |\ e_k)_{\abs{\phi}}, \ol e_j\wedge e_k>(0)=\sum_{q+1\leq t\leq n}\frac{1}{\abs{\lambda_t}
\abs{\lambda_j}\abs{\lambda_k}}\abs{\frac{\pr^3\phi}{\pr\ol z_k\pr z_t\pr\ol z_j}(0)}^2\\
&\times\Bigr(\frac{\abs{\lambda_j}+\abs{\lambda_k}}{\abs{\lambda_j}+
\abs{\lambda_k}+\abs{\lambda_t}}-\frac{\abs{\lambda_j}^2\abs{\lambda_k}^2}{(\abs{\lambda_j}+\abs{\lambda_k}+\abs{\lambda_t})^2}\bigr(\frac{1}{\abs{\lambda_j}+\abs{\lambda_t}}+\frac{1}{\abs{\lambda_k}+\abs{\lambda_t}}\bigr)^2\Bigr),
\end{split}
\end{equation} 
where $1\leq j,k\leq q$, 
\begin{equation} \label{s5-e10} 
\begin{split}
&< (\ddbar M^{-1}_\phi Qe_j\ |\ e_k)_{\abs{\phi}}, \ol e_j\wedge e_k>(0)\\
&=-\sum_{1\leq t\leq q}\frac{1}{\abs{\lambda_t}
\abs{\lambda_j}\abs{\lambda_k}}\frac{\abs{\lambda_j}}{\abs{\lambda_t}+\abs{\lambda_j}}\abs{\frac{\pr^3\phi}{\pr\ol z_k\pr z_t\pr\ol z_j}(0)}^2\\
&+\sum_{q+1\leq t\leq n}\frac{1}{\abs{\lambda_t}
\abs{\lambda_j}\abs{\lambda_k}}\frac{\abs{\lambda_k}}{\abs{\lambda_t}+\abs{\lambda_k}}\abs{\frac{\pr^3\phi}{\pr\ol z_k\pr z_t\pr\ol z_j}(0)}^2,
\end{split}
\end{equation} 
where $q+1\leq j\leq n$, $1\leq k\leq q$, and 
\begin{equation} \label{s5-e11} 
\begin{split}
&< (\ddbar M^{-1}_\phi Qe_j\ |\ e_k)_{\abs{\phi}}, \ol e_j\wedge e_k>(0)\\
&=\sum_{1\leq t\leq q}\frac{1}{\abs{\lambda_t}
\abs{\lambda_j}\abs{\lambda_k}}\frac{\abs{\lambda_k}}{\abs{\lambda_t}+\abs{\lambda_k}}\abs{\frac{\pr^3\phi}{\pr\ol z_k\pr z_t\pr\ol z_j}(0)}^2\\
&-\sum_{q+1\leq t\leq n}\frac{1}{\abs{\lambda_t}
\abs{\lambda_j}\abs{\lambda_k}}\frac{\abs{\lambda_j}}{\abs{\lambda_t}+\abs{\lambda_j}}\abs{\frac{\pr^3\phi}{\pr\ol z_k\pr z_t\pr\ol z_j}(0)}^2,
\end{split}
\end{equation} 
where $1\leq j\leq q$, $q+1\leq k\leq n$.
\end{prop}

Let $\Theta_\phi:C^\infty(\Complex^n;\, \Lambda^{1,0}T(\Complex^n))\To C^\infty(\Complex^n;\, \Lambda^{1,1}T^*(\Complex^n)\otimes\Lambda^{1,0}T(\Complex^n))$ be as in \eqref{s1-e19}. We recall that $\Theta_\phi=\left(\ddbar\theta_{j,k}\right)^n_{j,k=1}=\left(\Theta_{j,k}\right)^n_{j,k=1}$, where $\theta=h^{-1}\pr h=\left(\theta_{j,k}\right)^n_{j,k=1}$, $h=\left(\frac{\pr^2\phi}{\pr\ol z_j\pr z_k}\right)^n_{j,k=1}$. It is not difficult to see that 
\begin{equation} \label{s5-e12}
\Theta_{j,k}=\frac{1}{\lambda_j}\ddbar\pr(\frac{\pr^2\phi}{\pr\ol z_j\pr z_k})-\sum^n_{t=1}\frac{1}{\lambda_t\lambda_j}\ddbar(\frac{\pr^2\phi}{\pr\ol z_j\pr z_t})\wedge\pr(\frac{\pr^2\phi}{\pr\ol z_t\pr z_k}),
\end{equation}
$j,k=1,\ldots,n$. From \eqref{s5-e12}, it is straight forward to see that 
\begin{prop} \label{s5-p2} 
We have that 
\begin{equation} \label{s5-e13}
\begin{split}
&< (\Theta_\phi e_j\ |\ e_k)_{\abs{\phi}}, \ol e_j\wedge e_k>(0)=
\frac{1}{\abs{\lambda_j}\abs{\lambda_k}}\frac{\pr^4\phi}{\pr\ol z_j\pr z_j\pr\ol z_k\pr z_k}(0)\\
&+\sum^q_{t=1}\frac{1}{\abs{\lambda_t}\abs{\lambda_j}\abs{\lambda_k}}\abs{\frac{\pr^3\phi}{\pr\ol z_j\pr z_t\pr\ol z_k}(0)}^2-\sum^n_{t=q+1}\frac{1}{\abs{\lambda_t}\abs{\lambda_j}\abs{\lambda_k}}\abs{\frac{\pr^3\phi}{\pr\ol z_j\pr z_t\pr\ol z_k}(0)}^2,
\end{split}
\end{equation} 
where $q+1\leq j,k\leq n$, 
\begin{equation} \label{s5-e14} 
\begin{split}
&< (\Theta_\phi e_j\ |\ e_k)_{\abs{\phi}}, \ol e_j\wedge e_k>(0)=
-\frac{1}{\abs{\lambda_j}\abs{\lambda_k}}\frac{\pr^4\phi}{\pr\ol z_j\pr z_j\pr\ol z_k\pr z_k}(0)\\
&-\sum^q_{t=1}\frac{1}{\abs{\lambda_t}\abs{\lambda_j}\abs{\lambda_k}}\abs{\frac{\pr^3\phi}{\pr\ol z_j\pr z_t\pr\ol z_k}(0)}^2+\sum^n_{t=q+1}\frac{1}{\abs{\lambda_t}\abs{\lambda_j}\abs{\lambda_k}}\abs{\frac{\pr^3\phi}{\pr\ol z_j\pr z_t\pr\ol z_k}(0)}^2,
\end{split}
\end{equation} 
where $1\leq j,k\leq q$, 
\begin{equation} \label{s5-e15} 
\begin{split}
&< (\Theta_\phi e_j\ |\ e_k)_{\abs{\phi}}, \ol e_j\wedge e_k>(0)=
-\frac{1}{\abs{\lambda_j}\abs{\lambda_k}}\frac{\pr^4\phi}{\pr\ol z_j\pr z_j\pr\ol z_k\pr z_k}(0)\\
&-\sum^q_{t=1}\frac{1}{\abs{\lambda_t}\abs{\lambda_j}\abs{\lambda_k}}\abs{\frac{\pr^3\phi}{\pr\ol z_j\pr z_t\pr\ol z_k}(0)}^2+\sum^n_{t=q+1}\frac{1}{\abs{\lambda_t}\abs{\lambda_j}\abs{\lambda_k}}\abs{\frac{\pr^3\phi}{\pr\ol z_j\pr z_t\pr\ol z_k}(0)}^2,
\end{split}
\end{equation} 
where $q+1\leq j\leq n$, $1\leq k\leq q$, and 
\begin{equation} \label{s5-e16} 
\begin{split}
&< (\Theta_\phi e_j\ |\ e_k)_{\abs{\phi}}, \ol e_j\wedge e_k>(0)=
\frac{1}{\abs{\lambda_j}\abs{\lambda_k}}\frac{\pr^4\phi}{\pr\ol z_j\pr z_j\pr\ol z_k\pr z_k}(0)\\
&+\sum^q_{t=1}\frac{1}{\abs{\lambda_t}\abs{\lambda_j}\abs{\lambda_k}}\abs{\frac{\pr^3\phi}{\pr\ol z_j\pr z_t\pr\ol z_k}(0)}^2-\sum^n_{t=q+1}\frac{1}{\abs{\lambda_t}\abs{\lambda_j}\abs{\lambda_k}}\abs{\frac{\pr^3\phi}{\pr\ol z_j\pr z_t\pr\ol z_k}(0)}^2,
\end{split}
\end{equation} 
where $1\leq j\leq q$, $q+1\leq k\leq n$.
\end{prop}

As in section $1$, define  
\[R=\Theta_\phi-(\ddbar M_\phi^{-1})Q:\Lambda^{1,0}T(\Complex^n)\To \Lambda^{1,1}T^*(\Complex^n)\otimes\Lambda^{1,0}T(\Complex^n).\]
From \eqref{s5-e8}--\eqref{s5-e11} and \eqref{s5-e13}--\eqref{s5-e14}, it is not difficult to see that 
\begin{equation} \label{s5-e17}  
\begin{split} 
&\frac{1}{4}\sum^n_{j,k=1}(1+\delta_j(k)\frac{\abs{\lambda_j}-\abs{\lambda_k}}{\abs{\lambda_j}+\abs{\lambda_k}}<(Re_j\ |\ e_k)_{\abs{\phi}}, \ol e_j\wedge e_k>(0)\\
&=\frac{1}{4}\Bigr(-\sum_{1\leq j,k\leq q}\frac{1}{\abs{\lambda_j}\abs{\lambda_t}}\frac{\pr^4\phi}{\pr\ol z_j\pr z_j\pr\ol z_k\pr z_k}(0)+\sum_{q+1\leq j,k\leq n}\frac{1}{\abs{\lambda_j}\abs{\lambda_k}}\frac{\pr^4\phi}{\pr\ol z_j\pr z_j\pr\ol z_k\pr z_k}(0)\\
&-2\sum_{q+1\leq j\leq n,1\leq k\leq q}\frac{\abs{\lambda_j}-\abs{\lambda_k}}{\abs{\lambda_j}\abs{\lambda_k}(\abs{\lambda_j}+\abs{\lambda_k})}\frac{\pr^4\phi}{\pr\ol z_j\pr z_j\pr\ol z_k\pr z_k}(0)\Bigr)\\
&-\frac{1}{2}\sum_{q+1\leq j\leq n,1\leq k,s\leq q}\frac{\abs{\lambda_j}-\abs{\lambda_k}}{\abs{\lambda_j}\abs{\lambda_k}(\abs{\lambda_j}+\abs{\lambda_k})(\abs{\lambda_s}+\abs{\lambda_j})}\abs{\frac{\pr^3\phi}{\pr\ol z_s\pr z_j\pr z_k}(0)}^2\\
&+\frac{1}{2}\sum_{q+1\leq j,s\leq n,1\leq k\leq q}\frac{\abs{\lambda_j}-\abs{\lambda_k}}{\abs{\lambda_j}\abs{\lambda_k}(\abs{\lambda_j}+\abs{\lambda_k})(\abs{\lambda_s}+\abs{\lambda_k})}\abs{\frac{\pr^3\phi}{\pr\ol z_s\pr z_j\pr z_k}(0)}^2\\
&+\frac{1}{4}\sum_{q+1\leq j,k\leq n,1\leq s\leq q}\Bigr(\frac{\abs{\lambda_k}^2\abs{\lambda_j}^2}{\abs{\lambda_s}\abs{\lambda_j}\abs{\lambda_k}(\abs{\lambda_k}+\abs{\lambda_j}+\abs{\lambda_s})^2}\bigr(\frac{1}{\abs{\lambda_j}+\abs{\lambda_s}}+\frac{1}{\abs{\lambda_k}+\abs{\lambda_s}}\bigr)^2\\
&+\frac{1}{(\abs{\lambda_k}+\abs{\lambda_j}+\abs{\lambda_s})\abs{\lambda_j}\abs{\lambda_k}}\Bigr)\abs{\frac{\pr^3\phi}{\pr\ol z_s\pr z_j\pr z_k}(0)}^2\\
&+\frac{1}{4}\sum_{q+1\leq s\leq n,1\leq j,k\leq q}\Bigr(\frac{\abs{\lambda_k}^2\abs{\lambda_j}^2}{\abs{\lambda_s}\abs{\lambda_j}\abs{\lambda_k}(\abs{\lambda_k}+\abs{\lambda_j}+\abs{\lambda_s})^2}\bigr(\frac{1}{\abs{\lambda_j}+\abs{\lambda_s}}+\frac{1}{\abs{\lambda_k}+\abs{\lambda_s}}\bigr)^2\\
&+\frac{1}{(\abs{\lambda_k}+\abs{\lambda_j}+\abs{\lambda_s})\abs{\lambda_j}\abs{\lambda_k}}\Bigr)\abs{\frac{\pr^3\phi}{\pr\ol z_s\pr z_j\pr z_k}(0)}^2\\
&-\frac{1}{4}\sum_{1\leq j,k,s\leq q}\frac{1}{\abs{\lambda_j}\abs{\lambda_k}\abs{\lambda_s}}\abs{\frac{\pr^3\phi}
{\pr\ol z_s\pr z_j\pr z_k}(0)}^2\\
&-\frac{1}{4}\sum_{q+1\leq j,k,s\leq n}\frac{1}{\abs{\lambda_j}\abs{\lambda_k}\abs{\lambda_s}}\abs{\frac{\pr^3\phi}
{\pr\ol z_s\pr z_j\pr z_k}(0)}^2.
\end{split}
\end{equation} 

From \eqref{s5-e4}, we can check that 
\[(Qe_j\ |\ e_j)=\sum^n_{s=1}
\delta_j(s)\frac{\abs{\lambda_s}}{\abs{\lambda_j}(\abs{\lambda_j}+\abs{\lambda_s})}\frac{\pr^3\phi}{\pr\ol z_j\pr z_j\pr z_s}(0)dz_s\] 
at $0$. From this, it is straight forward to see that 
\begin{equation} \label{s5-e18} 
\begin{split} 
&-\sum^n_{j,k=1}\delta_j(k)\frac{\abs{\lambda_j}}{\abs{\lambda_j}+\abs{\lambda_k}}{\rm Re\,}
((Qe_j\ |\ e_j)\ |\ \pr M_\phi e_k\ |\ e_k)_{\abs{\phi}}(0)\\
&=-\frac{1}{2}\sum_{1\leq s,k\leq q,q+1\leq j\leq n}\frac{1}{(\abs{\lambda_j}+\abs{\lambda_k})
(\abs{\lambda_j}+\abs{\lambda_s})\abs{\lambda_k}}\\
&\times\Bigr(\frac{\pr^3\phi}{\pr\ol z_j\pr z_j\pr z_s}(0)\frac{\pr^3\phi}{\pr\ol z_k\pr z_k\pr\ol z_s}(0)+\frac{\pr^3\phi}{\pr\ol z_j\pr z_j\pr z_s}(0)\frac{\pr^3\phi}{\pr\ol z_k\pr z_k\pr\ol z_s}(0)\Bigr)\\
&-\frac{1}{2}\sum_{1\leq k\leq q,q+1\leq j,s\leq n}\frac{1}{(\abs{\lambda_j}+\abs{\lambda_k})
(\abs{\lambda_k}+\abs{\lambda_s})\abs{\lambda_j}}\\
&\times\Bigr(\frac{\pr^3\phi}{\pr\ol z_j\pr z_j\pr z_s}(0)\frac{\pr^3\phi}{\pr\ol z_k\pr z_k\pr\ol z_s}(0)+\frac{\pr^3\phi}{\pr\ol z_j\pr z_j\pr z_s}(0)\frac{\pr^3\phi}{\pr\ol z_k\pr z_k\pr\ol z_s}(0)\Bigr)
\end{split}
\end{equation} 
and 
\begin{equation} \label{s5-e19} 
\begin{split} 
&\frac{1}{2}\abs{\sum^n_{j=1}(Qe_j\ |\ e_j)}^2_{\abs{\phi}} \\
&=\frac{1}{2}\sum_{q+1\leq j,k\leq n,1\leq s\leq q}\frac{\abs{\lambda_s}}{\abs{\lambda_j}\abs{\lambda_k}
(\abs{\lambda_j}+\abs{\lambda_s})(\abs{\lambda_k}+\abs{\lambda_s})}
\frac{\pr^3\phi}{\pr\ol z_j\pr z_j\pr z_s}(0)\frac{\pr^3\phi}{\pr\ol z_k\pr z_k\pr\ol z_s}(0)\\ 
&+\frac{1}{2}\sum_{1\leq j,k\leq q,q+1\leq s\leq n}\frac{\abs{\lambda_s}}{\abs{\lambda_j}\abs{\lambda_k}
(\abs{\lambda_j}+\abs{\lambda_s})(\abs{\lambda_k}+\abs{\lambda_s})}
\frac{\pr^3\phi}{\pr\ol z_j\pr z_j\pr z_s}(0)\frac{\pr^3\phi}{\pr\ol z_k\pr z_k\pr\ol z_s}(0).
\end{split}
\end{equation} 
Combining \eqref{s5-e19}, \eqref{s5-e18} and \eqref{s5-e17} with \eqref{s4-e29}, Theorem~\ref{t-main} follows.


\begin{bibdiv}
\begin{biblist} 
\bib{BS05}{article}{
	title={Asymptotics for Bergman-Hodge kernels for high powers of complex line bundles},
	author={Berman, R.},
	author={Sj\"ostrand, J.},
	eprint={arXiv.org/abs/math.CV/0511158}
}  

\bib{Cat97}{article}{
	title={The Bergman kernel and a theorem of Tian},
	author={Catlin, D.},
	conference={
		title={Analysis and geometry in several complex variables},
		address={Katata},
		date={1997}
	},
	book={
		series={Trends in Math.},
		publisher={Birkhauser},
		address={Boston, MA},
		date={1999}
	},
	pages={1--23}
}  
\bib{Do:01}{article}{ 
   title={Scalar curvature and projective embeddings. {I}},
   author={Donaldson, S. K.}, 
   journal={J. Differential Geom.},
   volume={59},
   date={2001},
   pages={479--522} 
} 

\bib{Fine08}{article}{ 
   title={Calabi flow and projective embeddings},
   author={Fine, J.}, 
   journal={J. Differential Geom.},
   volume={84},
   date={2010},
   pages={489--523} 
} 

\bib{Fine10}{article}{
	title={Quantisation and the Hessian of Mabuchi energy},
	author={Fine, J.},
	eprint={1009.4543}
}  

\bib{Hor03}{book}{ 
    title={The analysis of linear partial differential operators I distribution
        theory and Fourier analysis},
    author={H\"{o}rmander, L.}, 
    series={Classics in Mathematics}, 
    publisher={Springer-Verlag}
    address={Berlin},
    date={2003} 
}

\bib{Hsiao08}{article}{
	title={Projections in several complex variables},
	author={Hsiao, C-Y.},
	journal={M\'em. Soc. Math. France, Nouv. S\'er.},
   number={123},
   date={2010},
   pages={131 p.} 
} 

\bib{Hsiao12}{article}{
	title={On the coefficients of the asymptotic expansion of the kernel of {Berezin}-{Toeplitz} quantization},
	author={Hsiao, C-Y.},
	journal={Annals of Global Analysis and Geometry},
    volume={42},
   date={2012},
   pages={207--245} 
} 

\bib{HM11}{article}{
	title={Asymptotics of spectral function of lower energy forms and Bergman kernel of semi-positive and big line bundles},
	author={Hsiao, C-Y.},
	author={Marinescu, G.}, 
	eprint={arXiv:1112.5464}
}  

\bib{Lu12}{article}{ 
   title={Thesis in preparation (under direction of {Ma} and {Marinescu})},
   author={Lu, W.}, 
} 

\bib{Lu00}{article}{ 
   title={On the lower order terms of the asymptotic expansion of
      {T}ian-{Y}au-{Z}elditch},
   author={Lu, Z.}, 
   journal={Amer. J. Math.},
   volume={122},
   number={2},
   date={2000},
   pages={235--273} 
} 

\bib{MM06}{article}{ 
   title={The first coefficients of the asymptotic expansion of the Bergman kernel
of the $spin^c$ Dirac operator},
   author={Ma, X.}, 
   author={Marinescu, G.}, 
   journal={Internat. J. Math},
   number={17},
   date={2006},
   pages={737-759} 
} 

\bib{MM07}{book}{ 
    title={Holomorphic {Morse} inequalities and {Bergman} kernels},
    author={Ma, X.}, 
    author={Marinescu, G.},
    series={Progress in Math.}, 
    publisher={Birkh{\"a}user}
    address={Basel},
    date={2007} 
}


\bib{MM08a}{article}{ 
   title={Generalized {Bergman} kernels on symplectic manifolds},
   author={Ma, X.}, 
   author={Marinescu, G.}, 
   journal={Adv. Math.},
   volume={217}
   number={4},
   date={2008},
   pages={1756--1815} 
} 

\bib{MM10}{article}{ 
   title={Berezin-{Toeplitz} quantization on {K\"ahler} manifolds},
   author={Ma, X.}, 
   author={Marinescu, G.}, 
   journal={J. reine angew. Math.},
   volume={662}
   date={2012},
   pages={1--56} 
} 

\bib{Tian}{article}{
	title={On a set of polarized {K\"ahler} metrics on algebraic
  manifolds},
	author={Tian, G.},
	journal={J. Differential Geom.},
	volume={32},
	date={1990},
	pages={99--130}
}  

\bib{Zel98}{article}{
	title={Szeg\"o kernels and a theorem of Tian},
	author={Zelditch, S.},
	journal={Internat. Math. Res. Notices},
	number={6},
	date={1998},
	pages={317--331}
}  
\end{biblist} 
\end{bibdiv}
\end{document}